\documentclass[final,12pt]{amsart} 
\usepackage{amsmath, bm} 
\usepackage{amssymb} 
\usepackage{graphicx} 
\numberwithin{equation}{section}

\def\e{\varepsilon}
\def\br{\breve}
 
\title[  INVARIANT TORI OF PERIODIC SYSTEMS  ] 
{ WEAKENED XVI HILBERT'S PROBLEM:\\ INVARIANT TORI OF THE PERIODIC SYSTEMS WITH NINE EQUILIBRIUM POINTS IN HAMILTONIAN UNPERTURBED PART }

\author[ V. V. Basov, A. S. Zhukov ]
{  Vladimir V. Basov$^{1*}$, Artem S. Zhukov$^{1}$}

\address{$^1$Faculty of Mathematics and Mechanics, St. Petersburg State University, 
	Universitetsky prospekt, 28, 198504, Peterhof, St. Petersburg, Russia}

\email{vlvlbasov@rambler.ru,  artzhukov1111@gmail.ru}



\begin{document}
	
\bigskip
\begin{abstract}	
			Two classes of two-dimensional time-periodic systems of ordinary differential equations with a small parameter $\e>0$ in the perturbed part, which is continuous and, for $\e=0$, analytic in zero, are studied. 
			Depending on the presence or absence of the common factor $\e$, these classes contain the system with "fast" or "slow" time.
			The unperturbed part of these systems is generated by the hamiltonian $H=2x^2-x^4+\gamma(2y^2-y^4)\ (\gamma\in(0,1]).$ 
			
			The universal method, meaning it can be applied to any hamiltonian, called the method of the generating tori splitting (GTS method), is developed and applied to the research of such systems. 
			For arbitrary system of any class, this method allows to find the sets of the initial values for the solutions of the corresponding unperturbed system, and for each such set, to provide the explicit conditions on the system perturbations independent of the parameter.
			Any such set, that satisfies the aforementioned conditions, determines the solution of the unperturbed system, which parametrizes the generating cycle. The obtained cycle is a generatrix of an invariant cylindrical surface. It is proven that the system has two-periodic invariant surface, homeomorphic to torus, if time is factored with the respect to the period, in the small with the respect to $\e$ neighbourhood of this surface. The formula and the asymptotic extension are provided for this surface, the number of properties is discovered. 
			
			An  example of the set of systems with eleven invariant tori and a perturbation, 
			which average value is a three-term polynomial of the third degree, is constructed as a demonstration of the practical use of the GTS method.
			
			The GTS method is an universal alternative to the so-called method
			of detection functions and the  Melnikov function method, which  are used in studies  concerning  the weakened XVI Hilbert's problem
			on the evaluation of a number of limit cycles of autonomous systems with the hamiltonian unperturbed part.
			The GTS method allows not only to evaluate the lower bound of the analogue of the Hilbert's cyclicity value, which determines the amount of the invariant tori in the periodic systems with "slow" time, but also to solve the same problem for the periodic systems of any even degree with the common factor $\e$ in its right-hand side. The results can be put into practice, while researching the systems of the ordinary differential equations of the second degree, which describe the oscillations of the weakly-coupled oscillators.
\end{abstract}			
			
\maketitle	
		
\section{Introduction} 
		
\subsection{Problem statement, results, methods.} \ 
The subject of this research is the two-dimensional periodic system with the hamiltonian unperturbed part and the small parameter
\begin{equation}\label{sv}
	\begin{cases} \dot x=\big(\gamma(y^3-y)+X^\nu(t,x,y,\e)\e\big)\e^\nu\\ \dot y=\big(\!-\!(x^3-x)+Y^\nu(t,x,y,\e)\e\big)\e^\nu \end{cases}\! 
	(\gamma\in(0,1],\, \e\in[0,\e_0],\, \nu=0,1), 
\end{equation}
where the functions $X^\nu,Y^\nu$ are continuous, $T$-periodic in $t;$ 
$$X^0=X_0^0(t,x,y)+X_1^0(t,x,y)\e+X_2^0(t,x,y,\e)\e^2,\ \ \ X^1=X_0^1(t,x,y)+X_1^1(t,x,y,\e)\e,$$ 
$Y^\nu$ has similar extension; the functions $X_0^\nu,Y_0^\nu,X_1^0,Y_1^0$ are uniformly with the respect to $t$ real-analytic in $x,y$ 
on $D_{\sigma,\sigma}^{x,y}=\{(t,x,y)\!:\,t\in \mathbb{R},$ $|x|,|y|<\sigma\},$ \,$\sigma>\sqrt{1+\gamma^{-1/2}};$\,   
the functions $X_{2-\nu}^\nu,Y_{2-\nu}^\nu\in C_{\,t,x,y,\e}^{0,1,1,\nu}(G_{\sigma,\sigma,\e_0}^{x,y,\e}),$ where 
$G_{\sigma,\sigma,\e_0}^{x,y,\e}=\{(t,x,y,\e)\!:\,t,x,y\in \mathbb{R},\,|x|,|y|\le \sigma,\,\e\in [0,\e_0]\}.$ 

\smallskip
{\it The continuous, $T$-periodic in $t$ function $Z(t,z_1,z_2)$ will be called uniformly with the respect to $t$ real-analytic in $z_1,z_2$ 
on $D_{z_1^0,z_2^0}^{z_1,\,z_2}=\{(t,z_1,z_2)\!:\,t\in \mathbb{R},|z_1|<z_1^0,\,|z_2|<z_2^0\}$ $(z_1^0,z_2^0>0),$ if for any $t\in \mathbb{R}$ 
and for any $|z_1|<z_1^0,\,|z_2|<z_2^0,$ the series $\displaystyle Z(t,z_1,z_2)=\sum\nolimits_{m,n=0}^\infty Z^{(m,n)}(t)z_1^m z_2^n$\,
with the real, continuous, and $T$-periodic in $t$ coefficients converges absolutely.}
		
\smallskip
Formula \eqref{sv} determines two different systems: one with $\nu=0,$ another with $\nu=1.$ 
Comparing these systems, we can say that the system with $\nu=1,$  which is usually called the standard system, has the "fast"\ time, 
because by reducing it to the system with $\nu=0$ we obtain period $T\e.$ In the case of the autonomous system, the change of the time variable allows to always consider the parameter $\nu=0.$
		
It is natural to refer to the autonomous system
\begin{equation}\label{snv}	\dot x=\gamma(y^3-y)\e^\nu,\quad \dot y=-(x^3-x)\e^\nu \qquad (\gamma\in(0,1],\ \nu=0,1). \end{equation}
as to the system of the first approximation or the unperturbed system with respect to \eqref{sv}.
		
System \eqref{snv} is Hamiltonian with the hamiltonian $$H(x,y)=(2x^2-x^4+\gamma(2y^2-y^4))\e^\nu/4.$$ 
Its phase plane consists of nine equilibrium points, closed orbits (cycles) and separatrices.
		
\begin{center} \includegraphics[scale=0.36]{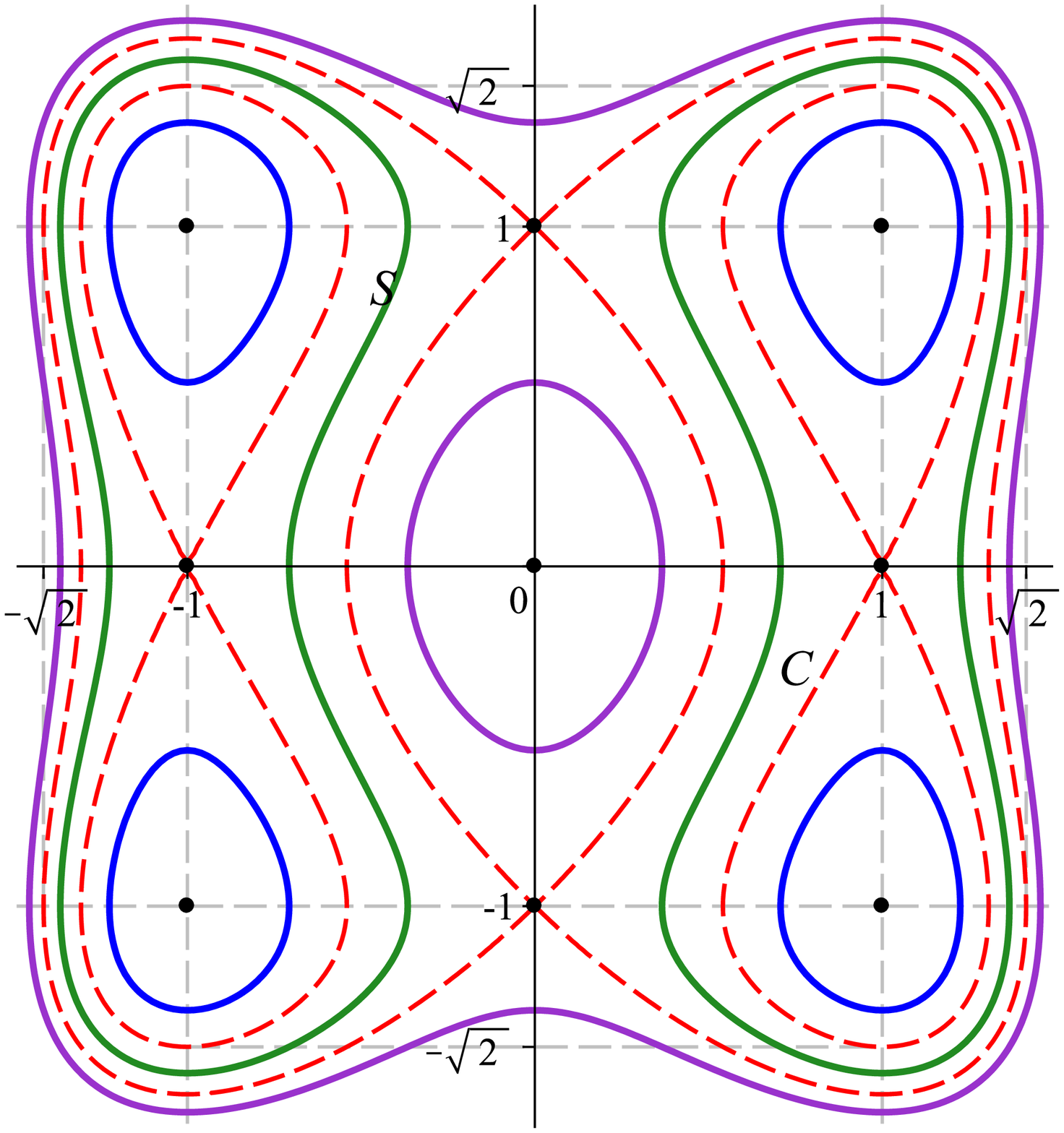}\quad \includegraphics[scale=0.36]{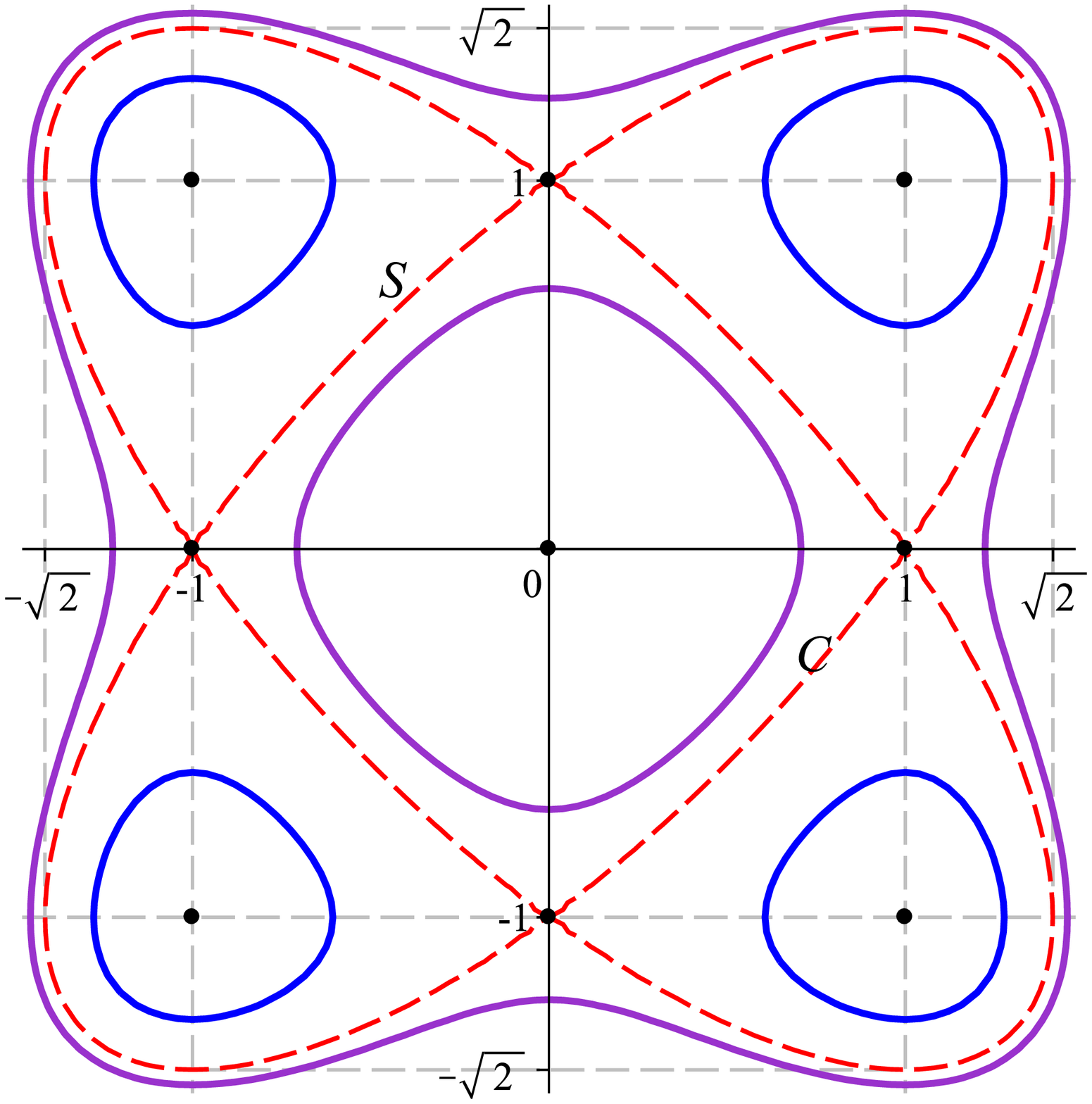}  \end{center} 
\begin{center} {\small \bf Fig.\,1.1. Phase portraits of the unperturbed system \eqref{snv} with $\gamma=1/2$ and $\gamma=1.$ } \end{center} 

The goal of this paper is to find, for system \eqref{sv} with any sufficiently small $\e>0,$ 
a certain number of two-dimensional invariant surfaces homeomorphic to the torus, which is obtained by factoring time with respect to the period. 
Each such surface is contained 
in a small with the respect to $\e$ neighborhood of the corresponding invariant surface of the system \eqref{snv} with the closed orbit, called the generating cycle (see definition\;5.2), as generatrix.
		
We explicitly write out conditions  (depending  on the parameter $\gamma$) on the unperturbed functions $X^\nu(t,x,y,0)$ and $Y^\nu(t,x,y,0),$ 
under which the perturbed system \eqref{sv}	has invariant surfaces described above, and obtain asymptotic expansions in powers of $\e$ 
for each of them (see theorem\;7.1).  
As an example, we provide a class of systems of the form \eqref{sv}, for which the average value of the perturbation for $\e=0$ 
does not have terms of order greater than three, and any perturbed system of this class for sufficiently small $\e>0$ has eleven invariant tori 
(see theorem\;8.1). 

As a result, the definition of Hilbert's cyclicity value can be generalized for the systems with a periodic perturbations. Its evaluation is also now possible in a sense of the evaluation, performed for the autonomous systems as part of the solution of the local Hilbert-Arnold problem, called  "weakened XVI Hilbert's problem" by V.\,I. Arnold in \cite{Arn}. The impossibility to perform change of the time variable, independent of the parameter, leads to the appearance of the two different classes of the researched systems: systems with so-called "fast" and "slow" time. In particular, these classes for the system \eqref{sv} are represented by the cases $\nu=1$ and $\nu=0.$ Let us notice, that reducing the periodic equations, which describe one or several weakly-coupled oscillators, to the Lienard systems, gives the realization to the case $\nu=1$ (see \cite{Bas4}).
		
Results of this paper are obtained using the method which we call the  generating tori splitting (GTS) method. The different stages of its application are described in sections 2--7.

Lets us overview the contents of this paper, which consists of the eight sections.

In the section 1 "Introduction" the short overview of the current results for the systems with the same or "similar" hamiltonians in the unperturbed systems is given. The majority of papers on this subject research the autonomous systems, therefore, instead of the invariant tori, the limit cycles are to be found, which simplifies the task. The main goal of these papers is the evaluation of the lower bound of the Hilbert's cyclicity value. The comparison of the results and possibilities of the applications of the existing methods, used for the autonomous systems, compared to the GTS method, are given.

In the section 2 "Parametrization of the orbits of the unperturbed system" the description of the orbits of the autonomous hamiltonian system \eqref{snv} is provided. There are three classes of the system's cycles: class $0]$ contains all cycles, similar to the purple cycles on the figure 1.1, class $1]$ --- green, class $2]$ --- blue. The parametrization of the cycles and separatrices is performed by using the solutions of the initial value problem of the system \eqref{cs}, which are the nonnormalized analogues of the generalized Lyapunov coordinates (see \cite{Lyap}, p.\,290). The intervals for the initial values are determined by the conditions \eqref{nd}. Also, section 2 contains the auxilliary lemmas describing the properties of the real-analytic two-periodic functions.

In the section 3 "The passing into a neighbourhood of a closed orbit" the construction and application to the system \eqref{sv} of the special polar change \eqref{pz}, which performs the passing into a neighbourhood of an arbitrary cycle of the system \eqref{snv}, observed from the equilibrium point it encircles, is given. Such change is possible, if the angular variable is changing monotonically. It was calculated that the angular variable does not change monotonically for all cycles of the class $1]$ (see lemma\;3.1). These cycles are found and excluded from the consideration by applying additional limitations \eqref{gamma*} on the initial values from \eqref{nd}.

In the section 4 "Primary radial averaging" the averaging of the right-hand side of the equation for the radial variable $r,$ independent of $\e,$  is performed for the special polar system \eqref{ps}, which has the arbitrarily chosen cycle of the unperturbed system as the equilibrium. As a result, the radial equation of the system \eqref{ps5}, obtained by performing the change \eqref{pave}, has the term, independent of $\e,$ with the order two or greater with the respect to $r.$ Then, after collecting the terms of the lesser order with the respect to $r$ and $\e,$ the system \eqref{ps5} is written in the form \eqref{fix}.

This allows us to write down the generating integral equation \eqref{pu}, by using the coefficients of the series $X_0(t,x,y),Y_0(t,x,y),$ in section 5. In the same section the admissable solutions are defined (see definition\;5.1). These solutions are denoted as $b_{kl}^\star$ and determine the initial values $0,b_{kl}^\star,l$ of the solutions parametrizing the so called generating cycles (see definition\;5.2) of the system \eqref{snv}.

Let us notice, that the definition of the admissable solutions in the GTS method requires not only standard nondegeneracy condition \eqref{ud} to hold, which is given constructively, but, also, due to the $T$-periodicity of the initial system, when $\nu=0,$ the Siegel's condition on $T$ and the calculated in \eqref{ome} periods $\omega_{kl}^\star$ of the solutions parametrizing these cycles.

Let us fix the arbitrary admissable solution $b_{kl}^\star,$ and at the same time the generating cycle, special polar system \eqref{ps}, which describes the motion in the small neighbourhood of this cycle, and the partially averaged system \eqref{fix}. The further considerations allow us to prove the existence of the invariant torus, which, for any small $\e$, the system \eqref{sv} preserves in a small vicinity of the generating cycle.

In the section 6 "The construction of the invariant surfaces" the averagings of radial and angular variables in terms, obtained from the perturbed part of the system \eqref{sv}, are performed in the system \eqref{fix}. After that the scaling is performed. The system \eqref{scsys2}, obtained after scaling, for any sufficiently small $\e$, (see theorem\;6.1) has two-dimensional invariant surface, which is parametrized by the continuous, $T$-periodic in $t$ and $\omega_{kl}^\star$-periodic in angular variable function $A_\e^\nu(t,\psi).$ The asymptotical stability of the invariant surface, the smoothness in $\psi,$ and continuity in $\e$ of the function $A_\e^\nu(t,\psi)$ are proven in lemmas 6.4\,--\, 6.6.

The ideas from the proof of the Hale's lemma in \cite{Hale} are used as the foundation of the provided proofs.

In the section 7 "Theoretical results" the results of the application of the GTS method to the system \eqref{sv} are provided. In particular, the two-dimensional invariant surfaces of the special polar system \eqref{ps} and the lesser terms of the asymptotic extensions in powers of the small parameter of these surfaces are written down (see lemma\;7.1 and corollary). After that, the substitution of the two-periodic parametrizations of the surfaces into the special polar change \eqref{pz} defines the formulas of the invariant tori, which are preserved by system \eqref{sv} for any sufficiently small $\e>0$ (see theorem\;7.1). Additionally, in the subsection 7.2 all obtained results are formulated for the autonomous systems to simplify its application in the autonomous case, when the case $\nu=1$ is not realized, and the Siegel's condition is not required.   

In the section 8 "Practical results" the analysis of the generating equation \eqref{pu} is performed, which allows to find its solutions for the specific periodic systems \eqref{sv} in both cases: $\nu=0$ and $\nu=1.$

For both $\nu=0$ and $\nu=1,$ the conditions \eqref{coeff} determine the set of systems \eqref{sv} with its periodic perturbations, independent of the parameter, having three terms of the order three or less. For such systems the eleven intervals, each containing one admissable solution of the generating equation \eqref{pu}, are found (see lemma\;8,1, theorem\;8.1). The generating cycles of the unperturbed system \eqref{snv} are plotted on the figure 8.2, by using the computing environment Maple. Thus, for $\nu=1$ the results of the paper \cite{LiQ} are repeated, but that paper only considers the autonomous systems and its perturbations had more terms. The results related to the existence of the eleven invariant tori for the case $\nu=0$ are new.

Additionally, the research related to the evaluation of the upper bound of the interval $(0,\e_*],$ where the limit cycles can be found, is performed. Theoretical results are confirmed visually by plotting with Maple the portraits of the integral curves for $\e_*=0.05,$ and with the help of the table of the intersections of spirals with the abscissa axis in the small neighbourhoods of the generating cycles for $\e_*=0.001.$ To illustrate the method the conditions \eqref{coeff} were applied to the autonomous system \eqref{sv}.
    
Let us note, that  in \cite{Bas4} the following simpler system has been studied with the GTS method:
\begin{equation}\label{loop} 
\dot x=\big(-y+\e X(t,x,y,\e)\big)\e^\nu,\ \ \dot y=\big(x^3-2\varpi x^2+\eta x+\e X_2(t,x,y,\e)\big)\e^\nu, 
\end{equation} 
where $\varpi=0,1;\ \eta=-1,0,1$ with $\varpi=0$ and $\eta\in\mathbb R^1$ with $\varpi=1,$
and the unperturbed part is Hamiltonian with \,$H=(x^4-8\varpi x^3/3+2\eta x^2+2y^2)\e^\nu/4.$
		
The system \eqref{loop} is a generalization of the system where  the  unperturbed part of the second equation 
is a homogeneous cubic polynomial with regard to $x,y$ and $\e,$ which was studied in \cite{Bas3} using an approach similar to the GTS method. 
In \cite{Bas3} the examples of the systems with three invariant tori were provided and the GTS method was demonstrated in full for the first time.
		
The GTS method is a modification of the method used by Y.\,N.\,Bibikov in \cite{Bib2} to study the  Duffing equation 
$\ddot x+x^3+\e bx=X(t,x,\dot x,\e).$ 
This equation characterizes the oscillations of the oscillator with an infinitesimal oscillation frequency with regard to the amplitude 
and is a special case of system \eqref{loop}.

\smallskip
{\bf Remark 1.} \ {\sl 
	In the majority of papers, devoted to the search of invariant surfaces in the systems with hamiltonian as an unperturbed part, these systems are considered autonomous. In this case, the preserved for any sufficiently small $\e$ invariant surfaces are cylindrical with the limit cycle as generatrix. It is located in the phase plane, in the small with the respect to $\e$ neighbourhood of a specific cycle of the unperturbed system, called generating cycle in this paper. }

\subsection{Connection to the weakened XVI Hilbert's problem.}  

The special cases of system \eqref{sv} are the subjects of the research with a similar setting. Let's focus on this research in this subsection.
		
Consider the autonomous system with hamiltonian unperturbed part
\begin{equation}\label{ha}
	\cfrac{du}{d\tau} = \cfrac{\partial \widetilde H(u,v)}{\partial v} + \varepsilon P(u,v),\quad 
	\cfrac{dv}{d\tau} = -\cfrac{\partial \widetilde H(u,v)}{\partial u} + \varepsilon Q(u,v), 
\end{equation}
where $\widetilde H,P,Q$ are polynomials with real coefficients of the corresponding degree $n,m_1,m_2,$ $\varepsilon\ge0$ is a small parameter.

The reader can acquaint themselves with the answers on the first question of the weakened XVI Hilbert's problem, named so by V.\,I. Arnold in \cite{Arn}, in \cite{Var},\,\cite{Ilyas}. They provide the upper bounds of the number $N(n, m) = N(n, m, \widetilde H, P, Q),$ which is for any sufficiently small $\e > 0$ the exact upper bound of the number of the limit cycles of the considered system.

The further questions of the local Hilbert-Arnold problem, as it became known later, are related to the evaluation of the lower bound of the number $N(n, m)$ or Hilbert's cyclicity value $H(n)=N(n-1,n-1)$ for the system \eqref{ha} and the classification of all possible relative locations of limit cycles.
		
The estimation of the $N(n, m)$ is usually performed using the Melnikov function method (see \cite{Mel,GH}).
This method is based on the Poincar\'e-Pontryagin theorem (see theorem\;6.1,\cite{Li}) and comes down to researching 
and calculating Abelian integrals. According to this method, the function 
\,$d_{\varepsilon}(h)=\varepsilon M_1(h)+\varepsilon^2 M_2(h)+\ldots + \varepsilon^k M_k(h) + \ldots\,$\,    
is considered, where the  parameter $h$ represents the value of the function $H(x, y)$ in a neighborhood of an equilibrium point 
and $M_k(h)$ are autonomous Melnikov functions of order $k.$
The number of zeros with their  multiplicities accounted for  the first nonvanishing Melnikov function  $M_k(h)\ (k=1,2,3,\ldots)$ determines 
the highest number of limit cycles bifurcating from the period annulus.
		
Using  the  Melnikov function method the results for some  Lienard systems with $n = 4$ and $m =3$ were obtained 
in \cite{ILY},\,\cite{Ili1},\,\cite{Du}, with $n = 7,\ m = 4$ --- in \cite{China}, with $n = 4$ and $m=\overline {3,100}$ --- in \cite{WTX}.

The synopsis of the existing results dating to 2003 is provided  by J.\ Li in \cite{Li}.
		
Another approach was chosen by J.\,Li and Q.\,Huang in \cite{LiQ}. 
By using the detection function method it was proved that $H(3)\ge 11$ for the system of type \eqref{ha} 
\begin{equation} \label{sva} 
 \cfrac{du}{d\tau} = v(1-\mu_1 v^2) + \varepsilon u(k_1 v^2 + k_2 u^2 - \lambda),\ 
 \cfrac{dv}{d\tau} = -u(1-\mu_2 u^2) + \varepsilon v(k_1 v^2 + k_2 u^2 - \lambda) 
\end{equation}
with \,$\widetilde H(u,v) = (\mu_1 v^4 - 2v^2  + \mu_2 u^4 - 2u^2)/4$\, and \,$0<\mu_1<\mu_2.$ 
		
		Using the results of \cite{LiQ} J.\,Li in \cite{Li} specified the lower bound for $H(n)$ for $n=2^k-1\ (k\ge3).$ 
		He considered the sequence of $\mathbb{Z}_q$-equivariant systems $(q\ge2)$  
		$$\cfrac{du}{d\tau} = \cfrac{\partial \widetilde H_k(u,v)}{\partial v} + \varepsilon P_k(u,v),\ \
		\cfrac{dv}{d\tau} = \cfrac{\partial \widetilde H_k(u,v)}{\partial u} + \varepsilon Q_k(u,v)$$ 
		with $\widetilde H_{k} (u,v) = \widetilde H_{k-1}(u^2 - \eta^{k-2},v^2 - \eta^{k-2}),$ 
		$P_{k} (u,v) = P_{k-1}(u^2 - \eta^{k-2},v^2 - \eta^{k-2}),$ 
		$Q_{k} (u,v) = Q_{k-1}(u^2 - \eta^{k-2},v^2 - \eta^{k-2}),$ where $\widetilde H_2,P_2,Q_2,$ were introduced in \eqref{sva}.
		
		It was done using the method of "quadrupling" of the system, suggested by the C.\,J, Christopher and N.\,G. Lloyd in \cite{CL}.
		
		Later on, the authors  of \cite{LLY} constructed the system \eqref{ha} with the more complicated Hamiltonian
		\,$\widetilde H(u,v)=(u^4+4(1-\lambda)u^3/3-2\lambda u^2+(v^2-k^2)^2-k^4)/4,$\, $P\equiv 0,$ $Q=v(a_1 + a_2 u$ $+ a_3 u^2 + a_4 v^2),$ 
		which  has thirteen limit cycles for certain values of parameters, i.\,e.\,$H(3)\ge 13.$  		

		\subsection{The comparison of the results obtained by different methods.} 
		
		Let us compare the results obtained in this paper with the existing ones.
		
		In \cite{LiQ}, for system \eqref{sva} with $\mu_1=1,\,\mu_2=2,\,k_1=-3,\,k_2=1$ the classification  specifying  the number of limit cycles and their location  depending on the value of the parameter $\lambda$ was provided. 
				In particular, it was shown that system \eqref{sva} with \,$-4.80305+O(\e)<\lambda<-4.79418+O(\e)$ has eleven limit cycles.
		
		The system \eqref{sva} is a special case of the system \eqref{sv}. Choosing these  values of parameters $\mu$ and $k$ 
and  using the change $\tau = (1/2)^{1/2}\,t,\ u = 2^{-1/2}y,\ v=x,$ it is transformed into the system \eqref{sv} with 
$$\nu =0,\ \gamma = 1/2,\ \ X(t,x,y,\e)=2^{-1/2}x(x^2-3y^2/2-\lambda),\ Y(t,x,y,\e)=2^{-1/2}y(x^2-3y^2/2-\lambda),$$  
i.\,e. in the system \eqref{sv} $X=X_0(x,y),$ $Y=Y_0(x,y).$
		
		Using the  GTS method it was figured out that the system \eqref{sva} has eleven admissible solutions of the generating equation \eqref{pu} 
(see definition 5.1). These solutions determine points on the abscissa  axis, and, according to  Theorem 1 stated below, a limit cycle that generates the invariant torus passing through a small neighborhood of each of these points. 
Thereby,  it gives another proof of the result of \cite{LiQ}  and its  generalization in the following sense: 

		a) it allows  to add  any function of $\e$ into the perturbation; 
		
		b) make systems \eqref{sv} time-periodic 
while preserving mentioned above functions $X$ and $Y$ as a mean values;
		
		c) it gives the possibility  to apply the  GTS method to standard systems, i.\ e., systems with $\nu=1,$ which were not considered in  \cite{LiQ}.
		
On the other hand, in Section 6.2 we describe a  set of systems \eqref{sv} with $\gamma = 1/2$ and time-periodic perturbations \eqref{coeff}, 
which also have eleven admissible solutions of the generating equation that determine invariant tori.
		
In the special case, when in the constructed autonomous system  
$$\nu=0;\quad X(t,x,y,\e)=0,\ \ Y(t,x,y,\e)=-3.314y-0.361y^3 + 4.493x^2 y,$$ 
the change \,$t = \sqrt{2}\tau,\ x = v,\ y=\sqrt{2}u$\, transforms it into system \eqref{sva} with 
$$\mu_1 = 1,\ \mu_2 = 2;\ \ P(u,v)=(4.493\sqrt{2}v^2 - 0.722\sqrt{2}u^2-3.314)u,\quad Q(u,v)=0,$$ 
also having eleven limit cycles, which were not considered in \cite{LiQ}. 

\subsection{GTS method features.} \ 
The GTS method allows to write down the generating  (or bifurcational) equation  for any time-periodic systems with a Hamiltonian unperturbed part. 
Any solution of the generating equation  which lies in the preadjusted boundaries, determines for any sufficiently small $\e>0$ a limit cycle in the autonomous case
or an invariant torus in the periodic case, assuming that: firstly, for $\nu=0$ Siegel's condition \eqref{zig} on the period of the perturbation and on the period of the encircling time of the generating cycle holds, and, secondly,
 the nondegeneracy condition \eqref{ud} holds, which is the inequality to zero of the certain constant or nondegeneracy of a constant matrix in polydimensional case. 
		
The main advantage of the GTS method is its versatility. In comparison to other mentioned methods, it allows:

		1. To consider not only autonomous but time-periodic systems determining invariant tori of different dimensions. Apparently, study of periodic systems with the application of Melnikov function method is halted by comptutation difficulties. One of the attempts to study such systems with the Melnikov function method known to authors is performed by M.\ Han \cite{Han}, but the studied system (1.12) has only time-periodic perturbations as a coefficients of $\e^2,$ all Melnikov function generated by the coefficients of $\e$ are autonomous.
		
2. To consider the so-called standard system, i.\,e. system \eqref{sv}, 
which has $\varepsilon$ as a factor in its right-hand side, or "fast" time. 
For such systems, the analogue of $H(n)$ can be introduced, and its lower bound can be found. It is already done in this paper. 
For such systems Siegel's condition is not required.
		
		3. To use in the system \eqref{sv} any smooth in $x, y$ perturbations $X_{2-\nu}^\nu(t,x,y,\e),Y_{2-\nu}^\nu(t,x,y,\e)$ and uniformly with the respect to $t$ analytic in $x, y$ perturbations $X_1^0(t,x,y),Y_1^0(t,x,y)$ in a specific neighbourhood of the point $x=0, y=0,$ because, as it was figured out, they don't affect the existence of the invariant tori.
		
		4. To write down asymptotic expansions of limit cycles and invariant tori with respect to the $\e.$
		
		5. To research higher dimensional periodic systems and to determine the bifurcations of the tori of various dimensions 
		that branch out from the equilibrium point for any sufficiently small parameter value.
				
				For instance  in \cite{BZ} the system of the dimension $2n\ (n\ge 2)$
		$$\dot x_i=-\gamma_{i}y_i\e^{\nu}+X_i(t,x,y,\e)\e^{1+\nu},\
		\dot y_i=\gamma_{i}(x_i^3-\eta_{i} x_i)\e^{\nu}+Y_i(t,x,y,\e)\e^{1+\nu},$$
		where $i=\overline{1,n},\ \nu=0,1,$ $\gamma_{i}\in (0,+\infty),$ $\eta_{i}=-1,0,1,$ $\e\in[0,\e_0)$	 was studied. 
		The unperturbed part of this system  has from one to  $3^n$ zeros.
		Using the GTS method  the conditions on the average values of the periodic perturbations of the system were obtained. 
		If these conditions hold, then the system preserves invariant tori of dimension $n+1.$ 
		
		As an example, 
		the four-dimensional system was provided. For this system, 
		six points were found such that unique three-dimensional surface homeomorphic to torus passes through the small neighbourhood of each point.

\section{Parametrization of the orbits of the unperturbed system} 
		
\subsection{Phase Portrait.} \  
Consider an autonomous hamiltonian system with nine equilibrium points similar to system \eqref{snv} defined by the equations
\begin{equation}\label{cs} 
	C'(\varphi)=\gamma(S^3(\varphi)-S(\varphi)),\ \ S'(\varphi)=-(C^3(\varphi)-C(\varphi))\quad (0<\gamma\le 1).
\end{equation}
		
For each $a>0,$ we consider the set $\Gamma_a$ of closed orbits on the plane $(C,S)$ 
determined by the integral of system \eqref{cs} given by the equation
\begin{equation}\label{oi} 	(C^2-1)^2+\gamma(S^2-1)^2=a. 	\end{equation}
		
Obviously, systems \eqref{cs} and \eqref{snv} have the same orbits, but for $\nu=0$ and  $\nu=1$ their parametrizations are  different.
		
For system \eqref{cs} five of nine singular points are centers. Points $(\pm 1,1),$ $(\pm 1,-1)$ are solutions of \eqref{oi} when $a=0,$ 
and $(0,0)$ is the solution of \eqref{oi} when $a=1+\gamma.$ The other singular points are saddle points; 
points $(\pm 1,0)$ are the solutions of \eqref{oi} when $a=\gamma,$ and $(0,\pm 1)$ are the solutions of \eqref{oi} when $a=1.$ 
		
It is sufficient  to describe orbits or parts of orbits that lie in the first quadrant, 
because equation \eqref{oi} is invariant with  respect to the change  $C \to -C$,  $S \to -S$. We denote such set of orbits by $\Gamma_a^*.$ 
		
The extremal values of the curves that belong to $\Gamma_a^*$ are:
\begin{equation}\label{exc}	\begin{matrix} 
		r^s=\sqrt{1+\gamma^{1/2}},\ \ l^s=\sqrt{1-\gamma^{1/2}},\ \ 
		r^i=\sqrt{1-(1-\gamma)^{1/2}},\ \ r^e=\sqrt{1+(1-\gamma)^{1/2}},\hfill\\ 
		u^e=\sqrt{1+\gamma^{-1/2}};\ \ l_1^1,l_2^1=\sqrt{1-a^{1/2}},\ \ lo_2^1=\sqrt{1-(a/\gamma)^{1/2}},\hfill \\
		r_0^{0i},l_1^0=\sqrt{1-(a-\gamma)^{1/2}},\ \ r_0^{0e},r_1^0=\sqrt{1+(a-\gamma)^{1/2}},\ \ 
		r_0^{1e},r_1^1,r_2^1=\sqrt{1+a^{1/2}},\hfill\\ 
		u_0^{0i}=\sqrt{1-((a-1)/\gamma)^{1/2}},\ u_0^{0e}=\sqrt{1+((a-1)/\gamma)^{1/2}},\ u_0^{1e},u_1^1,u_2^1=\sqrt{1+(a/\gamma)^{1/2}}, 
\end{matrix} \end{equation}		
where  $r,\ l,\ u,\ lo,\ s,\ i,\ e$ means right, left, upper, lower, separatrix, internal and  external, respectively,
the first five constants characterize $\Gamma_\gamma^*$ and $\Gamma_1^*,$ 
for the other constants the superscript determines the value of the other coordinate, 
while the  subscript determines the  number of a class which we will introduce below.
	
We first  consider the separatrix curves determined by equation \eqref{oi}, namely, 
the curves passing  through the singular points of the system \eqref{cs}.
		
For $\gamma<1$  the set $\Gamma_\gamma$ $(a=\gamma)$  consists of four closed curves (or two eight-loops). 
Two of those homoclinic orbits contact each other at the singular point $(1,0)$ and enclose points $(1,\pm 1),$ 
other two contact each others at the singular point $(-1,0)$ and enclose singular points $(-1,\pm 1).$
Thus, $\Gamma_\gamma^*$ is a top part of right eight-loop and has the following extremal points: $(r^s,1),\,(1,\sqrt 2),\,(l^s,1),\,(1,0).$
		
For $\gamma<1$ the  set $\Gamma_1$ $(a=1)$ consists of two closed curves: the internal curve $\Gamma_1^i$ and the external one $\Gamma_1^e.$ 
These curves contact each other at the singular points $(0,1)$ and $(0,-1).$  The set $\Gamma_1^i$ encloses the singular point $(0,0)$,  
and one of each eight-loops  mentioned above is located in the two areas between $\Gamma_1^i$ and $\Gamma_1^e.$ 
The set $\Gamma_1^{i*}$ has the extremal points $(r_i,0),\,(0,1),$ 
and $\Gamma_1^{e*}$ has the extremal points $(r^e,0),\,(\sqrt 2,1),\,(1,u^e),\,(0,1).$ 
		
For $a=\gamma=1$ the parameters $r^i=r^e=1$ and $\Gamma_\gamma$ coincides with $\Gamma_1,$ 
namely, $\Gamma_1^{i*}$ and $\Gamma_1^{e*}$ contact each other at the point $(1,0)$ creating the  top part of the right eight-loop.
		
Now, for $a\ne \gamma,1$ the set $\Gamma_a$ consists of closed orbits (cycles) of system \eqref{cs}.
The separatrices divide the set $\Gamma_a$ into three classes, which we denote by 0],\,1],\,2]	and define as follows:
				
0]\, $a>1.$ \ For  $1<a<1+\gamma$ the set $\Gamma_a$ consists of two cycles: 
the inner orbit $\Gamma_a^i$ which encloses $(0,0)$ and lies inside $\Gamma_1^i,$
and the outer  orbit $\Gamma_a^e$ which encloses $\Gamma_1^e.$
Then the points $(r_0^{0i},0),$ $(0,u_0^{0,i})$ are extremal for $\Gamma_a^{i*},$ 
the points $(r_0^{0e},0),$ $(r_0^{1e},1),$ $(1,u_0^{1e}),$ $(0,u_0^{0e})$ are extremal for $\Gamma_a^{e*},$ 
moreover, $r_0^{0i}\in (0,r^i),$ $r_0^{0e}\in (r^e,\infty).$  
Since  $\Gamma_a^i$ degenerates into	the  point $(0,0)$ for $a=1+\gamma,$ for $a\ge 1+\gamma,$ only the cycle $\Gamma_a^e$ exists.
		
As the  result, the class 0] naturally splits into two subclasses: 
\,$0^i]$ for cycles $\Gamma_a^i$ with $1<a<1+\gamma$\, and \,$0^e]$ for cycles $\Gamma_a^e$ with $a>1.$
		
1]\, $\gamma<a<1.$ The set $\Gamma_a$ consists of two cycles. The top half of  the right orbit $\Gamma_a^*,$ 
is located between $\Gamma_\gamma^*$ and $\Gamma_1^{i*}\cup\Gamma_1^{e*}$\, 
and	has the  following extremal points: $(r_1^0,0),\,(r_1^1,1),\,(1,u_1^1) $ $(l_1^1,1),\,(l_1^0,0),$ where $r_1^0\in (1,r^e).$
		
2]\, $0<a<\gamma.$ The  set $\Gamma_a$ consists of four cycles.
The orbit $\Gamma_a^*$ encloses $(1,1)$ and lies inside of $\Gamma_\gamma^*,$
its extremal points are $(r_2^1,1),\,(1,u_2^1),\,(l_2^1,1),\,(1,lo_2^1),$ where $r_2^1\in (1,r^s).$ 

\subsection{Parametrization of the cycles.} 
The parameter $a$ does not determine a specific cycle from $\Gamma_a,$ therefore it cannot be used  to parametrize the cycles. 
To define the parametrization  we use the extremal points: $(r_0^{0i},0),(r_0^{0e},0)$ in class 0], \,$(\pm r_1^0,0)$ in class~1], \,$(\pm r_2^1,1)$ and $(\pm r_2^1,-1)$ in class 2]. 
Each such point defines the corresponding cycle and the parameter $a$ from \eqref{oi} can be explicitly expressed through parameters $r$ mentioned above.

Thus, an arbitrary cycle of system \eqref{cs} is parameterized by the real analytic $\omega_{kl}$-periodic solution $CS_{kl}(\varphi)=(C(\varphi),S(\varphi))$ of the initial value problem 
with the initial conditions $0, b_{kl}, l,$ i.\,e.
\begin{equation}\label{nd} C(0)= b_{kl},\ \ S(0)=l \quad (k,l=0,\pm 1,\ (k,l)\ne (0,\pm 1)); \end{equation}
where  $b_{00}=\left[\begin{matrix} b_{00}^i=r_0^{0i}\in (0,r^i)\hfill \\ b_{00}^e=r_0^{0e}\in (r^e,r^\sigma)\end{matrix}\right.;$
$b_{10}=r_1^0\in (1,r^e),$ $b_{-10}=-r_1^0;$ $b_{1,\pm 1}=r_2^1\in (1,r^s),$ $b_{-1,\pm 1}=-r_2^1;$\,  
$r^\sigma = (1+(\gamma\sigma^2(\sigma^2 - 2))^{1/2})^{1/2}.$ 
 
In the formulas given above the $k$ determines the shift\ along the abscissa axis, the $l$ determines the "shift"\ along the ordinate axis 
and the number $|k|+|l|$ determines the class number related to the parameterized cycle; the constants $r$ are introduced in \eqref{exc}.

The right-hand side of the system \eqref{cs} is real and polynomial, therefore, according to Cauchey's theorem, the solution $CS_{kl}(\varphi)=CS_{kl}(\varphi, 0, b_{kl}, l))$ has the mentioned properties.

Denote by $a_{kl}$ the value of the parameter $a$ corresponding to the  initial value $b_{kl}$ and recalculate required extremal values from \eqref{exc}. A simple computation yields
\begin{equation}\label{ndb}	 \begin{matrix}
 		a_{kl}=(1-|l|)\gamma +(b_{kl}^2-1)^2;\ u_0^{1e}=(1+((\gamma + ((b_{00}^e)^2 - 1)^2)\gamma^{-1})^{1/2})^{1/2},\\ 
 		r_{0}^{1e} = (1 + (((b_{00}^e)^2 - 1)^2 + \gamma)^{1/2})^{1/2},\ r_{1}^{1} = (1 + ((b_{10}^2 - 1)^2 + \gamma)^{1/2})^{1/2},\\
 		l_1^0 = (2 - b_{01}^2)^{1/2},\ l_2^1 = (2 - b_{11}^2)^{1/2},\ \ l_{1}^{1} = (1 - ((b_{10}^2 - 1)^2 + \gamma)^{1/2})^{1/2}.\hfill 
\end{matrix}  \end{equation}

	{\bf Remark 2.1.} \ {\sl The limitation on the constant $\sigma:\ \sigma > \sqrt{1+\gamma^{-1/2}},$ which determines the domain of system \eqref{sv} is chosen in a such way so all class $0^e]$ cycles and, by extension, other cycles of unperturbed system \eqref{snv} are contained in the domain. 
		According to \eqref{nd}, the parameter $b_{00}^e\in(r^e, r^{\sigma}),$ therefore, according to \eqref{ndb}, $a_{00}^e=\gamma+((b_{00}^e)^2 - 1)^2 < \gamma + ((r^\sigma)^2 - 1)^2 = \gamma(1 + \sigma^2(\sigma^2 - 2)).$
		Moreover, for any class $0^e]$ cycle the following statement is correct:
		$$\max\,|C(\varphi)|=r_0^{1e}\le u_0^{1e}=\max\,|S(\varphi)|=(1+(a_{00}^e/\gamma)^{1/2})^{1/2}<\sigma.$$}

\subsection{Calculation of the  periods.} \
Let us introduce the auxiliary functions  
\begin{equation}\label{Spm} S_{\pm}(C^2(\varphi))=\sqrt{1\pm\gamma^{-1/2}(a_{kl}-(C^2(\varphi)-1)^2)^{1/2}}.\end{equation}

{\bf Proposition 2.1.} The period $\omega_{kl}$ of the real-analytic periodic solution $(C(\varphi),S(\varphi))$ of the initial value problem with the initial values from \eqref{nd} can be computed using the following formulae
\begin{equation}\label{ome}
\omega_{00}^i=4\varphi^i_-,\ \omega_{00}^e=4(\varphi^e_+ + \varphi^e_-),\  
		\omega_{\pm 10}=2(\varphi^l_- + \varphi^u_+ + \varphi^r_-),\ \omega_{1,\pm 1}=\omega_{-1,\pm 1}=\varphi^2_- + \varphi^2_+, 
\end{equation} 
where 
$\displaystyle \varphi^i_-=\int_{b_0^{0i}}^0 \zeta_- dC;$ \ $\displaystyle \varphi^e_+=\int_0^{r_0^{1e}} \zeta_+ dC,$ \
$\displaystyle \varphi^e_-=\int_{r_0^{1e}}^{b_0^{0e}} \zeta_- dC;$ \ $\displaystyle \varphi^l_-=\int_{l_1^0}^{l_1^1} \zeta_- dC,$ \
$\displaystyle \varphi^u_+=\int_{l_1^1}^{r_1^1} \zeta_+ dC,$ \ $\displaystyle \varphi^r_-=\int_{r_1^1}^{b_{10}} \zeta_- dC;$ \
$\displaystyle \varphi^2_-=\int_{b_{11}}^{l_2^1} \zeta_- dC,$ \ $\displaystyle \varphi^2_+=\int_{l_2^1}^{b_{11}} \zeta_+ dC,$\, 
the limits of integration are defined in \eqref{exc}, with  $a=a_{kl}$ from \eqref{ndb}, 
and $\zeta_{\pm}(C^2)=(\gamma(S_{\pm}^3(C^2)-S_{\pm}(C^2)))^{-1}.$

{\it Proof}.\, Describing the motion along the cycle where $\varphi$ is changing from $0$ to $\omega_{kl}$, for each class, 
function $S(\varphi)$ can be represented as the function $S(C(\varphi)).$ By using formula \eqref{oi}, constants from \eqref{nd} and \eqref{ndb} 
and the fact that in system \eqref{cs} $S'(0)>0$ for all cycles from the class $0^i]$ and $S'(0)<0$ for all cycles from the other classes, we have:
$$\begin{matrix}
	0^i]\ S(C(\varphi)) = \{ 
	S_-  \hbox{ when } C\!\downarrow_{-b_{00}^i}^{b_{00}^i},\,
	-S_- \hbox{ when } C\!\uparrow_{-b_{00}^i}^{b_{00}^i}\};\hfill \\
	0^e]\ S(C(\varphi)) = \{ 
	-S_- \hbox{ when } C\!\uparrow_{b_{00}^e}^{r_0^{1e}},\,
	-S_+ \hbox{ when } C\!\downarrow_{-r_0^{1e}}^{r_0^{1e}},\,
	-S_- \hbox{ when } C\!\uparrow_{-r_0^{1e}}^{-b_{00}^e},\,\hfill\\
	\hfill S_-  \hbox{ when } C\downarrow_{-r_0^{1e}}^{-b_{00}^e},\,
	S_+  \hbox{ when } C\!\uparrow_{-r_0^{1e}}^{r_0^{1e}},\,
	S_-  \hbox{ when } C\!\downarrow_{b_{00}^e}^{r_0^{1e}}\}; \\
	1]\ S(C(\varphi)) = \{ 
	-kS_- \hbox{ when } kC\!\uparrow_{kb_{k0}}^{kr_1^1},\,
	-kS_+ \hbox{ when } kC\!\downarrow_{kl_1^1}^{kr_1^1},\,
	-kS_- \hbox{ when } kC\!\uparrow_{kl_1^1}^{kl_1^0},\,\hfill\\
	\hfill kS_-  \hbox{ when } kC\!\downarrow_{kl_1^1}^{kl_1^0},\,
	kS_+  \hbox{ when } kC\!\uparrow_{kl_1^1}^{kr_1^1},\,
	kS_-  \hbox{ when } kC\!\downarrow_{kb_{k0}}^{kr_1^1}\}; \\
	2]\ S(C(\varphi)) = \{
	lS_-  \hbox{ when } kC\!\downarrow_{kl_2^1}^{kb_{kl}},\,
	lS_+  \hbox{ when } kC\!\uparrow_{kl_2^1}^{kb_{kl}}\ \ (kl=1);\hfill \\ 
	\hfill lS_-  \hbox{ when } kC\!\uparrow_{kl_2^1}^{kb_{kl}},\,
	lS_+  \hbox{ when } kC\!\downarrow_{kl_2^1}^{kb_{kl}}\ \ (kl=-1)\}. 
\end{matrix} $$

We can write down the first equation of system \eqref{cs} as  
$$d\varphi=(\gamma(S^3(\varphi)-S(\varphi)))^{-1}dC(\varphi).$$

Then we integrate this equation with the respect to $\varphi$ from $0$ to $\omega_{kl}$ and substitute the formula for $S(C(\varphi))$ in the right-hand side of this equation for each class. By reducing the amount of terms, considering that functions $\zeta_{\pm}$ are even with the respect to $C,$ we obtain the formulas \eqref{ome}. $\Box$

Moreover, we discover that $C(\omega^i_{00}/4)=C(\omega^e_{00}/4)=0,$ $S(\omega^i_{00}/4)=u_0^{0i},$ 
$S(\omega^e_{00}/4)=-u_0^{0e},$ $C(\omega_{10}/2)=l^0_1,$ $S(\omega_{10}/2)=0.$

\smallskip
{\bf $2.4.$  Real-analytic two-periodic functions.}
In this article we will consider real-analytic $\omega_{kl}$-periodic functions, other than the solution $C\!S_{kl}(\varphi)$ of system \eqref{cs}, and two-periodic in $t$ and $\varphi$ functions. Let us describe their properties.

\smallskip
{\it Denote continuous, $T$-periodic in $t$ and $\omega$-periodic in $\varphi$ function $\zeta(t,\varphi)$ as uniformly with respect to $t$ 
real-analytic with respect to $\varphi$ function, if for all $\varphi_*\in \mathbb{R}$ such parameter $r_{\varphi_*}>0$ exists, that the series 
\,$\displaystyle \zeta(t,\varphi)=\sum\nolimits_{m=0}^\infty \zeta^{(m)}(t) (\varphi-\varphi_*)^m$\, 
with real, continuous and $T$-periodic in $t$ coefficients converges absolutely for any $t\in\mathbb{R},$ 
for any such $\varphi\in \mathbb{C},$ that $|\varphi-\varphi_*|<r_{\varphi_*}.$ }

\smallskip
Introduce two sets for all $\varrho>0.$
\begin{equation}\label{BB} 
	 B_{\varrho}=\{({\rm Re}\,\varphi,{\rm Im}\,\varphi)\!:\,{\rm Re}\,\varphi\in \mathbb{R},\,|{\rm Im}\,\varphi|\le \varrho\},\ \  
		\mathcal B_{\varrho}=\{(t,\varphi)\!:\, t\in \mathbb{R},\,\varphi\in B_{\varrho}\}.
\end{equation}

{\bf Lemma 2.1.}  
{\it Assume that the function $\zeta(t,\varphi)$ is continuous, $T$-periodic in $t,$ $\omega$-periodic in $\varphi,$ uniformly with respect to $t$ real-analytic with respect to $\varphi.$ Assume that the number $M^r=\max_{\,t,\varphi\in \mathbb{R}}|\zeta(t,\varphi)|.$ Then 
	
	$1)$ for any $\delta>0$ such parameter $\varrho>0$ exists, that the function $\zeta(t,\varphi)$ is two-periodic, uniformly with respect to $t$ analytic function on $\mathcal B_{\varrho}$ and 
	$$\exists\,M\ \ (0<M\le M^r+\delta)\!:\ M=\max\nolimits_{\,(t,\varphi)\in {\mathcal B_\varrho}}|\zeta(t,\varphi)|; $$ 
	
	$2)$ for $\zeta(t,\varphi)$ exponential form of Fourier series' expansion with $\varphi\in[0,\omega]$
	\begin{equation}\label{varpi} \zeta(t,\varphi)=\sum_{n=-\infty}^{+\infty} \zeta_n(t) e^{2\pi i n \varphi/\omega},\quad 
		\zeta_n(t)=\omega^{-1}\int_0^\omega \zeta(t,\varphi)e^{-2\pi i n \varphi/\omega}d\varphi, \end{equation}
	where $\zeta_{-n}(t)=\overline{\zeta_n}(t),$ due to $\zeta(t,\varphi)$ being the real function, the following is correct: }
\begin{equation}\label{cofu} \forall\,t\in \mathbb{R}\!:\ \ |\zeta_n(t)|\le Me^{-2\pi|n|\varrho/\omega}. \end{equation}

{\it Proof.}\, 1) For each $\varphi\in \mathbb{C}$\, with ${\rm Re}\,\varphi\in [0,\omega]$ denote $K_{r_\varphi}(\varphi)$ 
as a convergence disk of the power series $\zeta(t,\varphi)$ with the convergence radius $r_\varphi>0$ and the center $({\rm Re}\,\varphi,0).$ 
Denote $K_{r_0}(\varphi_0),\ldots,K_{r_m}(\varphi_m)$ as a finite cover by such disks of the closed interval 
$\{{\rm Re}\,\varphi\in [0,\omega],$ ${\rm Im}\,\varphi=0\}.$

Assume $\varrho_0=\min\,\{{\rm Im}\,z_1,\ldots,{\rm Im}\,z_m\},$ where $z_\mu$ with ${\rm Im}\,z_\mu>0$\, \,$(\mu=\overline{1,m})$ 
are intersection points of the circles $|z-\varphi_{\mu-1}|=r_{\mu-1}$ and $|z-\varphi_\mu|=r_\mu.$ 
Then, choose $\varrho$ as any number, lesser than $\varrho_0.$ 

For all $t$ the function $\zeta(t,\varphi)$ is $\omega$-periodic on $B_\varrho,$  
because the analytic function $\zeta(t,\varphi+\omega)-\zeta(t,\varphi)\equiv 0$ 
with $\varphi\in\{{\rm Re}\,\varphi\in [0,\omega],\ {\rm Im}\,\varphi=0\}.$
Therefore, due to uniqueness, it is equal to zero for all $\varphi\in B_\varrho.$ Thus, the continuous two-periodic function $|\zeta(t,\varphi)|$ reaches its maximum value, denoted as $M,$ at the $\mathcal B_\varrho.$ 
Additionally, the constant $\varrho$ can be chosen so small, that $M$ is $\delta$-close to $M^r.$

2) For any fixed $t$ the contour integral with a contour
$$0\to \omega\to \omega+i\varrho\,{\rm sgn}\,n\to i\varrho\,{\rm sgn}\,n\to 0$$ 
from $B_\varrho$ of the function $\zeta(t,\varphi)e^{-2\pi i n \varphi/\omega}$\, is equal to zero, 
due to $\zeta(t,\varphi)$ being the $\omega$-periodic function.
The sum of integrals along the vertical lines of the same function is also equal to zero for the same reason.
Therefore, the sum of integrals along the horizontal lines is equal to zero as well. Thus, due to \eqref{varpi}, 
$$\omega\zeta_n=\!-\!\int_{\omega+i\varrho\,{\rm sgn}\,n}^{i\varrho\,{\rm sgn}\,n}\! \zeta(t,\varphi)e^{-2\pi i n \varphi/\omega}d\varphi= 
e^{2\pi|n|\varrho/\omega}\! \int_0^\omega\! \zeta(t,\psi+i\varrho\,{\rm sgn}\,n)e^{-2\pi i n \psi/\omega}d\psi.$$ 
The evaluation of the function $|\zeta_n(t)|$ in formula \eqref{cofu} is implied from this equality. \ $\Box$

\smallskip
It is obvious, that any decrease of the constant $\varrho$ does not affect the evaluations and "brings $M$ closer"\ to $M^r.$ 
Moreover, Lemma\;2.1 can be applied to any real-analytic, $\omega$-periodic functions, independent of $t,$  
for example, the functions $C(\varphi)$ and $S(\varphi).$ 

\smallskip
{\bf Lemma 2.2} \ {\it 
	Assume that $C\!S_{kl}(\varphi)=(C(\varphi),S(\varphi))$ is the solution of the system \eqref{cs} with initial values $0,b_{kl},l,$ 
	described in \eqref{nd}, then, for \,$\sigma$ from \eqref{sv} }
\begin{equation}\label{Mcs} \exists\,\varrho>0\!:\ \ M=\max_{\,\varphi\in B_{\varrho}}\{|C(\varphi)|,|S(\varphi)|\}<\sigma. \end{equation}

{\it Proof.}\, 
Taking into account remark 2.1, we have evaluations \,$M_s^r=\max\limits_{\varphi\in \mathbb{R}}|S(\varphi)|=u_0^{1e}<\sigma$ 
and \,$M_c^r=\max\limits_{\varphi\in \mathbb{R}}|C(\varphi)|\le M_s^r$ for the solution $C\!S_{00}(\varphi)$ with $b_{00}=b_{00}^e\in (r^e,r^\sigma),$ which parametrizes the arbitrary cycles of the class $0^e].$ 

Let us choose the constant $\delta$ from Lemma\;2.1 in such way, that inequality $M_s^r+\delta<\sigma$ is correct.  
Then, such constant $\varrho_s>0$ exists, that the function $S(\varphi)$ is analytic on
$B_{\varrho_s}$\, denoted in \eqref{BB}, $M_s=\max\limits_{\,\varphi\in B_{\varrho_s}}|S(\varphi)|<\sigma,$ 
and such constant $\varrho_c>0$ exists, that $\max\limits_{\,\varphi\in B_{\varrho_c}}|C(\varphi)|\le M_c<\sigma.$ 

The extremal values of the cycles from the other classes are lesser or equal to $u_0^{1e},$ therefore such constant $\varrho$ exists, that for the parametrization of the arbitrary cycle from these classes $|C\!S_{kl}(\varphi)|\le M_s.\ \ \Box$

\section{The passing into a neighbourhood of a closed orbit} 

\subsection{Monotonicity indicator of the angular variable.} \
For each $b_{kl}$ from \eqref{nd} we introduce the real-analytic $\omega_{kl}$-periodic function 
\begin{equation}\label{p} \alpha_{kl}(\varphi)=C'(\varphi)(S(\varphi)-l)-(C(\varphi)-k)S'(\varphi). \end{equation}
Using formula \eqref{oi} the  function can be written in a simpler form as 
$$\alpha_{kl}(\varphi)=a_{kl}-1-\gamma+C^2(\varphi)+\gamma S^2(\varphi)-k(C^3(\varphi)-C(\varphi))-\gamma l(S^3(\varphi)-S(\varphi)).$$

Differentiating this equality with the respect to the system \eqref{cs}, we have:
\begin{equation}\label{ader} \alpha'_{kl}(\varphi) = \gamma((S^3 - S)(2C - 3kC^2 + k) - (C^3 - C)(2S - 3lS^2 +l)). \end{equation}

Let us study how the function $\alpha_{kl}(\varphi)$ changes its sign along the orbits which go counter-clockwise 
when  $\varphi$ increases for the class $0^i]$, and go clockwise when $\varphi$ increases for the other classes.

The passage in a neighborhood of an arbitrary cycle is possible only if $\alpha_{kl}$ is a function of a fixed sign, 
because geometrically the sign of $\alpha_{kl}$ reflects the  monotonicity of the change of the angular variable, 
assuming we observe the movement along the cycle from the point $(k,l).$

This is the reason why the same function $\alpha_{00}$ cannot be used for the passing to cycles from the classes 1] and 2].
Assuming we observe the movement from the origin,  the  polar angle is not changing monotonically.
The subtraction of constants $k$ and $l$ in formula \eqref{p} means the  origin shifts to the point $(k,l)$ in system \eqref{cs}.

\subsection{Monotonicity  of the angular variable in class 0].} \
Let us show that $\alpha_{00}=a_{00}-1-\gamma+C^2(\varphi)+\gamma S^2(\varphi)$ is a function of a  fixed sign.

Substituting $k=l=0$ into \eqref{ader}, we have: \,$\alpha_{00}'=2\gamma C(\varphi)S(\varphi)(S^2(\varphi)-C^2(\varphi)).$

For the class $0^i]\ (a=a_{00}^i\in (1,1+\gamma)),$ the cycles are located in the first quadrant when $\varphi\in[0,\omega/4]$ 
and pass through the points $(r_0^{0i},0)$ when $\varphi=0$ and $(0,u_0^{0i})$ when $\varphi=\omega/4.$ 
Constants $r_0^{0i}=b_{00}^i$ and $u_0^{0i}$ are given in \eqref{exc}.

Note that  $\alpha_{00}'(\varphi)<0$ when $C_0(\varphi)>S_0(\varphi)$ and $\alpha_{00}'(\varphi)>0$ when $C_0(\varphi)<S_0(\varphi).$ 
Therefore, the function $\alpha_{00}(\varphi)$ takes maximum values at the endpoints $[0,\omega/4].$
However, $\alpha_{00}(0)=a-1-\gamma+(r_0^{0i})^2=(a-\gamma)^{1/2}((a-\gamma)^{1/2}-1)<0,$
$\alpha_{00}(\omega/4)=a-1-\gamma+\gamma (u_0^{0i})^2=(a-1)^{1/2}((a-1)^{1/2}-\sqrt\gamma)<0.$
 
Thus, due to the symmetry, we conclude that  $\alpha_{00}(\varphi)<0$ for any $\varphi.$

For the class $0^e]\ (a=a_{00}^e>1),$  the cycles are located in the first quadrant when $\varphi\in[3\omega/4,\omega]$ 
and pass through the point $(0,u_0^{0e})$ when $\varphi=3\omega/4$ and the point  $(r_0^{0e},0)$ when $\varphi=\omega.$ 
The constants $r_0^{0e}=b_{00}^e$ and $u_0^{0e}$ are given in \eqref{exc}.

We observe also that  $\alpha_{00}'(\varphi)>0$ when $C_0(\varphi)<S_0(\varphi)$ and $\alpha_{00}'(\varphi)<0$ when $C_0(\varphi)>S_0(\varphi).$ 
Therefore, the function $\alpha_{00}(\varphi)$ takes minimum values at the endpoints $[3\omega/4,\omega].$

However, $\alpha_{00}(3\omega/4)=a-1-\gamma+\gamma (u_0^{0e})^2>0$ and $\alpha_{00}(0)=a-1-\gamma+(r_0^{0e})^2>0.$ 
Thus, for each $\varphi,$ we conclude that  $\alpha_{00}(\varphi)>0$ due to the symmetry.

\subsection{Monotonicity of the angular variable in class 1].}
\ Assume $k=1,\,l=0$ in \eqref{p}, 
then $\alpha_{10}(\varphi)=a_{10}-\gamma+\gamma S^2-C^3+C^2+C-1,$ where $C\in (0,\sqrt{2}).$  
According to \eqref{oi} $S^2=1\pm \gamma^{-1/2}(a_{10}-(C^2-1)^2)^{1/2},$ 
therefore $\alpha_{10}^\pm=a_{10}-(C-1)^2(C-1)\pm \gamma^{-1/2}(a_{10}-(C^2-1)^2)^{1/2}$ and   $a_{10}-(C^2-1)(C-1))>a_{10}-(C^2-1)^2\ge 0.$
This implies that $\alpha_{10}^+>0$ $(S\ge 1).$  

Assume $S\in [0,1].$ According to \eqref{ndb}, $a_{10}=\gamma +(b_{10}^2-1)^2,$ therefore
$$\alpha_{10}^-=\alpha_{10}^-(b_{10},C,\gamma)=\theta+\gamma-(C-1)^2(C-1)-\gamma^{-1/2}(\theta+\gamma-(C^2-1)^2)^{1/2},$$ 
where $\theta=(b_{10}^2-1)^2,$\, $b_{10}\in (1,r^e)$ due to \eqref{nd}, \,$r^e=\sqrt{1+(1-\gamma)^{1/2}}.$ Then $\theta\in (0,1-\gamma).$ 

The equation $\alpha_{10}^-=0$ is equivalent to the following quadratic equation 
\begin{equation}\label{alm} 	\theta^2-2((C-1)^2(C+1)-\gamma/2)\theta+(C-1)^3 (C+1) (C^2-1+\gamma))=0. \end{equation}
Substituting its roots in the equality $b_{10}=(\sqrt \theta+1)^{1/2},$ two functions, \,$$b_{10}^{\mp}(C,\gamma)=\big(1+\big((C^2-1)(C-1)-\gamma/2\mp (\gamma/4)^{1/2}(\gamma-\tilde \gamma(C))^{1/2}\big)^{1/2}\big)^{1/2},$$ 
where the function $\tilde \gamma(C)=4C(C^2-1)(C-1)$ is such, that \,$\alpha_{10}(b_{10}^{\mp}(C,\gamma),C,\gamma)\equiv 0.$ 

Set
$$C_*=(\sqrt{17}-1)/8,\ \ \gamma_*=(51\sqrt{17}-107)/128,\ \ b_*=(1+(297-65\sqrt{17})^{1/2}/16)^{1/2}.$$

Using the equality $\tilde \gamma'(C)=4(C-1)(4C^2+C-1),$ we obtain: 
$$\begin{matrix} \gamma_*=\tilde \gamma(C_*)=\max\limits_{\,C\in [0,\sqrt{2}]}{\tilde \gamma(C)}\approx 0.807\quad (C_*\approx 0.39),\hfill\\
	b_*=b_{10}^{\mp}(C_*,\gamma_*)=(1+((C_*^2-1)(C_*-1)-\gamma_*/2)^{1/2})^{1/2}\approx 1.156.\end{matrix} $$ 

\includegraphics[scale=0.28]{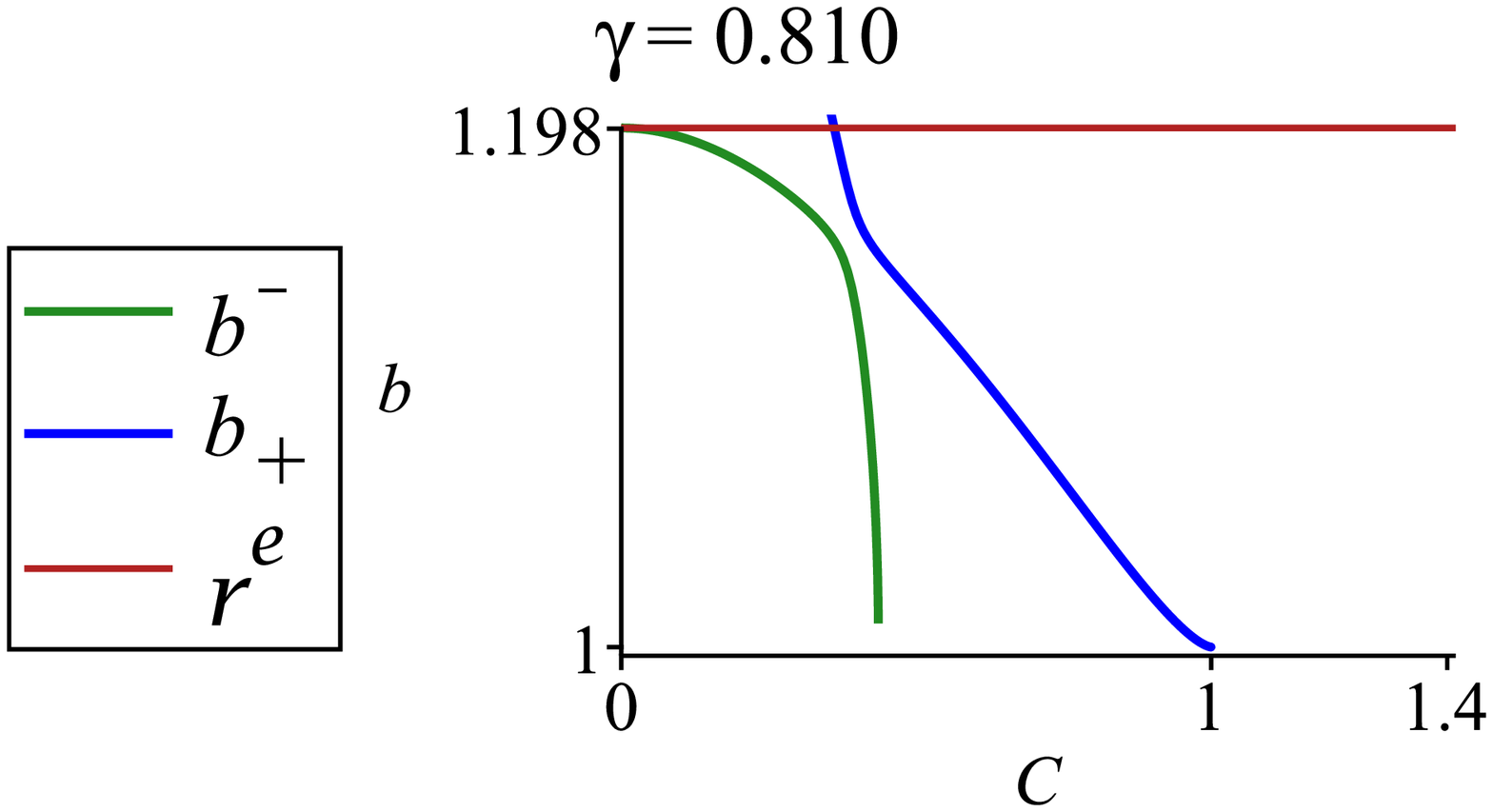} 
\includegraphics[scale=0.25]{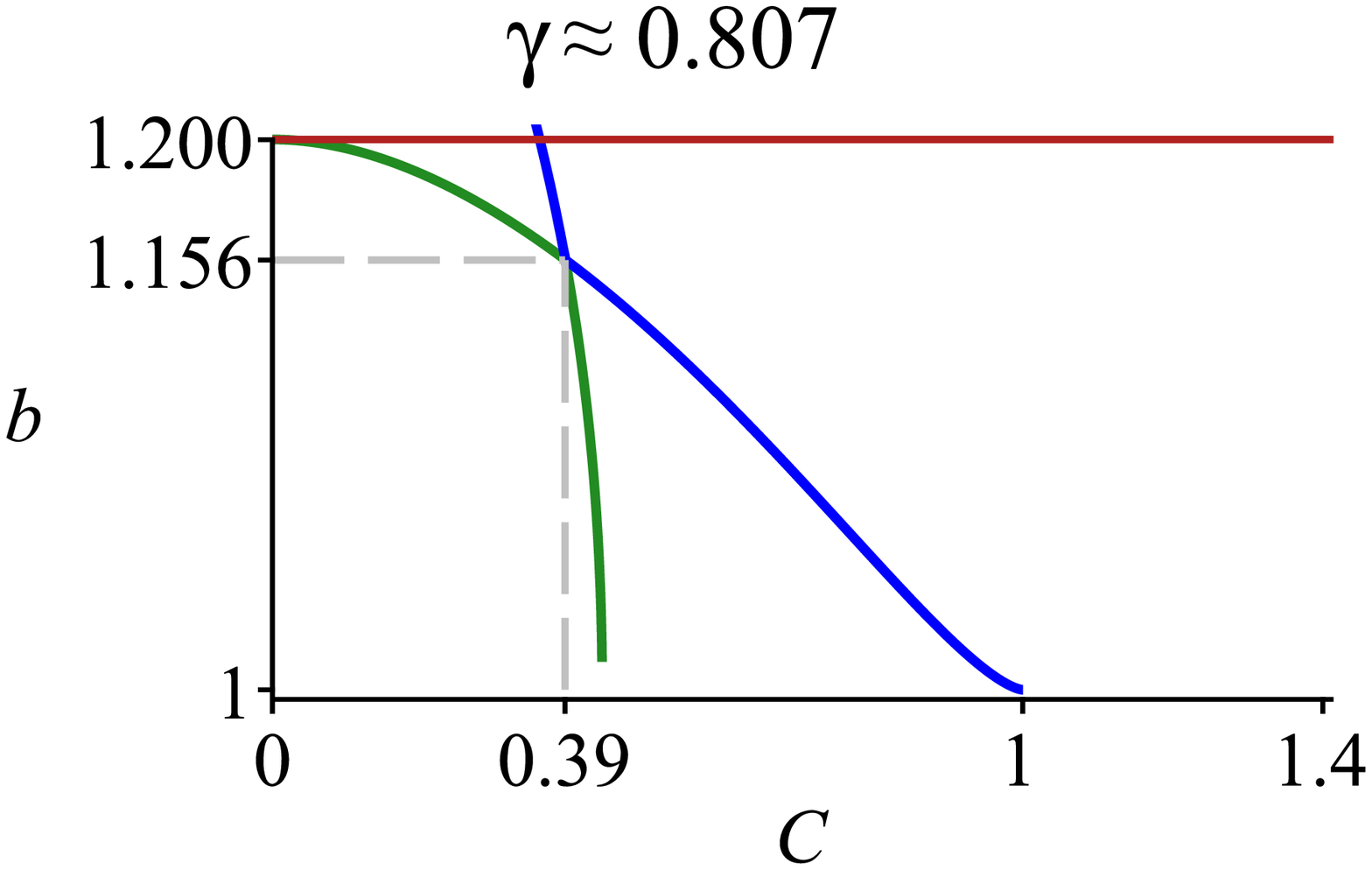} 
\includegraphics[scale=0.25]{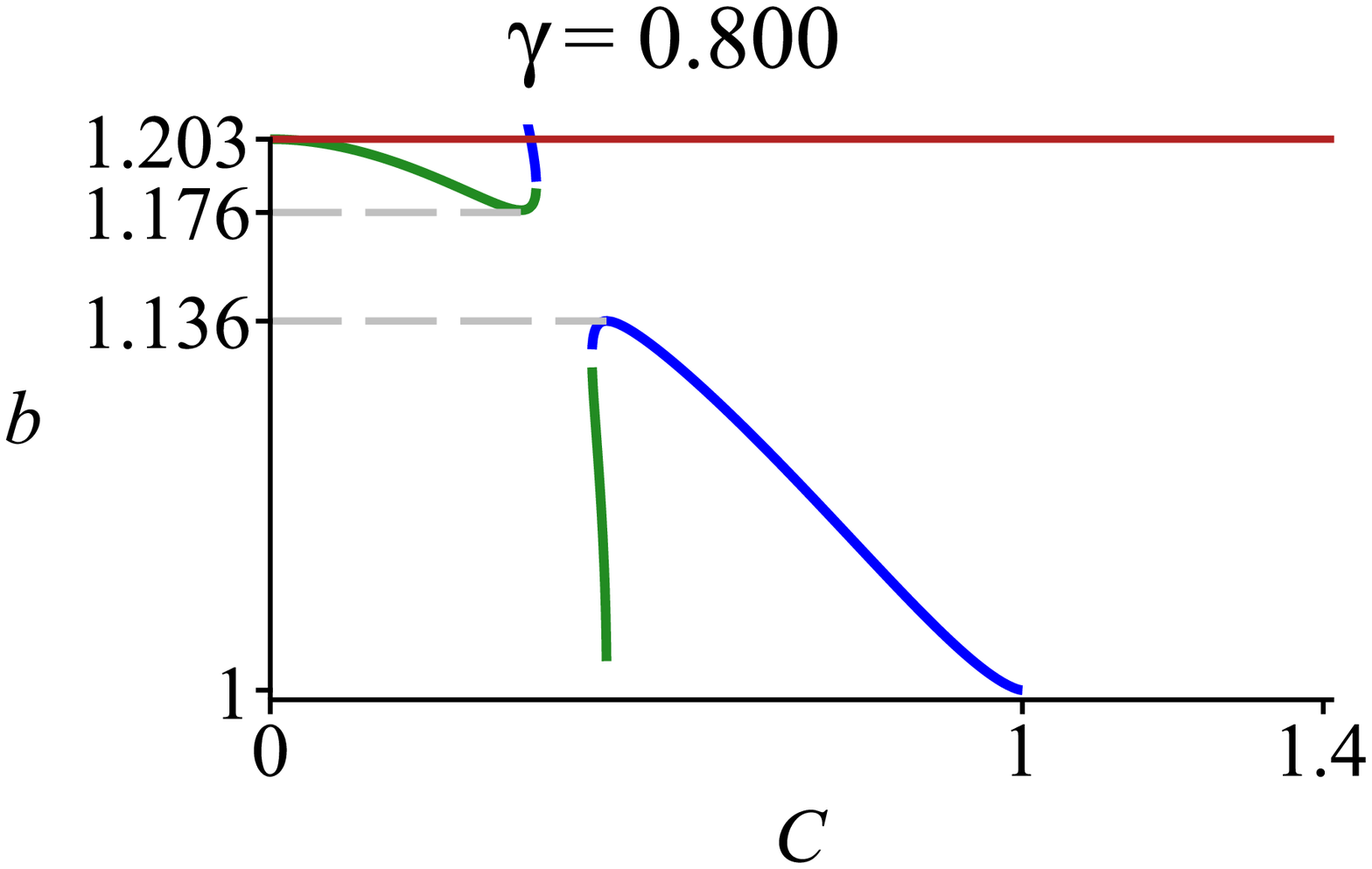} 

{\small \bf Fig.\,3.1. The curves of the arguments $b_{10},C,$ such that $\alpha_{10}^-(b_{10},C,\gamma)=0$ } 

\smallskip
Thus, for any parameters $\gamma\in [\gamma_*,1)$ the equation \eqref{alm} has at least one solution, 
because its discriminant $\gamma(\gamma - \tilde\gamma(C))\ge0,$ hence, for any $b_{10}\in(1,r^e),$ where $r^e=r^e(\gamma)=\sqrt{1+(1-\gamma)^{1/2}},$ 
the function $\alpha_{10}^-(b_{10},C,\gamma)$ has zero values (see fig.\,$3.1_1,\,3.1_2).$

Consider $\gamma\in (0,\gamma_*).$ Then, the equation $\tilde\gamma(C)=\gamma$ 
has two such solutions $C_1^\gamma < C_2^\gamma,$ that $(C_1^\gamma,b_{10}^-(C_1^\gamma,\gamma))$ and $(C_2^\gamma,b_{10}^+(C_2^\gamma,\gamma))$ 
are the contact points of the graphs of the functions $b_{10}^-(C,\gamma)$ and $b_{10}^+(C,\gamma)$ (see fig.\,$3.1_3,\,3.2$).

\includegraphics[scale=0.30]{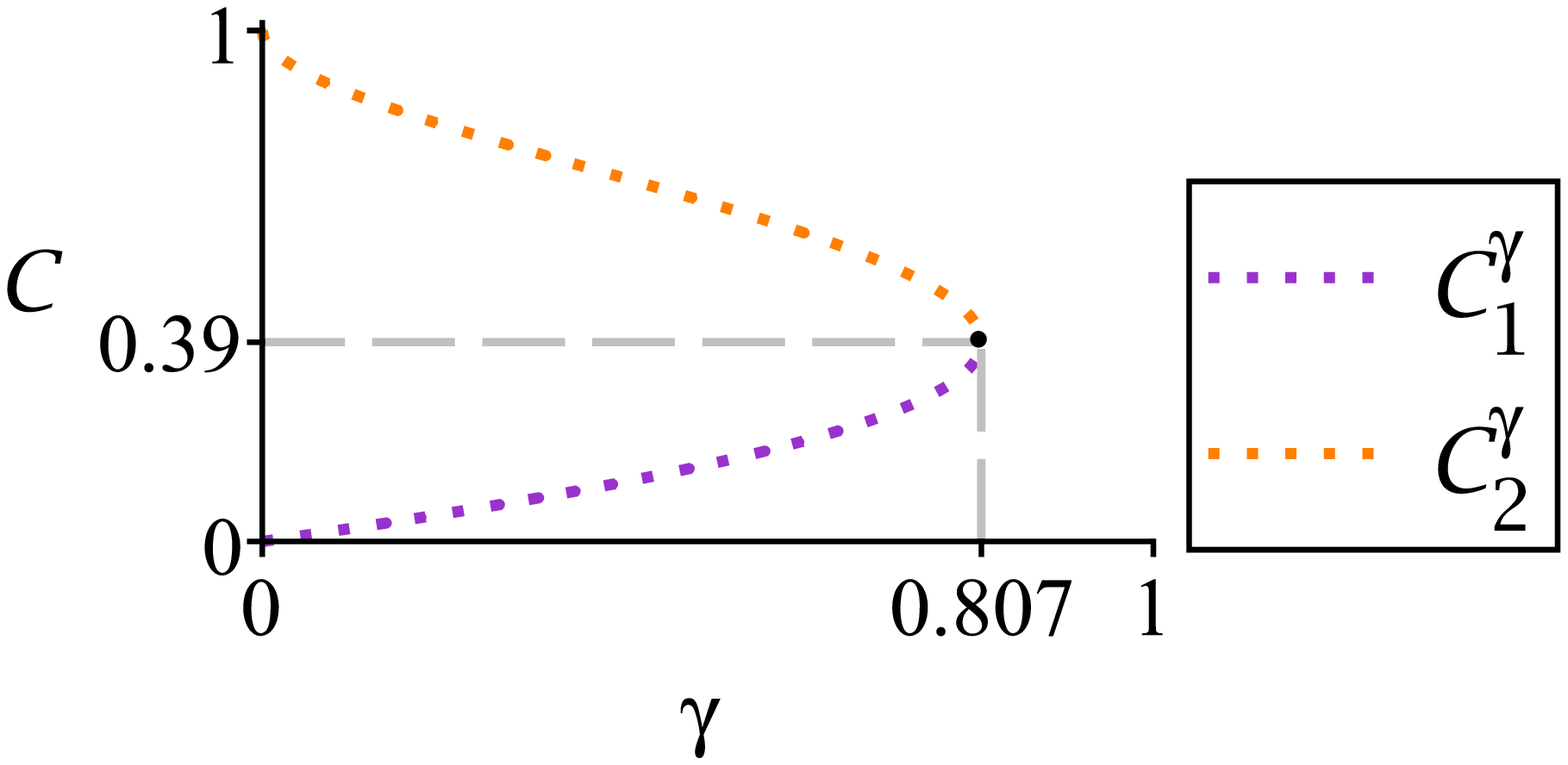} \qquad\qquad\quad 
\includegraphics[scale=0.39]{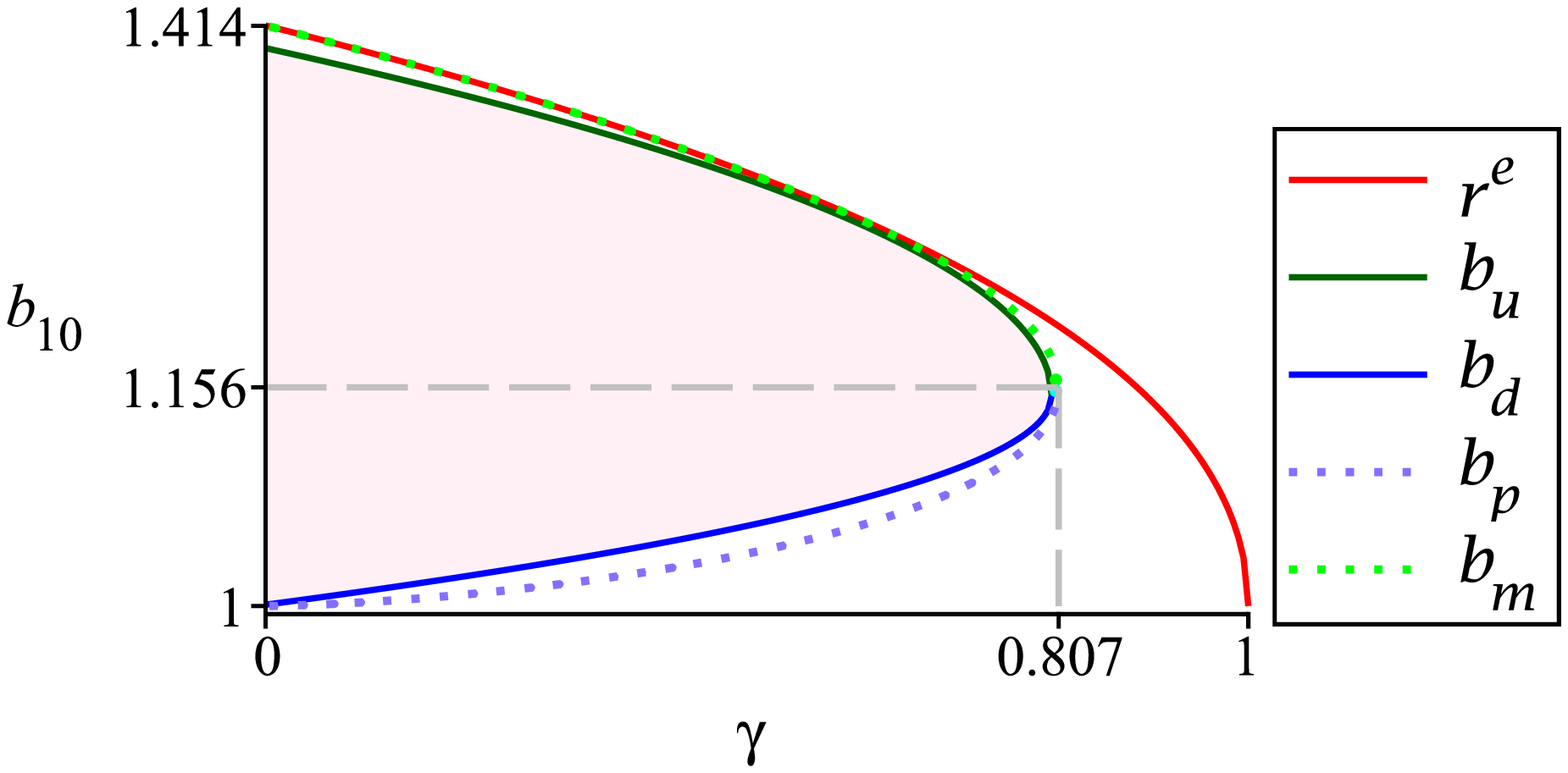}

{\small \bf 
	Fig.\,3.2. The curves of the contact $\phantom {aaaaaaaaa}$  Fig.\,3.3. The area of the positivity \\
	$\phantom {aaa}$ points of the functions $b_{10}^{\mp}(C,\gamma)$ $\phantom {aaaaaaaaaaaaaaaa}$ of the function $\alpha(b_{10},C,\gamma)$ } 

\smallskip 

Introduce the functions \,$b_{10}=b_m(\gamma)$ and \,$b_{10}=b_p(\gamma)\!:$ \,$b_m(\gamma_*)=b_p(\gamma_*)=b_{\mp}(C_*,\gamma_*),$ 
\begin{equation}\label{bpbm}
b_m=\min\limits_{C\in (0,C_1^\gamma]} b_-(C,\gamma),\ \ 
b_p=\max\limits_{C\in [C_2^\gamma,1)} b_+(C,\gamma)\ \hbox{ with }\ \gamma\in (0,\gamma_*).
\end{equation}
Their approximate values are found, using the computing environment MAPLE with the discrete change of the parameter \,$\gamma$\, equal to $10^{-3}.$

Particularly, $1<b_p(\gamma)<b_m(\gamma)<r^e(\gamma)$ with $\gamma\in (0,\gamma_*).$
Figure $3.1_1$ has constants $b_p\approx 1.136,\ b_m\approx 1.176$ with $\gamma=0.8<\gamma_*.$

Summing it up, assume $0<\gamma<\gamma_*,$ $b_{10}\in (b_p(\gamma),b_m(\gamma)),$ $(C(\varphi),S(\varphi))$ --- $\omega_{10}$-periodic solution of system's \eqref{cs} initial value problem with initial values $C(0)=b_{10},\,S(0)=0.$ This solution parametrizes specific cycle $\Gamma.$ 
Notice, that the part of $\Gamma,$ located in the first quadrant, is parametrized, when $\varphi\in [\omega_{10}/2,\omega_{10}]$
(see proof of lemma\;1).
It was proved, that for such $\varphi,$ function $\alpha_{10}(\varphi),$ set on the solution $(C(\varphi),S(\varphi)),$ is positive. 
If $b_{10}$ from \eqref{nd} is not from specific interval, or $\gamma\in [\gamma_*,1),$ then there are values of $\varphi,$ 
such that $\alpha_{10}(\varphi)$ changes its sign. 

\smallskip 
{\bf Lemma 3.1.} \ {\it 
	The following statements hold:\\
	$1)$ for each $\gamma\in (0,\gamma_*)$ with $\gamma_*=(51\sqrt{17}-107)/128,$ function $\alpha_{k0}(\varphi)$ 
	is positive, when $kb_{k0}\in(b_p(\gamma),b_m(\gamma))$ $(k=\pm 1),$\ $b_p,b_m$ are from \eqref{bpbm},\ $k=\pm 1;$\, 
	is equal to zero, when $kb_{k0}=b_p(\gamma)$ or $kb_{k0}=b_m(\gamma),$\,
	changes its sign, when  $kb_{k0}\in(1,b_p(\gamma))\cup (b_m(\gamma),r^e(\gamma));$\, 
	$2)$~for any $\gamma\in [\gamma_*,1),$\, the function \,$\alpha_{k0}(\varphi)$\, changes its sign for all $kb_{k0}\in(1,r^e(\gamma)).$ }

P r o o f. \ The result is already obtained for $k=1$ and $S\ge 0,$ 
that is for the part of the cycles from class 1], located in the first quadrant. 
It automatically expands on the case $S\le 0$ and left semiplane, when $k=-1,$ 
because integral \eqref{oi} doesn't change after substituting $C$ to $-C$ and $S$ to $-S.$ \ $\Box$

Because values of $b_p(\gamma)$ and $b_m(\gamma)$ are not given explicilty, 
we will use these limitations:
\begin{equation} \label{gamma*} 0<\gamma\le 0.8,\quad kb_{k0}\in (b_d(\gamma),b_u(\gamma))\ \ (k=\pm 1)\quad  (d - down,\ u - up),\end{equation}
where \,$b_u=1.15+0.28(0.8-\gamma)^{1/2},\ b_d=1.15-(0.8-\gamma)^{1/2}.$ 
It's also obtained for any $\gamma\in (0,0.8],$ that $b_{10}=b_d(\gamma),\ b_{10}=b_u(\gamma).$

If it's necessary, the conditions \eqref{gamma*} can be specified by approximating functions $\gamma,b_d,b_u$ to functions $\gamma_*,b_p,b_m.$

\subsection{Monotonicity of the angular variable in the class 2].} 
Let us show that $\alpha_{kl}(\varphi)$ $(k,l=\pm1)$ from $(\ref{p})$ is positive.
Because $k^2 =l^2 =1,$ formula \eqref{ader} can ber written as 
$$\alpha_{kl}'(\varphi)=\gamma (C(\varphi)-k) (S(\varphi)-l)( kC(\varphi)- lS(\varphi))(3kl C(\varphi) S(\varphi)+k C(\varphi)+ l S(\varphi) +1).$$

Let $lk=1.$ Since  $lS>0>-( C+k)/(3 C+k),$ by multiplying the  inequality  above by $k(3C+k) ( >0)$ we get that $3 k l C S+k C+l S+1>0$. 

Considering that the  movement along the orbit is clockwise the analysis of  the sign of  $\alpha_{kl}$'s shows
that the function $\alpha_{kl}(\varphi)$ has  local minimums   at the points $\varphi_1,\varphi_2,\varphi_3$ determined by the equations
$ C(\varphi_1)=k,$ $S(\varphi_1)=l\cdot  lo_2^1;$\, $ S(\varphi_2)=l,$ $C(\varphi_2)=k l_2^1;$\, $kC(\varphi_3)= lS(\varphi_3)>1.$ 

According to formulas $(\ref{exc})$ $l_2^1=\sqrt{1-a^{1/2}},\ \ lo_2^1=\sqrt{1-(a/\gamma)^{1/2}},$ hence,
$a=\gamma(1-(lo^1_2)^2)^2=(1-(l^1_2)^2)^2,$ then,
$\alpha_{kl}(\varphi_1)=a+\gamma((lo^1_2)^2 - 1)(1-lo^1_2 )=\gamma lo^1_2(lo^1_2 - 1)^2(lo^1_2+1)>0.$
Similarly, $\!\alpha_{kl}(\varphi_2)=a+((l^1_2)^2 - 1)(1-l^1_2 )=l^1_2(l^1_2 - 1)^2(l^1_2+1)>0.$ 
Finally, $ \alpha_{kl}(\varphi_3)=(1+\gamma)(S^2 - 1)S(S-l)>0.$ 

The property that  $\alpha_{kl}(\varphi)>0$ when $kl=-1$ can be proved analogously.

Thus, generalizing the results from the section 3.2-3.4, we formulate

\smallskip 
{\bf Lemma 3.2.} \ {\it
	The real-analytic $\omega_{kl}$-periodic function $\alpha_{kl}(\varphi),$ denoted in \eqref{p} and set on the solutions $CS_{kl}(\varphi)$ of the system \eqref{cs}, which parametrize the cycles of the system \eqref{snv}, is positive for all classes of the cycles, except the cycles from the class $0^i],$ where it is negative. Moreover, there is an additional limitation \eqref{gamma*} that affects the positivity of the functions $\alpha_{10}(\varphi)$ and $\alpha_{-1,0}(\varphi)$ from the class $1].$}

\smallskip 
\includegraphics[scale=0.25]{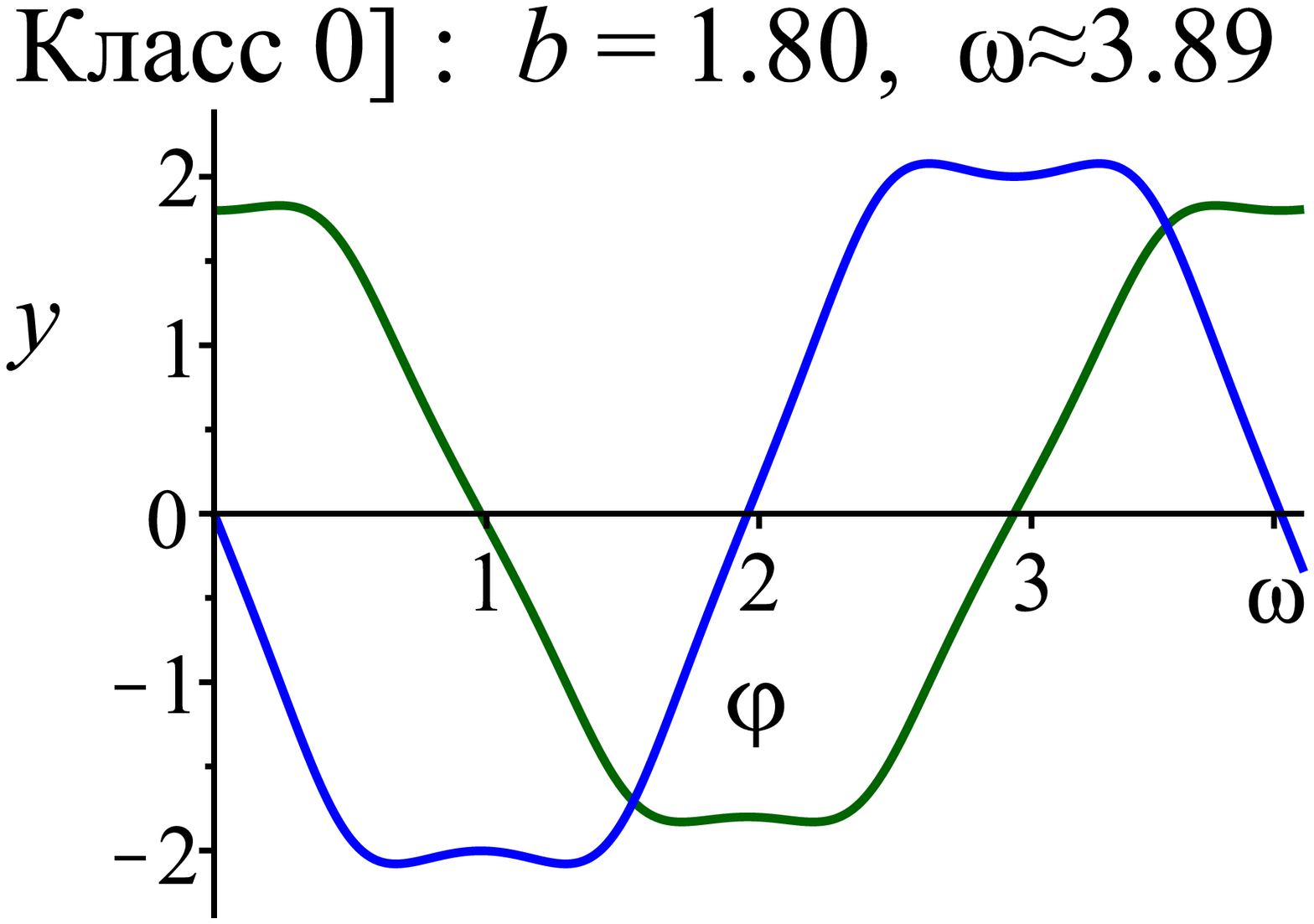} \includegraphics[scale=0.25]{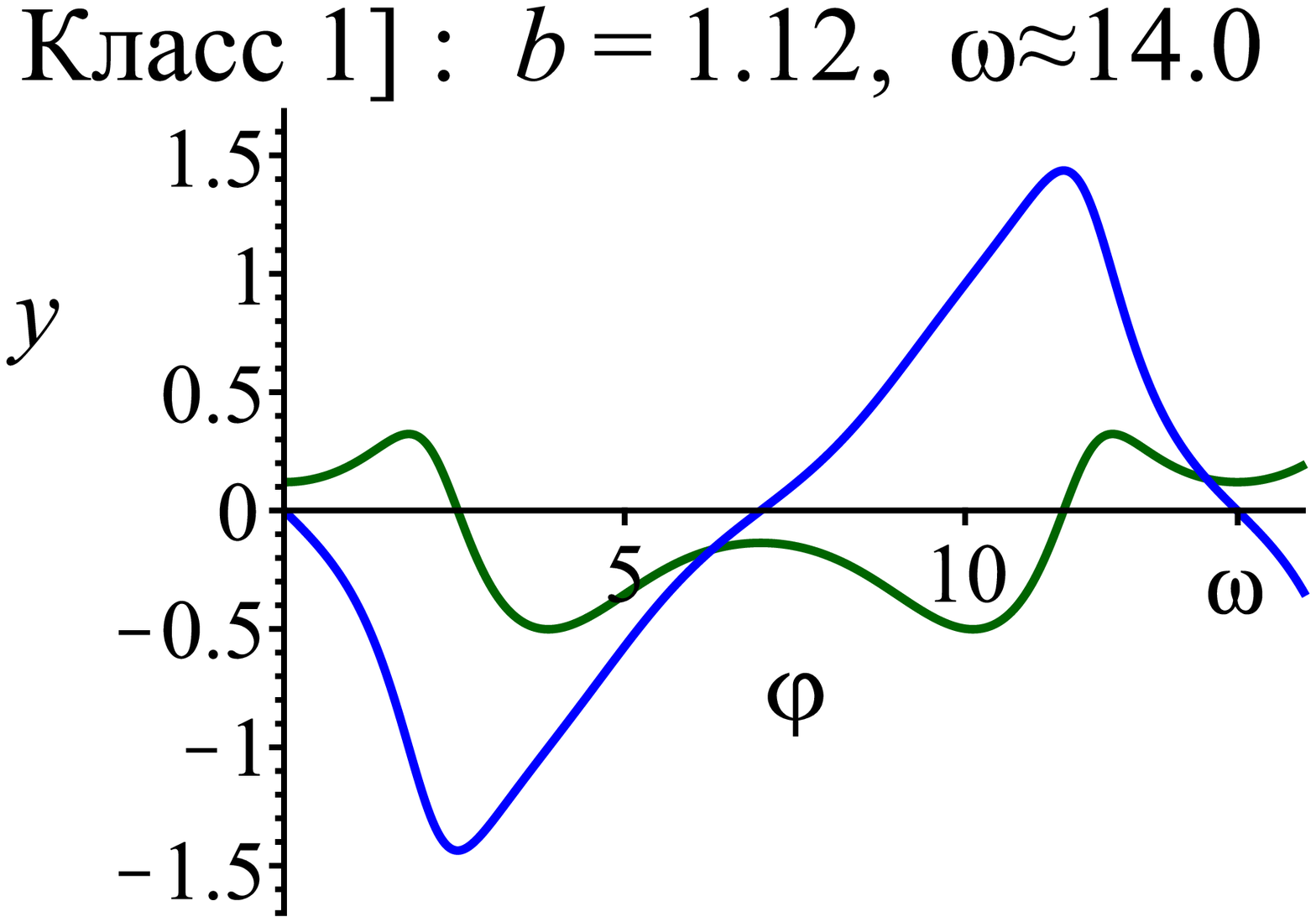} \includegraphics[scale=0.29]{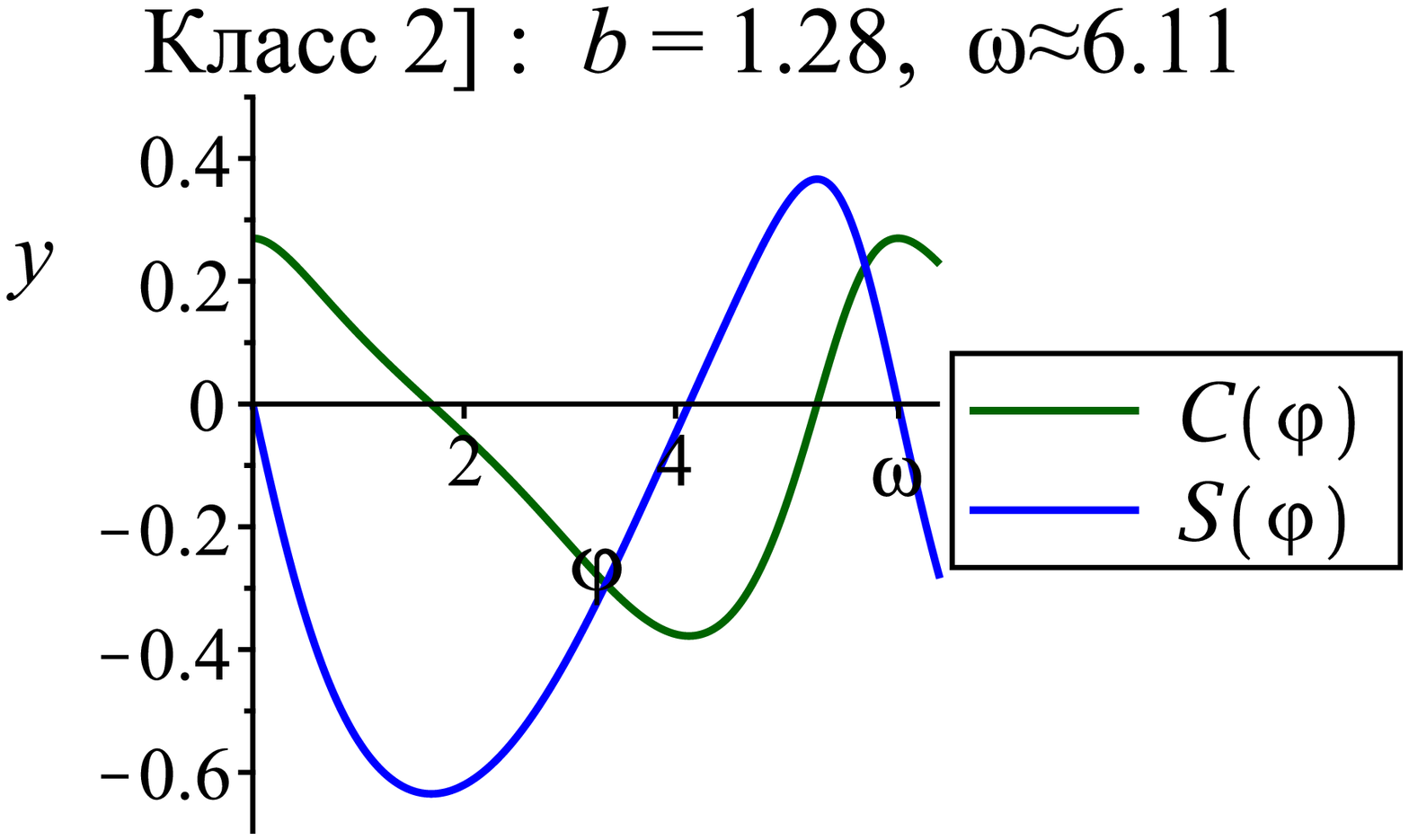} 

\includegraphics[scale=0.25]{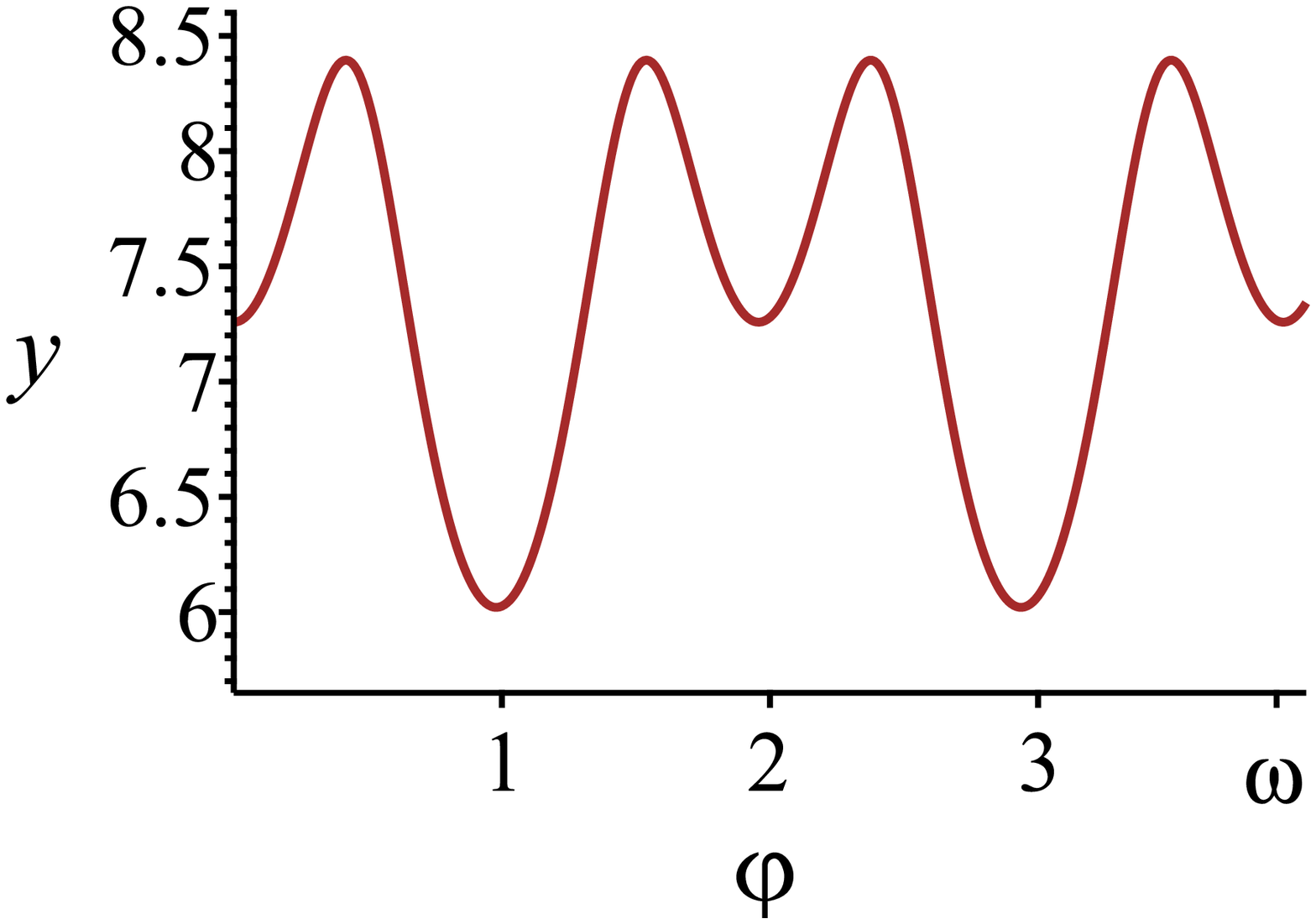} \includegraphics[scale=0.25]{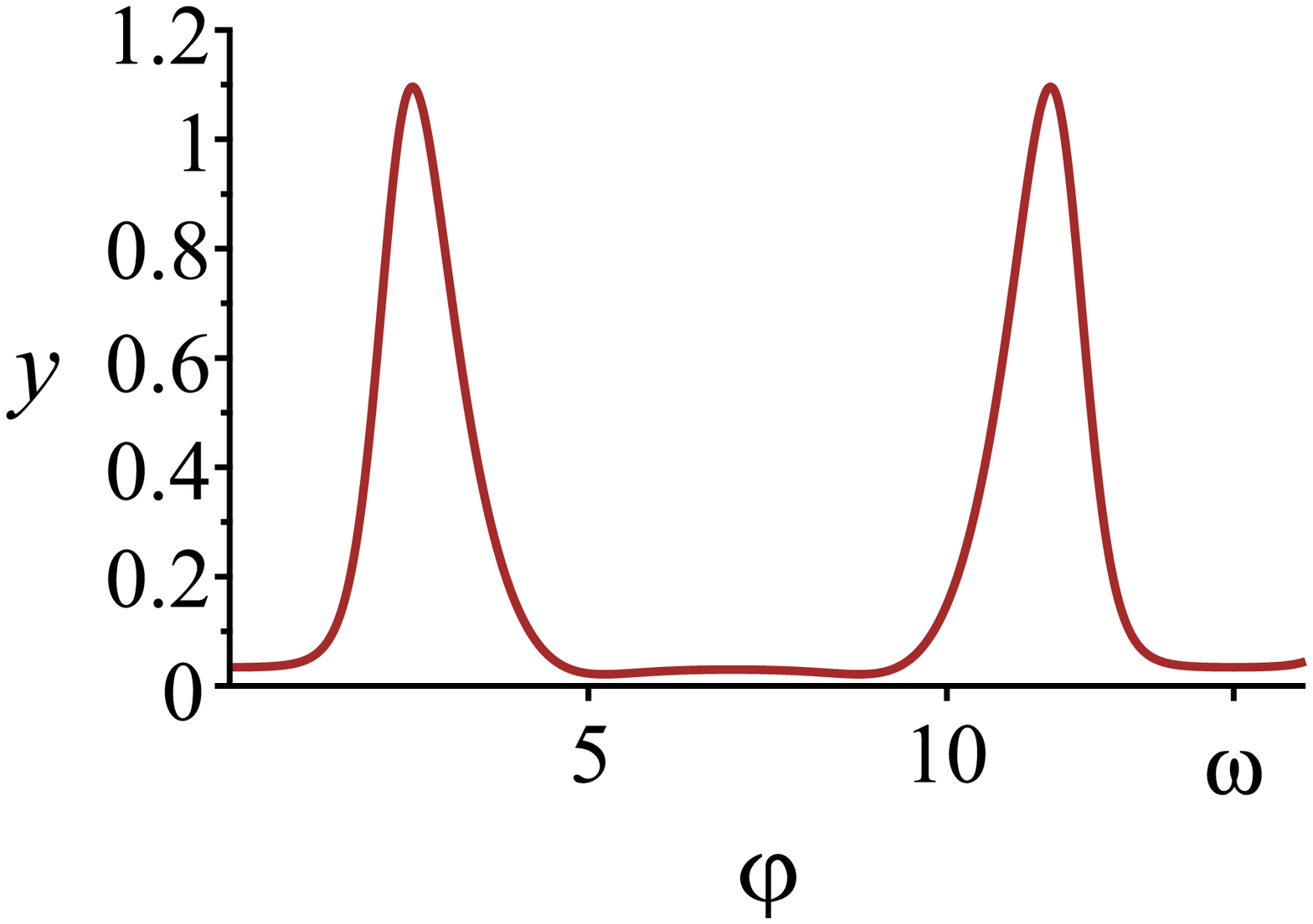} \includegraphics[scale=0.28]{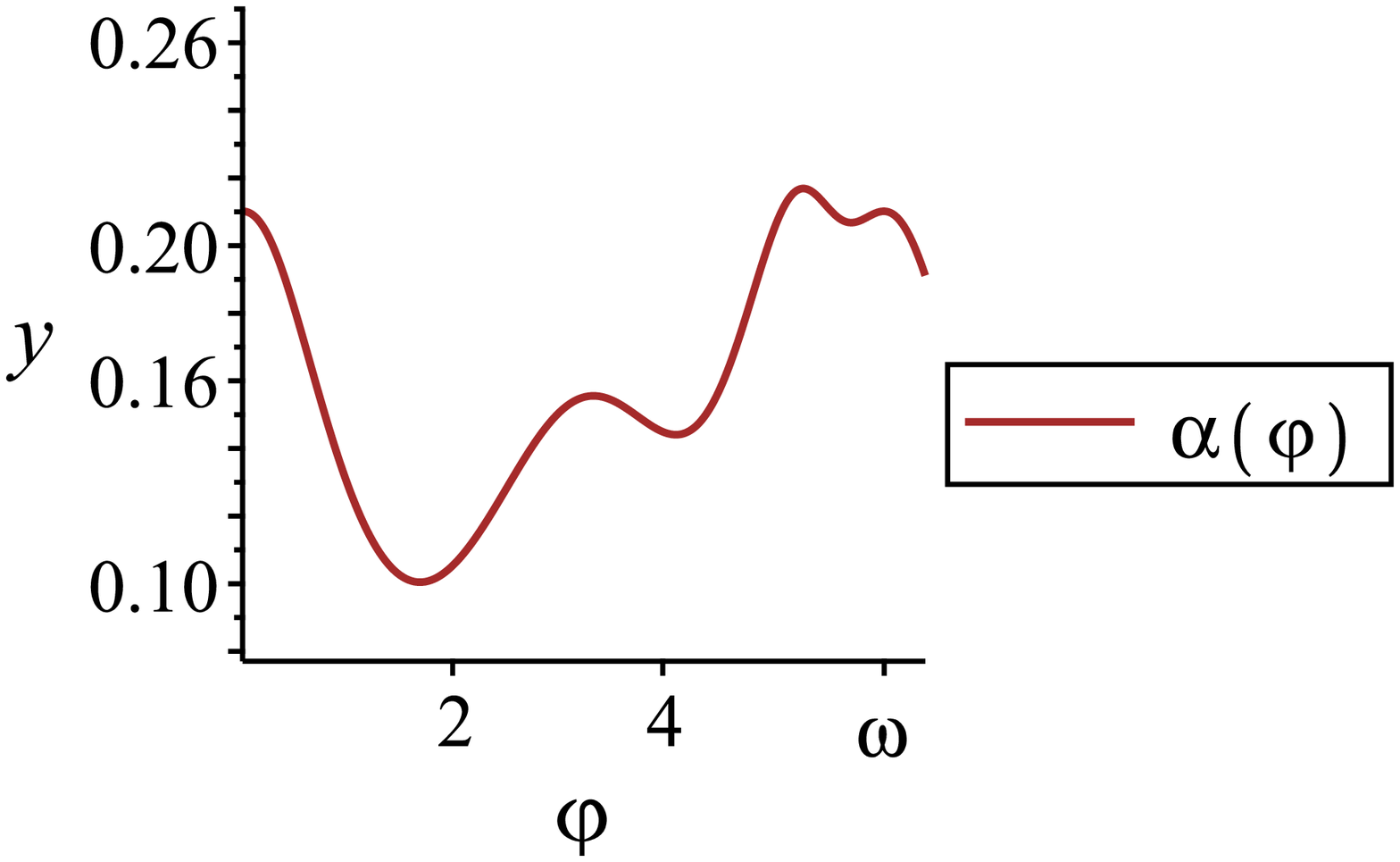} 

\begin{center} {\small \bf Fig.\,3.4. \ The graphs of the functions $C,S,\alpha$ with $\gamma=1/2.$  } \end{center}

\subsection{The passing to the special polar system.} \  
Let us consider the special polar change in system \eqref{sv}:
\begin{equation}\label{pz} 
	x=C(\varphi)+(C(\varphi)-k)\rho,\ \ y=S(\varphi)+(S(\varphi)-l)\rho\quad (|\rho|\le \rho_*<1),
\end{equation}
where $k,l=0,\pm 1,\ (k,l)\ne (0,\pm 1);$
vector-function ${C\!S}_{kl}(\varphi)=(C(\varphi),\,S(\varphi))$ is real-analytic $\omega_{kl}$-periodic solution of system \eqref{cs} with initial values $C(0)= b_{kl},\,S(0)=l,$ bounds of $b_{kl}$ are provided in \eqref{nd} and specified for $b_{10},b_{-1,0}$ in \eqref{gamma*}.

The change \eqref{pz} is the composition of the three following changes.

1) Change $x=x_k+k,\ y=y_l+l$ 
performs the shift of the origin into the equilibrium points $(k,l)$ of unperturbed system \eqref{snv}, 
which allows to track the motion along the cycles of the classes 1] and 2], while being inside of them as it happens for the class 0], when $k=l=0.$ 
Let us notice that the similar shift $C=C_k+k,$ $S=S_l+l$ in system \eqref{cs} transforms it into the system 
$C_k'(\varphi)=\gamma(S_l^3(\varphi)-S_l(\varphi)),\ S_l'(\varphi)=-(C_k^3(\varphi)-C_k(\varphi)).$

2) Non normalized polar change $x_k=C_k(\varphi)\tilde\rho,\ y_l=S_l(\varphi)\tilde\rho\ \ (\tilde\rho>0)$ 
is performed in the neighbourhood of the point $(k,l)$ and with $\varrho=1$ determines a cycle, 
parametrized by the chosen solution $(C(\varphi),S(\varphi))$ of the unperturbed system. 

Let us notice that, in fact, the non normalized analogues of the generalized sine and cosine are used in this change. These functions were introduced by the A.\,M. Lyapunov as a solutions of system ${\rm Cs}'\theta=-d\,{\rm Sn}\,\theta,$ ${\rm Sn}'\theta={\rm Cs}^{2n-1}\,\theta$ 
with initial values ${\rm Cs}\,0=1,$ ${\rm Sn}\,0=0.$ 
Similar normalization, for example, was used in \cite{Bib2}. 
Refusal to do such normalization allows to obtain bifurcation (generating) equation of another, integral nature.

3) Affine change $\tilde\rho=\rho+1$ performs the shift into the neighbourhood of the cycle, which passes through the point $(b_{kl},l)$ and is obtained by assuming $\rho=0.$ 

\smallskip
{\bf Proposition 3.1.} \ {\it  	
	Assume that in the substitute \eqref{pz} 
	\begin{equation}\label{rho*} \rho_*\in (0,\,\min\,\{\rho_{kl},\rho_0\}), \end{equation} 
	where $\rho_{00}^i = \min \{ 1, r^i / b_{00}^i - 1\} ;$ $\rho_{00}^e = \min\{ 1 - r^e / b_{00}^e, r^\sigma / b_{00}^e - 1\} ;$ 
	$\rho_{k0} = \min\{  (b_{k0} - kb_d) / (b_{k0} - k), (kb_u - b_{k0}) / (b_{k0} - k)\}<1$ $(k\ne 0);$
	$\rho_{kl} = \min\{ 1, (kr^s - b_{kl}) / (b_{kl} - k)\} \ (k,l\ne 0);$ \ $\rho_0=(\sigma-M)/(M+1)$ with the constant $M$ from \eqref{Mcs}.\,  
	Then
	
	\,$1)$ for $|\rho|\le \rho_*<\rho_{kl},$ the cycle of the unperturbed system \eqref{snv}, passing through the point $(\rho,l),$ 
	is of the same class $|k|+|l|$ as a chosen cycle parametrized by the solution 	$C\!S_{kl}(\varphi),$  
	which can be obtained by assuming $\rho=0$ in \eqref{pz};\, 
	
	$2)$ for $|\rho|\le \rho_*<\rho_0$ and $\varphi\in B_{\varrho}$ from \eqref{Mcs}, the following is correct: 
	$|C(\varphi)+(C(\varphi)-k)\rho|<\sigma$\, and \,$|S(\varphi)+(S(\varphi)-l)\rho|<\sigma,$ therefore, 
	the substitute \eqref{pz} is possible in the right-hand side of the system \eqref{sv}. }

{\it Proof.}\,  
1) For class $0^e]$ we will check that $x(0,\rho)=b_{00}^e(1+\rho)\in(r^e,r^\sigma)$ with $|\rho|<\rho_{00}^e.$\, 
According to \eqref{rho*}, \,$r^e / b_{00}^e \le 1-\rho_{00}^e$ and $\rho_{00}^e + 1\le r^\sigma / b_{00}^e,$ 
therefore \,$r^e \le b_{00}^e(1 - \rho_{00}^e) < b_{00}^e(1+\rho) < b_{00}^e(1 + \rho_{00}^e) \le r^\sigma.$

For class $1]$ we will check that $kx(0,\rho)=kb_{k0}+k(b_{k0}-k)\rho\in(b_d, b_u)$ with $|\rho|<\rho_{k0},$ 
where $b_d,b_u$ are described in \eqref{gamma*}, $k=\pm 1.$ 

According to \eqref{rho*}, \,$(kb_d - b_{k0})/(b_{k0} - k)\le -\rho_{k0}$ and $(kb_u - b_{k0})/(b_{k0}-k)\ge \rho_{k0}.$
Then $b_d=kb_{k0}+k(kb_d - b_{k0})\le kb_{k0} - k(b_{k0}-k)\rho_{k0}< kb_{k0}+k(b_{k0}-k)\rho,$ 
$b_u=kb_{k0}+k(kb_u-b_{k0})\ge kb_{k0}+k(b_{k0}-k)\rho_{k0}>kb_{k0}+k(b_{k0}-k)\rho.$

The inequalities for the classes $0^i]$ and $2]$ can be proven analogously.

2)\ According to lemma 2.2, $M=\max_{\varphi\in B_{\varrho}}\{|C(\varphi)|,|S(\varphi)|\}<\sigma,$ therefore, for example, $|S(\varphi)+(S(\varphi)-l)\rho| < M + (M+1)\rho_0 = \sigma.$\ $\Box$

\smallskip 
{\bf Lemma 3.3.} {\it 
	Assuming the limitations \eqref{rho*} hold, the change \eqref{pz} transforms system \eqref{sv} into system
	\begin{equation}\label{ps}
		\begin{cases} 
			\dot \rho=(-\alpha_{kl}'(\varphi)\alpha_{kl}^{-1}(\varphi)\rho-p_{kl}(\varphi)\rho^2+\breve p_{kl}(\varphi)\rho^3+
			\e P_{kl}^{\nu \rho}(t,\varphi,\rho,\e))\e^\nu, \\
			\dot\varphi=(1+\alpha_{kl}(\varphi)q_{kl}(\varphi)\rho+\breve q_{kl}(\varphi)\rho^2+\e \Phi_{kl}^{\nu \rho}(t,\varphi,\rho,\e))\e^\nu, \end{cases}
	\end{equation}
	where \,$q_{kl}(\varphi)=\alpha_{kl}^{-2}( (C-k)^3 (2C+k) + \gamma (S-l)^3 (2S+l)),$  
	$p_{kl}(\varphi)=3\gamma \alpha_{kl}^{-1}CS((S^2-1)(C-k)^2-(C^2-1)(S-l)^2),$ 
	$\breve p_{kl} = \gamma\alpha_{kl}^{-1}((C^3 - C)(S - l)^3 - (S^3 - S)(C - k)^3),$ 
	$\breve q_{kl} = \alpha_{kl}^{-1}((C-k)^4 + \gamma(S-l)^4),$ and all these functions are $\omega_{kl}$-periodic and real-analytic; 
	$P_{kl}^{\nu \rho}=\alpha_{kl}^{-1}(C'Y(\ae^\rho)-S'X(\ae^\rho)),$  
	$\Phi_{kl}^{\nu \rho}=\alpha_{kl}^{-1}(\varphi)(1+\rho)^{-1}((S-l)X(\ae^\rho)-(C-k)Y(\ae^\rho));$ \ 
	$\ae^\rho=(t,C+(C-k)\rho,S+(S-l)\rho,\e).$ }

\smallskip
{\it Proof.}\, 
Differentiating change \eqref{pz} with the respect to $t$ and solving both equations with the respect to \,$\dot \rho$ and $\dot \varphi$\, 
by using $(\ref{p})$, we have:
$$\alpha_{kl}(\varphi)\dot \rho=C'(\varphi)\dot y-S'(\varphi)\dot x,\quad 
(\rho+1)\alpha_{kl}(\varphi)\dot \varphi=(S(\varphi)-l)\dot x-(C(\varphi)-k)\dot y.$$
Substituting the right-hand sides of the system $(\ref{sv}),$ we obtain 
$$\begin{matrix}\alpha_{kl}\dot \rho=(C'(-(C+(C-k)\rho)^3+C+(C-k)\rho+\e Y)- \hfill \\
	\hfill -S'(\gamma((S+(S-l)\rho)^3-(S+(S-l)\rho))+\e X))\e^\nu , \\
	(\rho+1)\alpha_{kl}\dot \varphi=((S-l)(\gamma((S+(S-l)\rho)^3-(S+(S-l)\rho))+\e X)- \\
	\hfill -(C-k)(-(C+(C-k)\rho)^3+C+(C-k)\rho+\e Y))\e^\nu .\end{matrix}$$

Simple calculations show that this system can be written as \eqref{ps} due to the fact that $|\rho|\le\rho_*<1.$ \ $\Box$

\section{Primary radial averaging} 

\smallskip
{\bf $\quad\,4.1.$ Averaging of terms of the unperturbed part of the radial equation.} \ 
Now it's necessary to average the first two terms in the right-hand side of the radial equation of the special polar system \eqref{ps}. 

Consider the polynomial with the respect to $r$ and real-analytic with the respect to $\varphi$ with $\varphi\in B_{\varrho}$ averaging change
\begin{equation}\label{pave}
	\rho= \alpha_{kl}^{-1}(\varphi)(r+ \beta_{kl}(\varphi)r^2) \quad (|r|\le r_* < 1),
\end{equation}
where the function $\displaystyle \beta_{kl}(\varphi)=\int_{0}^{\varphi} (\xi_{kl}(s)-\overline{\xi_{kl}})\,ds$\, and   
$\xi_{kl}(\varphi)=\alpha_{kl}^{-1}(\varphi)(\alpha_{kl}'(\varphi) q_{kl}(\varphi)-p_{kl}(\varphi)),$ 
$\displaystyle \overline \xi_{kl}=\frac{1}{\omega_{kl}}\int_{0}^{\omega_{kl}} \xi_{kl}(s)\,ds$ --- 
the average value of the function $\xi_{kl}(\varphi).$

\smallskip 
{\bf Proposition 4.1.} \ {\it Assume there is a limitation on the constant $r_*$ from \eqref{pave} 
	\begin{equation}\label{r*} 0<r_*<\min\,\{ (4\beta^*)^{-1},r_0 \}, \end{equation} 
	where $\beta^* =\max_{\,\varphi\in B_{\varrho}} |\beta_{kl}(\varphi)|,$ $r_0=\alpha_*(1+\beta^*)^{-1}\rho_*,$ 
	$\alpha_*=\min_{\,\varphi\in B_{\varrho}} |\alpha_{kl}(\varphi)|>0,$ $\rho_*$ from $\eqref{rho*}.$
	The following statements hold: 
	\,$1)$ $1/2<|1+2\beta_{kl}r|<3/2$ with $|r|\le r_*<(4\beta^*)^{-1};$\, 
	\,$2)~|\alpha_{kl}^{-1}(\varphi)(r+\beta_{kl}(\varphi)r^2)|<\rho_*$ with $|r|\le r_*<r_0,$ 
	therefore, the right-hand side of the system \eqref{ps} allows the substitute \eqref{pave}. }

\smallskip 
{\bf Lemma 4.1.} {\it The real-analytic $\omega_{kl}$-periodic change with respect to $\varphi$ change \eqref{pave} with $\overline{\xi_{kl}}=0$ 
	transforms the system \eqref{ps} into the system
	\begin{equation}\label{ps5}
		\begin{cases} 
			\dot r=(Q_{kl}(\varphi,r)r^3+ R_{kl}^\nu(t,\varphi,r,\e)\e)\e^\nu, \\
			\dot\varphi=(1+q_{kl}(\varphi)r+\mathcal Q_{kl}(\varphi,r)r^2+ \Phi_{kl}^\nu(t,\varphi,r,\e)\e)\e^\nu, \end{cases}
	\end{equation}
where $\mathcal Q_{kl}=\beta_{kl}q_{kl}+\alpha_{kl}^{-2}\breve q_{kl}(1+\beta_{kl}r)^2,\  
Q_{kl}=\big(\alpha_{kl}^{-1}\alpha'_{kl}\mathcal Q_{kl}-(\alpha_{kl}^{-1}\alpha'_{kl}\beta_{kl}-\beta'_{kl})(q_{kl}+\mathcal Q_{kl})
-\alpha_{kl}^{-2}(\alpha_{kl}p_{kl}(2\beta_{kl}+\beta_{kl}^2r)-\breve p_{kl}(1+\beta_{kl}r)^3)\big)(1+2\beta_{kl}r)^{-1};$ 
$\Phi_{kl}^\nu=\Phi_{kl}^{\nu\rho}(\ae^r),$ where $\ae^r=(t,\varphi,\alpha_{kl}^{-1}(r+\beta_{kl}r^2),\e);$ \ 
$R_{kl}^\nu=\big(P_{kl}^\nu 
+(\alpha^{-1}_{kl}\alpha'_{kl}r+(\alpha^{-1}_{kl}\alpha'_{kl}\beta_{kl}-\beta'_{kl})r^2)\Phi_{kl}^\nu\big)(1+2\beta_{kl}r)^{-1},$  
where $P_{kl}^\nu(t,\varphi,r,\e)=P_{kl}^{\nu\rho}(\ae^r),$ 
under assumption, that $r_*$ is such that the inequation \eqref{r*} and inequation 
\begin{equation}\label{r*2} (2r_*)^{-1}>\max\limits_{|r|\le 1,\,\varphi\in B_{\varrho}} |q_{kl}(\varphi)+\mathcal Q_{kl}(\varphi,r)r|. \end{equation}
hold. }

\smallskip
{\it Proof.}\, Substituting the change \eqref{pave} into the second equation of the system \eqref{ps}, 
we obtain the second equation of the system \eqref{ps5}.  

Let us differentiate the change \eqref{pave} with the respect to $t$ and system \eqref{ps} and, then, multiply the resulting equality by $\alpha_{kl}.$ We have 
$$\begin{matrix} 
	(1+2\beta_{kl}r)\dot r-(\alpha_{kl}^{-1}\alpha'_{kl}r+
	(\alpha_{kl}^{-1}\alpha'_{kl}\beta_{kl}-\beta'_{kl})r^2)(1+q_{kl}r+\mathcal Q_{kl}r^2+\Phi_{kl}^\nu\e)\e^\nu= \\ 
	(-\alpha_{kl}^{-1}\alpha'_{kl}r-(\alpha_{kl}^{-1}\alpha'_{kl}\beta_{kl}+\alpha_{kl}^{-1}p_{kl})r^2- 
	\alpha_{kl}^{-2}(\alpha_{kl}p_{kl}(2\beta_{kl}+\beta_{kl}^2r)-\breve p_{kl}(1+\beta_{kl}r)^3)r^3+P_{kl}^\nu\e)\e^\nu. \end{matrix}$$

According to \eqref{pave}, the average value \,$\overline \xi_{kl}=\alpha_{kl}^{-1}(\alpha'_{kl}q_{kl}-p_{kl})-\beta'_{kl},$ therefore
$$\dot r=\big(\overline \xi_{kl}(1+2\beta_{kl}(\varphi)r)^{-1}r^2+Q_{kl}(\varphi,r)r^3+R_{kl}^\nu(t,\varphi,r,\e)\e\big)\e^\nu.$$

Assume $\overline \xi_{kl}\ne 0,$ then the composition of changes \eqref{pz} and \eqref{pave}, if applied not to the system \eqref{sv}, 
but to hamiltonian unperturbed system \eqref{snv}, transforms it into the "shortened" system
$$\dot r=(\overline \xi_{kl}(1+2\beta_{kl}(\varphi)r)^{-1}\!+Q_{kl}(\varphi,r)r)r^2\e^\nu,\quad  
\dot\varphi=(1+q_{kl}(\varphi)r+\mathcal Q_{kl}(\varphi,r)r^2)\e^\nu,$$ 
which describes the motion in the neighbourhood of the arbitrarily chosen cycle. 

In the change \eqref{pave} \,$|r|\le r_*<1,$ therefore, according to \eqref{r*2} and \eqref{r*}, we obtain: 
$$\e^{-\nu}\dot\varphi>1-(2r_*)^{-1}|r|\ge 1/2,\ \ \e^{-\nu}|\dot r|>(2|\overline \xi_{kl}|/3-|Q_{kl}(\varphi,r)r|)r^2>|\overline \xi_{kl}|r^2/3$$ 
under the additional limitation \,$|\overline \xi_{kl}|(3r_*)^{-1}>\max\limits_{|r|\le 1,\,\varphi\in B_{\varrho}} |Q_{kl}(\varphi,r)|.$

As a result, the shortened system and the unperturbed system \eqref{snv} along with it don't have any orbits in the small neighbourhood of the arbitrary cycle. This contradiction means that $\overline \xi_{kl}=0.$ 

Therefore, the change \eqref{pave} transform the system \eqref{ps} into the system \eqref{ps5}. \ $\Box$

\smallskip 
{\bf Remark 4.1.} \ {\sl The function $\xi_{00}(\varphi)$ from \eqref{pave} can be explicitly integrated.
	We have: $p_{00}=3\gamma \alpha^{-1}C(\varphi)S(\varphi)(S^2(\varphi)-C^2(\varphi))=3\alpha_{00}^{-1}\alpha_{00}'(\varphi)/2$ 
	and for $q_{00}$ the chain of equalities can be written: 
	$q_{00}\buildrel \eqref{ps} \over = 
	2\alpha_{00}^{-2} (C^4 + \gamma S^4)\buildrel \eqref{oi} \over = 
	2\alpha_{00}^{-2} (a_{00} - 1 - \gamma + 2C^2 + 2\gamma S^2)\buildrel \eqref{p} \over = 
	2\alpha_{00}^{-2} (2\alpha_{00} - (a_{00} - 1 - \gamma))\buildrel \eqref{ndb} \over =
	2\alpha_{00}^{-2} (2\alpha_{00} - (b_{00}^4 - 2b_{00}^2)).$  
	Then\, $\xi_{00}=\alpha_{00}^{-2}(5/2-2(b_{00}^4 - 2b_{00}^2)\alpha_{00}^{-1})\alpha_{00}',$ 
	 \,$\beta_{00}=(b_{00}^4 - 2b_{00}^2)\alpha_{00}^{-2}(\varphi)-5\alpha_{00}^{-1}(\varphi)/2+(3b_{00}^2-1)(b_{00}^3-b_{00})^{-2}/2,$ 
	i.\,e. $\beta_{00}=\beta_{00}(\alpha_{00}(\varphi)).$ } 

\smallskip
{\bf $4.2.$ The structure of the perturbed part of the system $\bf{\eqref{ps5}}.$ } 
Let us specify, what are the coefficients of the lesser powers of $r$ and $\e$ in the right-hand side of the system \eqref{ps5}.
Denote $\ae^{\circ}=(t,C,S,0)$ and $Z^{\circ}(t,\varphi)=Z(t,\varphi,0,0).$ Then, according to \eqref{ps} and \eqref{sv},
\begin{equation}\label{pf2} \begin{matrix} 
		R_{kl}^{\nu\circ}(t,\varphi)=P_{kl}^{\nu\circ}(t,\varphi)=P_{kl}^{\nu\rho}(t,\varphi,0,0)=C'Y^\nu(\ae^{\circ})-S'X^\nu(\ae^{\circ})=\hfill \\
		\qquad\qquad\quad C'Y_0^\nu(t,C,S)-S'X_0^\nu(t,C,S),\hfill \\
		\Phi_{kl}^{\nu\circ}(t,\varphi)=\Phi_{kl}^{\nu\rho}(t,\varphi,0,0)=\alpha_{kl}^{-1}((S-l)X^\nu(\ae^{\circ})-(C-k)Y^\nu(\ae^{\circ}))=\hfill \\ 
		\qquad\qquad\quad \alpha_{kl}^{-1}((S-l)X_0^\nu(t,C,S)-(C-k)Y_0^\nu(t,C,S));\hfill  \\
		(\!{P_{kl}^{\nu}}_r')^{\!\circ}(t,\varphi)=\alpha_{kl}^{-1}\big(C'((C-k){Y^\nu}_{\!x}'(\ae^{\circ})+(S-l){Y^\nu}_{\!y}'(\ae^{\circ}))-\hfill \\
		\qquad\qquad\qquad S'((C-k){X^\nu}_{\!x}'(\ae^{\circ})+(S-l){X^\nu}_{\!y}'(\ae^{\circ}))\big)=\hfill \\ 
		\qquad\qquad\qquad \alpha_{kl}^{-1}\big(C'((C-k){Y_0^\nu}_{\!x}'(t,C,S)+(S-l){Y_0^\nu}_{\!y}'(t,C,S))-\hfill \\
		\qquad\qquad\qquad S'((C-k){X_0^\nu}_{\!x}'(t,C,S)+(S-l){X_0^\nu}_{\!y}'(t,C,S))\big),\hfill \\ 
		(\!{R_{kl}^{\nu}}_r')^{\!\circ}(t,\varphi)=
		(\!{P_{kl}^{\nu}}_r')^{\!\circ}+\alpha^{-1}_{kl}\alpha'_{kl}\Phi^{\nu\circ}_{kl}-2\beta_{kl}R^{\nu\circ}_{kl},\hfill \\ 
		(\!{R_{kl}^{\nu}}_\e')^{\!\circ}(t,\varphi)=
		(\!{P_{kl}^{\nu}}_\e')^{\!\circ}(t,\varphi)=C'{Y^\nu}_{\!\e}'(\ae^{\circ})-S'{X^\nu}_{\!\e}'(\ae^{\circ})=\hfill \\ 
		\qquad\qquad\qquad 
		\left[ \begin{matrix} C'Y_1^0(t,C,S)-S'X_1^0(t,C,S)\text{ for }\nu=0,\\ 
			C'Y_1^1(\ae^{\circ})-S'X_1^1(\ae^{\circ})\hfill \text{ for }\nu=1.\end{matrix}\right.\hfill  
\end{matrix} \end{equation}

{\bf Lemma 4.2.} \ {\it Such constant $\varrho>0$ exists, that all functions, introduced in \eqref{pf2}, are continuous, two-periodic, uniformly with the respect to $t$ real-analytic with the respect to $\varphi$ on $\mathcal B_{\varrho},$ where $B_{\varrho}$ is from \eqref{BB}, except the function $(\!{R_{kl}^1}_\e')^{\!\circ}\in C_{t,\varphi}^{0,1}(\mathbb{R}^2).$ }

\smallskip
{\it Proof.}\, Consider the functions
\begin{equation}\label{mathcalXY} 
	\mathcal X_{\iota}^\nu(t,\varphi,r)=X_{\iota}^\nu(t,\breve x(\varphi,r),\breve y(\varphi,r)),\  
	\mathcal Y_{\iota}^\nu(t,\varphi,r)=Y_{\iota}^\nu(t,\breve x(\varphi,r),\breve y(\varphi,r)), \end{equation} 
where \,$\iota=0,1,$ $(\iota,\nu)\ne (1,1);$ the arguments \,$\breve x=C+(C-k)\alpha_{kl}^{-1}r(1+\beta_{kl}r),$ 
$\breve y=S+(S-l)\alpha_{kl}^{-1}r(1+\beta_{kl}r)$ are the result of the compositon of changes \eqref{pz} and \eqref{pave}. 

According to the proposition\;4.1 and lemma\;4.1, $|\breve x(\varphi,r)|,|\breve y(\varphi,r))|<\sigma$ with $|r|\le r_*$ and $\varphi\in B_{\varrho},$ therefore the functions, introduced in \eqref{mathcalXY}, are 
continuous, two-periodic, uniformly with the respect to $t$ real-analytic with the respect to \,$\varphi,\,r$\,on $D_{\varphi,r}=\{(t,\varphi,r)\!:\, t\in \mathbb{R},\,\varphi\in B_{\varrho},\,|r|\le r_*\}$ as a composition of the analytic functions.

Let us write them down in the following form 
$$\mathcal X_{\iota}^\nu(t,\varphi,r)=\sum_{j=0}^\infty \mathcal X_{\iota j}^\nu(t,\varphi)r^j,\ \ 
\mathcal Y_{\iota}^\nu(t,\varphi,r)=\sum_{j=0}^\infty \mathcal Y_{\iota j}^\nu(t,\varphi)r^j,$$ 
where $\mathcal X_{\iota j}^\nu(t,\varphi),\mathcal Y_{\iota j}^\nu(t,\varphi)$ are continuous, two-periodic, 
uniformly with the respect to $t$ real-analytic with the respect to $\varphi$ functions on $\mathcal B_{\varrho}$ (see \,lemma\;2.1), 
the series $\mathcal X^\nu,\mathcal Y^\nu$ converge absolutely, when  $|r|<r_*.$
It is obvious that the derivatives ${\mathcal X_0^\nu}_r^{\,'}(t,\varphi,r)$ and ${\mathcal Y_0^\nu}_r^{\,'}(t,\varphi,r)$ have the same properties.

According to \eqref{mathcalXY}, $\mathcal X_{\iota 0}^\nu(t,\varphi)=X_{\iota}^\nu(t,C(\varphi),S(\varphi)),$ 
$\mathcal Y_{\iota 0}^\nu(t,\varphi)=Y_{\iota}^\nu(t,C(\varphi),S(\varphi)),$ and also
\,$\mathcal X_{01}^\nu(t,\varphi)=
{X_0^\nu}_x^{\,'}(t,C(\varphi),S(\varphi))\breve x'(\varphi,0)+{X_0^\nu}_y^{\,'}(t,C(\varphi),S(\varphi))\breve y'(\varphi,0),$ \\
$\mathcal Y_{01}^\nu(t,\varphi)=
{Y_0^\nu}_x^{\,'}(t,C(\varphi),S(\varphi))\breve x'(\varphi,0)+{Y_0^\nu}_y^{\,'}(t,C(\varphi),S(\varphi))\breve y'(\varphi,0),$\\
therefore in the formulae \eqref{pf2} 
\,$(\!{P_{kl}^\nu}_r ')^{\!\circ}(t,\varphi)=C'(\varphi)\mathcal Y_{01}^\nu(t,\varphi)-S'(\varphi)\mathcal X_{01}^\nu(t,\varphi).$ \ $\Box$

\smallskip
As a result, the system \eqref{ps5} can be presented in the following form 
\begin{equation}\label{fix}
	\begin{cases} \dot r=\big(Q_{kl}(\varphi,r)r^3+\big(R_{kl}^{\nu \circ}\!(t,\varphi)+
		(\!{R_{kl}^\nu}_r')^{\!\circ}\!(t,\varphi)r+(\!{R_{kl}^\nu}_\e')^{\!\circ}\!(t,\varphi)\e+R_{kl}^{\nu *}(t,\varphi,r,\e)\big)\e\big)\e^\nu, \\ 
		\dot\varphi=\big(1+q_{kl}(\varphi)r+\mathcal Q_{kl}(\varphi,r)r^2+\big(\Phi_{kl}^{\nu \circ}(t,\varphi)+\Phi_{kl}^{\nu *}(t,\varphi,r,\e)\big)\e\big)\e^{\nu}, 	\end{cases}
\end{equation}
where functions \,$R_{kl}^{\nu *},\Phi_{kl}^{\nu *}$\, are continuous, and \,$R_{kl}^{\nu *}(t,\varphi,r,\e)=O((|r|+\e)^2),$ 
$\Phi_{kl}^{\nu *}(t,\varphi,r,\e)=O(|r|+\e),$ \,$q_{kl}$ is from \eqref{ps}.  

\smallskip
{\bf Proposition 4.2.} \ 
{\it In the system \eqref{fix} $R_{kl}^{\nu *},\Phi_{kl}^{\nu *}\in C_{t,\varphi,r,\e}^{0,1,1,0}(G_{\mathbb{R},r_*,\e_0}^{\varphi,r,\e}),$ 
	where the set \,$G_{\mathbb{R},r_*,\e_0}^{\varphi,r,\e}=\{(t,\varphi,r,\e)\in \mathbb{R}^4\!:\,|r|\le r_*,\,\e\in [0,\e_0]\}.$ } 

\smallskip
Indeed, the functions $R_{kl}^\nu,\Phi_{kl}^\nu,$ introduced in \eqref{ps5}, are of class $C_{t,\varphi,r,\e}^{0,1,1,0}(G_{\mathbb{R},r_*,\e_0}^{\varphi,r,\e})$ 
and continuous in totality of its arguments. 

\section{The generating equation and the generating cycles} 

\smallskip
\subsection{The goal of further transformations.} \ 
The formula \eqref{fix} describes the set of systems obtained from \eqref{sv}. 
Each such system is fixated by the choice of the initial value $b_{kl},$ which meets the conditions \eqref{nd} and \eqref{gamma*}.
Now, we need to distinguish such $b_{kl}$ that, for each one of them, for any sufficiently small $\e>0,$ the corresponding system \eqref{fix} has the two-dimensional invariant surface, and this surface has a non-empty intersection with the arbitrarily chosen neighbourhood of the point $(b_{kl},l).$

To achieve that, it is required to perform the number of an averaging and scaling changes, including the change of the angular variable. These changes yield the system, to which the theorem\;7.1 on the existence of the invariant surface about the existence of the invariant surface can be applied. 

\smallskip
\subsection{Decompositions of two-periodic functions, the Siegel's condition.} \
Following the mentioned plan,  we average the  functions $R_{kl}^o,\,(\!{R_{kl}}'_r)^o$ and $(\!{R_{kl}}'_\e)^o$ in the system \eqref{fix} in such way that the average value of function $R_{kl}^o(t,\varphi)$ can be equal to zero with the correct choice of $b_{kl}.$
However, for this system  $\dot\varphi=1+\ldots$ when $\nu=0,$ $\dot\varphi=\e+\ldots$ for  $\nu=1.$ 
Therefore the  following   averaging changes and their existence conditions will be different for each $\nu.$ 

For continuous, $T$-periodic in $t,$ real analytic and $\omega$-periodic in $\varphi$ 
functions $\eta^\nu(t,\varphi)$ we use the following decomposition depending on the value of parameter $\nu:$
$$\eta^\nu(t,\varphi)=\overline\eta^\nu+\hat\eta^\nu(\varphi)+\tilde\eta^\nu(t,\varphi)\qquad (\nu=0,1),$$
where $\displaystyle \overline\eta^\nu=\frac{1}{\omega T}
\int_0^\omega\int_0^T\eta^\nu(t,\varphi)\,dt\,d\varphi$ is the average value of the  function $\eta^\nu,$
$\displaystyle \hat\eta^0=0$ and $\displaystyle \hat\eta^1=\frac{1}{T}\int_0^T\eta^1(t,\varphi)\,dt-\overline\eta^1.$ 

Then, the  function $\tilde\eta^\nu(t,\varphi)$ has  zero average value with respect to $t,$ 
which implies periodicity of the function $\displaystyle\int_{t_*}^t \tilde\eta^\nu(\tau,\varphi)\,d\tau,$ 
which also has  zero average value by the virtue of the choice of the constant $t_*\in[0,T].$

From now on, we will denote the  derivative with respect to $t$ of any function which has $t$ and $\varphi$ as its arguments by a dot, and the derivative with respect to $\varphi$ by a prime.

\smallskip
{\bf Lemma 5.1.} {\it Assume $\tilde\eta(t,\varphi)$ is continuous, $T$-periodic in $t,$ uniformly with the respect to $t$ real-analytic 
	and \,$\omega$-periodic in $\varphi$ functions, the following condition, called Siegel's condition, holds for the periods $\omega$ and $T$
	\begin{equation}\label{zig} |m\omega+nT|>\vartheta(|m|+|n|)^{-\tau}\quad (\vartheta>0,\ \tau\ge 1,\ m,n\in\mathbb{Z},\ m^2 + n^2\neq 0). 
	\end{equation} 
	Then the equation
	\begin{equation}\label{zigeq} \dot{\tilde\chi}(t,\varphi)+\tilde\chi'(t,\varphi)=\tilde\eta(t,\varphi) \end{equation}
	has a singular solution $\tilde\chi(t,\varphi),$ which has the same properties as the function $\tilde\eta(t,\varphi).$ 
	It's Fourier series expansion is term-wise differentiable with the respect to $t$ and $\varphi.$ } 

\smallskip
{\it Proof.}\, 
According to lemma\;2.1, such constant $\varrho>0$ exists, that \,$\tilde\eta$ is two-periodic function on $\mathcal B_{\varrho}.$ For this function such constant $M_{\eta}>0$ exists, 
that $M_{\eta}=\max_{\,(t,\varphi)\in {\mathcal B_{\varrho}}}|\tilde\eta(t,\varphi)|,$ 
and the coefficients of the expansion \,$\displaystyle \eta_n(t)=\omega^{-1}\int_0^\omega \tilde\eta(t,\varphi)e^{-2\pi i n \varphi/\omega}d\varphi$ 
from \eqref{varpi} of the function $\tilde\eta$ into a Fourier series 
$\displaystyle \tilde\eta(t,\varphi)=\sum_{n=-\infty}^{+\infty}\eta_n(t) e^{2\pi i n \varphi/\omega}$ 
with $\varphi\in[0,\omega]$ meet the condition \eqref{cofu} \,$|\eta_n(t)|\le M_{\eta}e^{-2\pi|n|\varrho/\omega}$ for any $t\in \mathbb{R};$ 
moreover, $\eta_0(t)\equiv 0.$

\smallskip
Consider the linear equation
\begin{equation}\label{lnu} \dot{\chi}_n(t)=(-2\pi in/\omega) \chi_n(t) + \eta_n(t)\quad (n\in \mathbb{Z},\ n\ne 0). \end{equation}

It's singular $T$-periodic solution is the following function:
\begin{equation}\label{chin}
	\chi_n^p (t)=  d_n^{-1}\int_{t-T}^{t} e^{2\pi i n(s-t)/\omega}\eta_n(s)ds,\ \ d_n = 1 - e^{-2\pi i nT/\omega}.
\end{equation}

Using the expansion $1=e^{2\pi i m}\ (m\in\mathbb{Z})$ and Euler's formula, we obtain: 
$$\forall\,m\in\mathbb{Z}:\ \ |d_n|=2|\sin(\pi(m+nT/\omega))|.$$

Denote $$m_n = -[nT/\omega + 1/2],$$ 
then the argument $|\pi(m_n+nT/\omega)|\le \pi/2,$ because $0\le m_n+nT/\omega+1/2<1.$ 
For the same reason the evaluation $|m_n| \le (T/\omega+1)|n|$ is correct.

Due to $|\sin x|\ge (2/\pi)|x|$ when $|x|\le \pi/2,$ we deduce that 
$$|d_n|\ge 4\omega^{-1}|m_n\omega+nT|>4\vartheta\omega^{-1}(|m_n| + |n|)^{-\tau}\ge 4\vartheta\omega^{-1}((T/\omega+2)|n|)^{-\tau},$$ 
because $n\ne 0$ and Siegel's condition \eqref{zig} holds, in particular, for $m=m_n.$  

As a result, according to \eqref{cofu} and formula \eqref{chin}, we obtain  
$$|\chi_n^p(t)| < M_* ((T/\omega+2)|n|)^{\tau}e^{-2\pi|n|\varrho/\omega},\quad M_* = M_{\eta}T\omega(4K)^{-1}.$$
Then such constant $M_{\chi}>0$ exist, that $M_*((T/\omega+2)|n|)^{\tau}e^{-2\pi|n|(\varrho/4)/\omega}\le M_{\chi}.$ 

Henceforth, the evaluation similar to \eqref{cofu} is correct for the function $\chi_n^p(t):$
\begin{equation}\label{chiest} 
	\forall\,t\in \mathbb{R}\!:\ \ |\chi_n^p(t)| < M_{\chi} e^{-2\pi|n|(3\varrho/4)/\omega} \quad (n\ne 0). 
\end{equation} 
Therefore the following evaluation is correct for the function $\chi_n^p(t)e^{2\pi in\varphi/\omega}$ on $\mathcal B_{\varrho/2}$ 
$$|\chi_n^p(t)e^{2\pi in\varphi/\omega}| < M_{\chi}e^{-2\pi |n|(\varrho/4)/\omega},$$ 
This evaluation guarantees, that for any \,$t_*\in \mathbb{R}$\, the function \,$\chi_n^p(t_*)e^{2\pi in\varphi/\omega}$\, is uniformly analytic with the respect to $\varphi$ on $B_{\varrho/2}.$

According to the Weierstrass' theorem about the uniform convergence of the functional series, the continuous,
$T$-periodic in $t$ and $\omega$-periodic in $\varphi,$ taking real values for any $\varphi$ function
\begin{equation}\label{chitphi} 
	\tilde\chi(t,\varphi) = \sum_{n=-\infty}^{\infty}\ \chi_n^p(t)e^{2\pi i n\varphi/\omega}\qquad (\chi^p_0(t)\equiv 0) 
\end{equation}
is analytic with the respect to $\varphi$ on $\mathcal B_{\varrho/2}.$

The coefficients of the series, obtained by the term-wise differentiation of the series \eqref{chitphi} with the respect to $\varphi,$ 
satisfy the exponential evaluation independent of $t$ and $\varphi,$ therefore it's sum is equal to the partial derivative $\tilde\chi'(t,\varphi)$ 
on $\mathcal B_{\varrho/2}.$

According to \eqref{chiest}, \eqref{cofu} and \eqref{lnu}, the evaluation
$|\dot{\chi}_n^p(t)| < (2\pi nM_{\chi}/\omega + M_{\eta})e^{-2\pi |n|(3\varrho/4)/\omega}$ is correct. 
Then, for any $(t,\varphi)\in \mathcal B_{\varrho/2},$ $|\dot{\chi}_n^p(t)e^{2\pi in\varphi/\omega}| < 
(2M_{\chi}\pi n/\omega + M_{\eta})e^{-2\pi |n|(\varrho/4)/\omega},$  i.\,e. coefficients of the series, 
obtained by the term-wise differentiation of the series $\eqref{chitphi}$ with the respect to $t,$ 
satisfy the exponential evaluation independent of $t$ and $\varphi,$ 
therefore it's sum is equal to the partial derivative $\dot{\tilde\chi}(t,\varphi)$ on $\mathcal B_{\varrho/2}.$

The direct substitution shows that the function  $\tilde\chi(t,\varphi)$ is a solution of the equation \eqref{zigeq} 
on $\mathcal B_{\varrho/2}.$ \ $\Box$

\subsection{The choice of generating cycles.} \  The function 
\begin{equation}\label{rio} 
	R_{kl}^{\nu \circ}(t,\varphi)=C'(\varphi)Y_0(t,C(\varphi),S(\varphi))-S'(\varphi)X_0(t,C(\varphi),S(\varphi)), \end{equation}
introduced in \eqref{pf2} will  play the key role in our analysis. 
In \eqref{rio} \,$(C(\varphi),S(\varphi))$ is a real analytic $\omega_{kl}$-periodic solution of the  initial value problem 
of the system \eqref{cs} with the initial values $C(0)=b_{kl},\,S(0)=l,$ the  parameter $b_{kl}$ is any number from \eqref{nd} and \eqref{gamma*}
and the  period $\omega_{kl}$ is calculated in \eqref{ome}. 

Using the mentioned above decomposition, we write out $R^{o}_{kl}(t,\varphi)$ as the following sum:
$$R_{kl}^{\nu \circ}=\overline {R_{kl}^{\nu \circ}}+\widehat {R_{kl}^{\nu \circ}}(\varphi)+\widetilde {R_{kl}^{\nu \circ}}(t,\varphi),$$
where 
\begin{equation}\label{eve} 
	\overline {R_{kl}^{\nu \circ}}=\overline {R_{kl}^{\nu \circ}}(b_{kl})=
	{1\over T\omega_{kl}}\int_0^{\omega_{kl}}\int_0^T R_{kl}^{\nu \circ}(t,\varphi)dt\,d\varphi.
\end{equation} \
is  the average value of $R_{kl}^o.$ 

In case $\nu=0$ consider the following equation
\begin{equation}\label{gn}
	\dot{\tilde g}_{kl}^0(t,\varphi)+\tilde g_{kl}^{0'}(t,\varphi)=\widetilde R_{kl}^o(t,\varphi)
\end{equation}
According to the lemma 5.1, this equation has the singular solution $\tilde g_{kl}^0(t,\varphi)$ with the same properties 
as the function $\widetilde R_{kl}^o(t,\varphi).$ 

Now we can formulate the nondegeneracy condition as 
\begin{equation} \label{ud}  K_{kl}^\nu\ne 0;\qquad 
	K_{kl}^0=\overline{(\!{R_{kl}^0}'_r)^{\circ}-\tilde g_{kl}^{0\,'}q_{kl}},\ \ 
	K_{kl}^1=\overline{(\!{R_{kl}^1}'_r)^{\circ}-\widehat {R_{kl}^{1 \circ}}q_{kl}}. \end{equation}

According to \eqref{pf2}, the function ${\tilde g}_{kl}^0(t,\varphi)$ and the constants $K_{kl}^\nu$ depend only on the perturbations $X_0,Y_0$ of the system \eqref{sv}.

{\bf Definition 5.1.} \ {\it The solution of the generating (bifurcation) equation
	\begin{equation}\label{pu} \overline {R_{kl}^{\nu \circ}}(b_{kl})=0\end{equation} 
	with $\overline {R_{kl}^{\nu \circ}}$ from \eqref{eve} is an admissable for the system \eqref{sv} solution and is denoted as $b_{kl}^{\star},$ if it meets the conditions \eqref{nd} and \eqref{gamma*}, after the choice $b_{kl}=b_{kl}^{\star}$ the nondegeneracy condition \eqref{ud} holds and, 
	for $\nu=0,$ the Siegel's condtion \eqref{zig} for the periods $T$ and $\omega_{kl}$ from \eqref{ome} holds.}

\smallskip
From now on, let us fix some admissable $b_{kl}^{\star}$ for the system \eqref{sv} solution of the generating equation \eqref{pu}. Fixing this solution chooses the $\omega_{kl}$-periodic solution of the initial value problem $C\!S_{kl}(\varphi,0,b_{kl},l)$ of the system \eqref{cs}, which parametrizes the cycle related to the class $|k|+|l|.$

All subsequent changes and systems are also fixed by the choice of the $b_{kl}^{\star}.$ In particular, the conditions \eqref{zig} and \eqref{ud} hold, and in the system \eqref{fix} the function $R_{kl}^{\circ}(t,\varphi),$ introduced in \eqref{pf2}, has zero average value.

\smallskip
Designations. {\sl From now on, the subindex $_{kl}$ of the functions and the constants, 
which represents the fixation of the generating parameter $b_{kl}^{\star},$ will be changed to subindex $_{\star}$ for brevity. 
For example, $R_{kl}^{\circ}=R_{\star}^{\nu \circ},\,K_{kl}^\nu=K_{\star}^\nu.$ }

\smallskip
{\bf Definition 5.2.} \ 
{\it For any admissable parameter $b_{kl}^*,$ cycle $GC_{\star}=\{(x,y)\!:\,x=C(\varphi),\,y=S(\varphi)\ \ (\varphi\in \mathbb{R})\}$\, of the unperturbed system \eqref{snv}, parametrized by the solution of the initial value problem of system \eqref{cs} with the initial values $b_{kl}^*,l,$ is called the generating cycle.	This cycle is a generatrix of the cylindrical invariant surface 
	$C\!I\!S_{\star}=\{(x,y,t)\!:\,x=C(\varphi),\,y=S(\varphi),\,t\in \mathbb{R}\ \ (\varphi\in \mathbb{R})\}$\, 	of the system \eqref{snv}.  }

\smallskip
{\bf Remark 5.1.} \ {\sl 
	We will  be interested in systems \eqref{sv} with a nonempty set of admissible solutions.
	For each  $b_{kl}^{\star},$ it will be proved that the perturbed system \eqref{sv} retains a two-periodic invariant surface homeomorphic to a torus, which is obtained by factoring the time $t$ with	respect to the period $T$ in a small  neighborhood of the cylindric surface $C\!I\!S_{\star},$ for each sufficiently small $\e.$}

\section{The construction of invariant surfaces} 

\subsection{Secondary radial averaging.} \
Let us average the coefficients of $r\e$ and annul the coefficient of $\e^2$ in the first equation of the system \eqref{fix}.

\smallskip
{\bf Proposition 6.1}  {\it Such constant $\varrho>0$  exists, that the solutions of the equations
	\begin{equation} \label{fgh} \begin{matrix} 
			\dot{\tilde h}_{\star}^0+\tilde h_{\star}^{0\,'}=(\!{R_{\star}^0}'_r)^{\!\circ}-\tilde g_{\star}^{0\,'}q_{\star}-K_{\star}^0,\ \ 
			\tilde g_{\star}^0\hbox{ from } \eqref{gn},\hfill \\
			\tilde f_{\star}^{0\,'}+\dot {\tilde f}_{\star}^0=(\!{R_{\star}^0}'_r)^{\!\circ}(\overline g_{\star}^0+\tilde g_{\star}^0)+
			(\!{R_{\star}^0}'_\e)^{\!\circ}-\tilde g_{\star}^{0\,'}(\Phi_{\star}^{0 \circ}+(\overline g_{\star}^0+\tilde g_{\star}^0)q_{\star});
			\phantom{\cfrac cc}\\ 
			\hat g_{\star}^{1\,'}=\widehat {R_{\star}^{1 \circ}},\ \dot{\tilde g}_{\star}^1=\widetilde {R_{\star}^{1 \circ}},\ \    
			\hat h_{\star}^{1'}=\widehat {(\!{R_{\star}^1}'_r)^{\!\circ}}-\widehat {R_{\star}^{1 \circ}}q_{\star},\ \dot{\tilde h}_{\star}^1=
			\widetilde {(\!{R_{\star}^1}'_r)^{\!\circ}},\hfill \\  
	\end{matrix} \end{equation}
	where 
	\,$\overline g_{\star}^0=\big(\overline {\tilde g_{\star}^{0\,'}\Phi_{\star}^{0 \circ}}-\overline {(\!{R_{\star}^0}'_\e)^{\!\circ}}-
	\overline {\tilde g_{\star}^0((\!{R_{\star}^0}'_r)^{\!\circ}-\tilde g_{\star}^{0\,'}q_{\star})}\big)/K_{\star}^0,$\, 
	\,$\overline g_{\star}^1=\big(\overline {\widehat {R_{\star}^{1 \circ}}\Phi_{\star}^{1 \circ}}-\overline {(\!{R_{\star}^1}'_\e)^{\!\circ}}-
	\overline {\hat g_{\star}^1((\!{R_{\star}^1}'_r)^{\!\circ}-\widehat {R_{\star}^{1 \circ}}q_{\star})}\big)/K_{\star}^1,$\, 
	exist, are singular and are continuous, $T$-periodic and continuously differentiable in $t,$ uniformly with the respect to $t$ real-analytic and 
	\,$\omega_{\star}$-periodic in $\varphi$ on $\mathcal B_{\varrho}$ from \eqref{BB}, and solutions of the equations
	\begin{equation} \label{f1} \hat f_{\star}^{1'}=\widehat{\mathfrak F}_{\star}^1,\ \ \dot{\tilde f}_{\star}^1=\widetilde {\mathfrak F}_{\star}^1,
	\end{equation}
	where \,$\mathfrak F_{\star}^1=\overline g_{\star}^1((\!{R_{\star}^1}'_r)^{\!\circ}-\widehat R_{\star}^{1 \circ}q_{\star})-
	\big(\widehat R_{\star}^{1 \circ}\Phi_{\star}^{1 \circ}-
	(\!{R_{\star}^1}'_\e)^{\!\circ}-\hat g_{\star}^1((\!{R_{\star}^1}'_r)^{\!\circ}-\widehat R_{\star}^{1 \circ}q_{\star})+\tilde g_{\star}^{1'}\big),$ 
	are also singular and are of class $C^1(\mathbb{R}^2).$ }

{\it Proof.}\, Lemma 5.1 can be applied to the equations $(\ref{fgh}^0),$ which are analogous to the equations \eqref{gn}, because the constants $K_{\star}^0$ and $\overline g_{\star}^0$ turn the average values of the right-hand sides of the equations into zero. The properties of the solutions of the equations $(\ref{fgh}^1)$ are obvious.

The solutions of the equations \eqref{f1} are continuously differentiable, because the function $(\!{R_{\star}^1}'_\e)^{\!\circ},$ which is continuously differentiable only in $\varphi$ (see Lemma 4.2),
is the part of the function $\mathfrak F_{\star}^1.$

The uniqueness of each solution from each equation from \eqref{fgh},\,\eqref{f1} is implied from their average values being equal to zero. \ $\Box$

\smallskip
{\bf Remark 6.1.} \ {\sl The real-analytic functions $X_1^0(t,x,y)$ and $Y_1^0(t,x,y),$ are separated in the system \eqref{sv}, specifically to have the possibility to find the two-periodic function $\tilde f_{\star}^0$ from \eqref{fgh}. These functions are part only of the function $(\!{R_{\star}^0}'_\e)^{\!\circ},$ introduced in \eqref{pf2}.}

\smallskip
Consider the two-periodic in $t$ and $\varphi$ change of the angular variable
\begin{equation}\label{save} 
	r=u+G_{\star}^\nu(t,\varphi,\e)\e+H_{\star}^\nu(t,\varphi,\e)u\e+F_{\star}^\nu(t,\varphi,\e)\e^2\quad (\nu=0,1), \end{equation}
where
\,$G_{\star}^0=\overline g_{\star}^0+\tilde g_{\star}^0(t,\varphi),$ $H_{\star}^0=\tilde h^0_{\star}(t,\varphi),$ 
$F_{\star}^0=\tilde f^0_{\star}(t,\varphi);$
$G_{\star}^1=\overline g^1_{\star}+\hat g^1_{\star}(\varphi)+\tilde g^1_{\star}(t,\varphi)\e,$  
$H_{\star}^1=\hat h^1_{\star}(\varphi)+\tilde h^1_{\star}(t,\varphi)\e,$ $F_{\star}^1=\hat f^1_{\star}(\varphi)+\tilde f^1_{\star}(t,\varphi)\e.$ 

\smallskip
{\bf Proposition 6.2.} \ {\it Such constants $u_*>0$ and $\e_*\ (0<\e_*\le \min\{\e_0,1\})$ exist, that in \eqref{save} 
	\,$|u+G_{\star}^\nu(t,\varphi,\e)\e+H_{\star}^\nu(t,\varphi,\e)u\e+F_{\star}^\nu(t,\varphi,\e)\e^2|\le r_*$ 
	for any \,$t\in \mathbb{R},\,\varphi\in B_{\varrho},\,|u|\le u_*,\,\e\in [0,\e_*].$ } 

\smallskip 
Indeed, according to lemma\;2.1 and proposition\;6.1, the absolute values of the functions $G_{\star}^\nu,H_{\star}^\nu,F_{\star}^0$ reach their maximum at the set $\mathcal B_{\varrho}.$ 

As a result of the choosing specific $\varrho,u_*,\e_*$, the substitution \eqref{save} is possible without the change of the properties of the functions in the system \eqref{fix}.

\smallskip 
{\bf Lemma 6.1.} \ {\it The change \eqref{save} transforms the system \eqref{fix} into the following system 
	\begin{equation}\label{savesys} 
		\begin{cases} \dot u=(K_{\star}^\nu u\e+U_{\star}^\nu(t,\varphi,u,\e))\e^\nu,\\
			\dot\varphi=(1+\Theta_{\star}^{\nu}(t,\varphi)\e+q_{\star}(\varphi)u+\breve \Phi_{\star}^\nu(t,\varphi,u,\e))\e^\nu,\end{cases} 
	\end{equation}
	where \,$\Theta_{\star}^0=\Phi_{\star}^{0 \circ}+(\overline g_{\star}^0+\tilde g_{\star}^0)q_{\star},$\, 
	$\Theta_{\star}^1=\Phi_{\star}^{1 \circ}+(\overline g_{\star}^1+\hat g_{\star}^1)q_{\star},$ \,$q_{\star}$ is from \eqref{ps},
	continuous, two-periodic in $t$ and $\varphi$ functions $U_{\star}^\nu,\breve \Phi_{\star}^\nu\in 
	C_{t,\varphi,u,\e}^{0,1,1,0}(G_{\mathbb{R},u_*,\e_*}^{\varphi,u,\e}),$ 
	where $G_{\mathbb{R},u_*,\e_*}^{\varphi,u,\e}=\{(t,\varphi,u,\e)\in \mathbb{R}^4\!:\,|u|\le u_*,\,\e\in [0,\e_*]\},$ 
	and \,$U_{\star}^\nu=O((|u|+\e)^3),$ $\breve \Phi_{\star}^\nu=O((|u|+\e)^2).$} 

{\it Proof.}\,  Substituting the change \eqref{save} into the second equation of the system \eqref{fix}, 
we obtain: \,$\dot\varphi=(1+q_{\star}(u+G_{\star}^\nu\e)+\Phi_{\star}^{\nu \circ}\e+O((|u|+\e)^2))\e^\nu.$ 
This allows to find $\Theta_{\star}^\nu.$ The rest of the terms make the function  $\breve \Phi_{\star}^\nu(t,\varphi,u,\e).$ 

Now, differentiating the change (\ref{save}) with the respect to $t,$ with the respect to the system $(\ref{fix})$ and (\ref{savesys}) and dividing the parts of equality by $\e^\nu,$ we obtain the following equality
$$\begin{matrix}
	\e(R_{\star}^{\nu \circ}+(\!{R_{\star}^\nu}'_r)^{\!\circ}(u+G_{\star}^\nu\e)+(\!{R_{\star}^\nu}'_\e)^{\!\circ}\e)+O((|u|+\e)^3)= \\ 
	\phantom{\cfrac cc}=K_{\star}^\nu u\e+G_{\star}^{\nu\,'}\e(1+\Theta_{\star}^\nu\e+q_{\star}u)+H_{\star}^{\nu'}u\e+F_{\star}^{\nu'}\e^2+
	(\dot G_{\star}^\nu+\dot H_{\star}^\nu u+\dot F_{\star}^\nu \e)\e^{1-\nu},\end{matrix}$$ 
where all terms, having a sum of degree of $u$ and $\e$ are three or higher, are collected
in the continuous with the respect to $t,\varphi,u,\e$ from $G_{\mathbb{R},u_*,\e_*}^{\varphi,u,\e}$ function $O((|u|+\e)^3).$

For $\nu=0,$ equalizing the coefficients of $\e,u\e,\e^2,$ we obtain the first equations from \eqref{fgh}.
For $\nu=1,$ equalizing the coefficients of the same powers as the case of $\nu=0,$ we obtain the equations 
$R_{\star}^{1 \circ}=\hat g_{\star}^{1'}+\dot{\tilde g}_{\star}^1,$ 
$\hat h_{\star}^{1'}=\widehat {(\!{R_{\star}^1}'_r)^{\!\circ}}-\widehat R_{\star}^{1 \circ}q_{\star}$ and
$\dot{\tilde h}_{\star}^1=\widetilde {(\!{R_{\star}^1}'_r)^{\!\circ}},$
$(\overline g_{\star}^1+\hat g_{\star}^1)(\!{R_{\star}^1}'_r)^{\!\circ}+(\!{R_{\star}^1}'_\e)^{\!\circ}=\tilde g_{\star}^{1'}+
\hat g_{\star}^{1'}(\Phi_{\star}^{1 \circ}+(\overline g_{\star}^1+\hat g_{\star}^1)q_{\star})+\hat f_{\star}^{1'}+\dot{\tilde f}_{\star}^1,$ 
which can be reduced to the equations $(\ref{fgh}^1)$ and \eqref{f1}. \ $\Box$

\smallskip
{\bf $6.2.$ Angular averaging.} \
Let us perform one more averaging: we average the continuous, $T$-periodic in $t,$ uniformly with the respect to $t$ real-analytic and $\omega_{\star}$-peridoic with the respect to $\varphi$ function $\Theta_{\star}^\nu(t,\varphi)$ in the second equation of the system (\ref{savesys}).

\smallskip
{\bf Proposition 6.3.} \ {\it Such constant $\varrho>0$ exists, that the solutions of the equations
	\begin{equation} \label{delta} 
		\tilde \delta_{\star}^{0'}+\dot{\tilde\delta}_{\star}^0=\Theta_{\star}^0(t,\psi)-\overline {\Theta_{\star}^0};\quad 
		\hat \delta_{\star}^{1'}=\widehat\Theta_{\star}^1(\psi),\ \ \dot{\tilde\delta}_{\star}^1=\widetilde\Theta_{\star}^1(t,\psi) \end{equation}
	are continuous, $T$-periodic and continuously differentiable in $t,$ 
	uniformly with the respect to $t$ real-analytic and \,$\omega_{\star}$-periodic in $\psi$ 
	on $\mathcal B_{\varrho}^{\psi}=\{(t,\psi)\!: t\in \mathbb{R},\,\psi\in B_{\varrho}^{\psi}\}$\,  
	with \,$B_{\varrho}^{\psi}=\{({\rm Re}\,\psi,{\rm Im}\,\psi)\!:{\rm Re}\,\psi\in \mathbb{R},\,|{\rm Im}\psi|\le \varrho\}.$ } 

Indeed, the function $\Theta_{\star}^0(t,\psi),$ introduced in \eqref{savesys}, meets the requirements of lemma\;5.1. Everything else is obvious.

\smallskip
Consider two-periodic in $t$ and $\psi$ change of the angular variable 
\begin{equation}\label{phi} \varphi=\psi+\Delta_{\star}^\nu(t,\psi,\e)\e, \end{equation}
where \,$\Delta_{\star}^0=\tilde\delta_{\star}^0(t,\psi),$ $\Delta_{\star}^1=\hat\delta_{\star}^1(\psi)+\tilde\delta_{\star}^1(t,\psi)\e.$

\smallskip
{\bf Proposition 6.4.} \ {\it Such constant $\e_*>0$ exists, that 
	\,$|\psi+\Delta_{\star}^\nu(t,\psi,\e)\e|\le \varrho$ for any $(t,\psi)\in \mathcal B_{\varrho}^{\psi}$ and $\e\in [0,\e_*].$ } 

\smallskip
As a result of choosing specific $\varrho,\e_*,$ the substitution \eqref{phi} is possible without the change of the properties of the functions in the system \eqref{savesys}.

\smallskip
{\bf Lemma 6.2.} \ {\it The continuous, $T$-periodic in $t,$ uniformly with the respect to $t$ real-analytic and \,$\omega_{\star}$-periodic in~$\psi$ change of the angular variable transfroms the system $\eqref{savesys}$ into the following system
	\begin{equation}\label{phisys}
		\dot u=(K_{\star}^\nu u\e+\breve U_{\star}^\nu(t,\psi,u,\e))\e^\nu,\ \ 
		\dot\psi=(1+\overline {\Theta_{\star}^\nu}\e+q_{\star}(\psi)u+\breve \Psi_{\star}^\nu(t,\psi,u,\e))\e^\nu, \end{equation} 
	where continuous, two-periodic in $t$ and $\psi$ functions \,$\breve U_{\star}^\nu,\breve \Psi_{\star}^\nu\in 
	C_{t,\psi,u,\e}^{0,1,1,0}(G_{\mathbb{R},u_*,\e_*}^{\,\psi,u,\e}),$ 
	\,$G_{\mathbb{R},u_*,\e_*}^{\,\psi,u,\e}=\{(t,\psi,u,\e)\in \mathbb{R}^4\!:\,|u|\le u_*,\,\e\in [0,\e_*]\};$ 
	moreover, \,$\breve U_{\star}^\nu=O((|u|+\e)^3),$ $\breve \Psi_{\star}^\nu=O((|u|+\e)^2).$ }  

\smallskip
{\it Proof.}\, In the first equation of the system \eqref{phisys} $\breve U_{\star}^\nu(t,\psi,u,\e)=U_{\star}^\nu(t,\psi+\Delta_{\star}^\nu\e,u,\e).$

Differentiating the change \eqref{phi} with the respect to $t,$ with the respect to
the systems \eqref{savesys} and \eqref{phisys} and dividing the parts of equality by $\e^\nu,$ we obtain equality
$$\begin{matrix}  \phantom{\cfrac{c}{c}} 1+\Theta_{\star}^{\nu}(t,\psi+\Delta_{\star}^\nu\e)\e+q_{\star}(\psi+\Delta_{\star}^\nu\e)u+
	\breve \Phi_{\star}^\nu(t,\psi+\Delta_{\star}^\nu\e,u,\e)= \\
	=\dot \Delta_{\star}^\nu\e^{1-\nu}+(1+{\Delta_{\star}^\nu}'\e)(1+\overline {\Theta_{\star}^\nu}\e+q_{\star}(\psi)u+
	\breve \Psi_{\star}^\nu(t,\psi,u,\e)),
\end{matrix}$$
where 
$\Theta_{\star}^{\nu}(t,\psi+\Delta_{\star}^\nu\e)=\Theta_{\star}^\nu(t,\psi)+O(\e),$ $q_{\star}(\psi+\Delta_{\star}^\nu\e)=q_{\star}(\psi)+O(\e).$
Therefore, equalizing the functions that serve as coefficients of $\e,$ for $\nu=0,$ we obtain the first of the equations \eqref{delta}, 
for $\nu=1,$ we obtain the equation 
$\Theta_{\star}^1(t,\psi)=\hat \delta_{\star}^{1'}(\psi)+\dot{\tilde\delta}_{\star}^1(t,\psi)+\overline\Theta_{\star}^1,$ 
which dissolves into the two last equations from \eqref{delta}. 

The rest of the terms in the equality make the function $\breve \Psi_{\star}^\nu(t,\psi,u,\e),$ which has the properties, described in the theorem. \ 
$\Box$

\smallskip
Notice, that the inverse to the change $(\ref{phi})$ can be written as
\begin{equation}\label{psi}
	\psi=\varphi+\Omega_{\star}^\nu(t,\varphi,\e)\e,
\end{equation} 
where \,$\Omega_{\star}^0=-\tilde\delta_{\star}^0(t,\varphi)+O(\e),$ \ $\Omega_{\star}^1=-\hat\delta_{\star}^1(\varphi)+O(\e)$ are
$T$-periodic in $t,$ $\omega_{\star}$-periodic and uniformly with the respect to $t$ real-analytic with the respect to $\varphi$ function.

\smallskip
{\bf $6.3.$ The scaling.} \ 

{\bf Lemma 6.3.} \ {\it The change
	\begin{equation}\label{scale} u=v\e^{3/2}\qquad (0<\e\le \e_*) \end{equation}
	transforms the system $(\ref{phisys})$ into the following system
	\begin{equation}\label{scsys}
		\dot v=( K_{\star}^\nu v\e+V_{\star}^\nu(t,\psi,v,\e)\e^{3/2})\e^\nu,\ \ 
		\dot\psi=(\mathcal K_{\star\e}^\nu+\Psi_{\star}^\nu(t,\psi,v,\e)\e^{3/2})\e^\nu,
	\end{equation}
	where $\mathcal K_{\star\e}^\nu=1+\overline{\Theta_{\star}^\nu}\e>0,$ $V_{\star}^\nu=O((|v|\e^{3/2}+\e)^3)\e^{-3},$ $\Psi_{\star}^\nu=q_{\star}(\psi)v+O((|v|\e^{3/2}+\e)^2)\e^{-3/2}$ 
	are continuous, two-periodic in $t$ and $\psi$ function of class \,$C_{t,\psi,v,\e}^{0,1,1,0}(G_{\mathbb{R},v_*,\e_*}^{\psi,v,\e}),$ 
	where \,$G_{\mathbb{R},v_*,\e_*}^{\psi,v,\e}=\{(t,\psi,v,\e)\in \mathbb{R}^4\!:\,|v|\le v_*,\,\e\in [0,\e_*]\}$ 
	and \,$v_*=u_*$ from $(\ref{phisys}).$  }

\smallskip
Introduce such constants $M=M_{\e_*}>0,$ and $L=L_{\e_*}>0,$ that
\begin{equation}\label{lip} \begin{matrix}
		\max_{G_{\mathbb{R},v_*,\e_*}^{\psi,v,\e}} \{|V_{\star}^\nu|,|\Psi_{\star}^\nu|\} \le M;\\  
		\forall\,(t,\psi_1,v_1,\e),(t,\psi_2,v_2,\e)\in G_{\mathbb{R},v_*,\e_*}^{\psi,v,\e}\!: \\
		\hfill |V_{\star}^\nu (t,\psi_2,v_2,\e)-V_{\star}^\nu (t,\psi_1,v_1,\e)|\le L(|\psi_2 - \psi_1| + |v_2 - v_1|),\\  
		\hfill |\Psi_{\star}^\nu (t,\psi_2,v_2,\e)-\Psi_{\star}^\nu (t,\psi_1,v_1,\e)|\le L(|\psi_2 - \psi_1| + |v_2 - v_1|), \end{matrix} \end{equation}
i.\,e. the Lipschitz' conditions with the respect to $\psi$ and $v$ with the global constant $L$ hold.

By decreasing, if needed, the value of the constant $\e_*,$ which, obviously, doesn't affect the values of the constants $M$ and $L,$ the following inequalities can be made correct 
\begin{equation}\label{fineq} \mathcal M=M/|K_{\star}^\nu|\le v_*\e_*^{-1/2},\ \ \mathcal L=4L/|K_{\star}^\nu|\le \e_*^{-1/2}, \end{equation}
moreover, the evaluations $$L(1+\mathcal L\e^{1/2})\e^{1/2}\le |K_{\star}^\nu|/2,\ \   2L|K_{\star}^\nu|^{-1}\e^{1/2}\le 1/2,\ \ 
2L(1+\mathcal L\e^{1/2})|K_{\star}^\nu|^{-1}\le \mathcal L$$ are implied for any $\e\in [0,\e_*]$ from the inequality $(\ref{fineq}_2).$

Perform one more scaling, by setting 
\begin{equation}\label{scale2} v=z\e^{1/2}\quad (|z|\le z_{\e},\ 0<\e\le \e_*),\end{equation} 
where $z_{\e}=v_*\e^{-1/2}.$ Then the system \eqref{scsys} can be written as
\begin{equation}\label{scsys2}
	\dot z=(K_{\star}^\nu z+V_{\star}^\nu(t,\psi,z\e^{1/2},\e))\e^{\nu+1},\ \ 
	\dot\psi=(\mathcal K_{\star\e}^\nu+\Psi_{\star}^\nu(t,\psi,z\e^{1/2},\e)\e^{3/2})\e^\nu.
\end{equation} 

\smallskip
{\bf $6.4.$ The existence and the properties of the invariant surface.} \

\smallskip
{\bf Theorem 6.1.} \ {\it For any $\e\in (0,\e_*],$ the system \eqref{scsys2} has 
	the two-dimensional invariant surface 
	\begin{equation}\label{isep} z=A_\e^\nu(t,\psi)\quad (t,\psi\in \mathbb{R}), \end{equation} 
	which is parametrized by the continuous, $T$-periodic in $t,$ $\omega_{\star}$-periodic function $A_\e^\nu.$ The Lipschitz' condition with the respect to $\psi$ with the global constant $\mathcal L$ holds for this function.
	Moreover, \,$\max_{\,t,\psi\in \mathbb{R}} |A_\e^\nu(t,\psi)|\le \mathcal M,$ where $\mathcal M$ and $\mathcal L$ are introduced in $(\ref{fineq}).$ }

\smallskip
{\it Proof.}\,  
1) Assume that in the system \eqref{scsys2} the constant $K_{\star}^\nu>0.$ Otherwise, it can be achieved by substituting $t$ for $-t.$

2) Denote as $\mathcal F$ the metric space of the real, continuous, $T$-periodic in $t,$ $\omega_\star$-periodic in $\psi$ 
functions $F(t,\psi)$ with the metric \,$\rho(F_1,F_2)=\max\nolimits_{\,t,\psi\in \mathbb{R}}|F_2(t,\psi)-F_1(t,\psi)|.$
Also assume the the evaluation  $|F(t,\psi)|\le \mathcal M$ and the Lipschitz' condition with the respect to $\psi$ with the global constant $\mathcal L$ hold for these functions. 

The metric space $(\mathcal F,\rho)$ is complete, because the metric spaces of the continuous on the compact set $[0,T]\times [0,\omega_*]$ functions and two-periodic, continuous on  $\mathbb {R}^2$ functions are complete. And \,$\mathcal F$ is a closed subset of the latter space.

Indeed, for any sequence $F_n\in\mathcal F,\ F_n\rightarrow F$ and for any pair of points $(t,\psi_1),(t,\psi_2)$ 
consider the numeric sequences $a_n = |F_n(t,\psi_1)-F_n(t,\psi_2)|,\ b_n = \mathcal L|\psi_1 - \psi_2|.$ 
Due to Lipschitz' condition, $a_n\le b_n$ for any $n.$ Then, according to the theorem on the passing to the limit in the evaluations, we obtain
\,$|F(t,\psi_1)-F(t,\psi_2)|\le \mathcal L|\psi_1 - \psi_2|.$
The boundness of the absolute value of the limit function by the same constant $\mathcal M$ is obvious.

\smallskip 
3) Assume $\psi^t=\psi(t,t_0,\psi^{t_0},\e)$ is a solution of the initial value problem of the equation 
\begin{equation}\label{sceq}
	\dot\psi = (\mathcal K_{\star\e}^\nu+\Psi_{\star}^\nu(t,\psi,F(t,\psi)\e^{1/2},\e)\e^{3/2})\e^\nu\qquad (F\in \mathcal F),
\end{equation}
i.\,e. \,$\displaystyle \psi^t=\psi^{t_0}+\int_{t_0}^t\! 
\big(\mathcal K_{\star\e}^\nu+\Psi_{\star}^\nu(s,\psi^s,F(s,\psi^s)\e^{1/2},\e)\e^{3/2}\big)\e^\nu ds$ \ $(t_0,\psi^{t_0}\in \mathbb{R}).$

Denote
$$h^t=|\psi_2^t - \psi_1^t|,$$
where $\psi_j^t$ is a solution of the equation \eqref{sceq} c $F=F_j\ \ (j=1,2).$\, Then 
$$\begin{matrix} 	\displaystyle h^t\le h^{t_0}+\bigg|\int_{t_0}^t\!\big|\Psi_{\star}^\nu(s,\psi_2^s,F_2(s,\psi_2^s)\e^{1/2},\e)-
	\Psi_{\star}^\nu(s,\psi_1^s,F_1(s,\psi_1^s)\e^{1/2},\e)\big|\e^{\nu + 3/2}ds\bigg|\buildrel (\ref{lip}_2) \over \le  \\ 
	\displaystyle	\le h^{t_0}+L\bigg|\int_{t_0}^t (h^s+|F_2(s,\psi_2^s)-F_1(s,\psi_1^s)|\e^{1/2})ds\bigg|\e^{\nu + 3/2}\le \\ 
	\displaystyle	\le h^{t_0}+L(1+\mathcal L\e^{1/2})\bigg|\int_{t_0}^t h^s\,ds\bigg|\e^{\nu + 3/2}+L\rho(F_1,F_2)|t-t_0|\e^{\nu + 2}, \end{matrix}$$
because
\begin{equation}\label{F2-F1} |F_2(s,\psi_2^s)-F_1(s,\psi_1^s)|\le \rho(F_1,F_2)+\mathcal Lh^s, \end{equation} 
because $|F_2(s,\psi_2^s)-F_1(s,\psi_1^s)|\le |F_2(s,\psi_2^s)-F_1(s,\psi_2^s)|+|F_1(s,\psi_2^s)-F_1(s,\psi_1^s)|.$

\smallskip 
By applying the reinforced Gronwall's lemma, we obtain inequality
\begin{equation} \label{psih} h^t\le h^{t_0}e^{L(1+\mathcal L\e^{1/2})|t-t_0|\e^{\nu + 3/2}}+
	\rho(F_1,F_2)(1+\mathcal L\e^{1/2})^{-1}(e^{L(1+\mathcal L\e^{1/2})|t-t_0|\e^{\nu + 3/2}}-1)\e^{1/2}. \end{equation}
	 
4) Consider the mapping $S_\e^\nu(F),$ dependant on the parameter $\e$ and set on the space $(\mathcal F,\rho),$ described with the following formula 
	\begin{equation} \label{SF} S_\e^\nu(F)(t_0,\psi^{t_0})=
		-\!\int\nolimits_{t_0}^{\infty} e^{-K_{\star}^\nu (s-t_0)\e^{\nu+1}}\e^{\nu+1} V_{\star}^\nu(s,\psi^s,F(s,\psi^s)\e^{1/2},\e)ds.\end{equation}
	
	Let us show, that, for any $\e\in (0,\e_*],$ the mapping $S_\e^\nu(F)\in \mathcal F,$ i.\,e. the mapping is from $\mathcal F$ to $\mathcal F,$
	and is uniformly-contracting.
	
a) The mapping $S_\e^\nu(F)$ is $T$-periodic in $t_0.$\, Indeed, 
$$\begin{matrix} \displaystyle S_\e^\nu(F)(t_0+T,\psi^{t_0})=
	-\!\int\nolimits_{t_0+T}^{\infty} e^{-K_{\star}^\nu (s-t_0-T)\e^{\nu+1}}\e^{\nu+1} V_{\star}^\nu(s,\tilde\psi^s,F(s,\tilde\psi^s)\e^{1/2},\e)ds= \\
	\displaystyle =-\!\int\nolimits_{t_0}^{\infty} 	e^{-K_{\star}^\nu (\tau-t_0)\e^{\nu+1}}\e^{\nu+1} 
	V_{\star}^\nu(\tau+T,\tilde\psi^{\tau+T},F(\tau+T,\tilde\psi^{\tau+T})\e^{1/2},\e)d\tau,	\end{matrix}$$ 
where, by defintition, $\tilde\psi^s=\psi(s,t_0+T,\tilde\psi^{t_0+T},\e),$ and $\tilde\psi^{t_0+T}=\psi^{t_0}.$

Due to the function $\Psi_{\star}^\nu$ being $T$-periodic in $t,$ the function \,$\tilde\psi^{t+T}=\psi(t+T,t_0+T,\psi^{t_0},\e)$ is a solution of the equation \eqref{sceq} 
with the same initial values $t_0,\psi^{t_0},$ as a solution $\psi^t=\psi(t,t_0,\psi^{t_0},\e),$ hence, these solutions are the same.
Now the required equality is implied from the functions $V_{\star}^\nu$ and $F$ being $T$-periodic for any $t_0\in \mathbb{R}^1.$ 

b) Let us show, that the mapping $S_\e^\nu(F)$ is $\omega_\star$-periodic in $\psi^{t_0}.$ 
Due to the functions $\Psi_{\star}^\nu$ and $F_\e$ being periodic, the function $\br\psi^t = \psi(t,t_0,\psi^{t_0},\e)+\omega_\star$ 
is a solution of the equation \eqref{sceq}, moreover, it has the same initial values $t_0,\psi^{t_0}+\omega_\star,$ 
as a solution $\check\psi^t=\psi(t,t_0,\psi^{t_0}+\omega_\star,\e),$ hence, they are the same. The rest of the proof is obvious.

c) \,$\displaystyle |S_\e^\nu(F)(t_0,\psi^{t_0})|\le M\e^{\nu+1}\int_{t_0}^{\infty}\! 
e^{-K_{\star}^\nu (s-t_0)\e^{\nu+1}}ds\le M(K_{\star}^\nu)^{-1}\buildrel (\ref{fineq})\over =\mathcal M.$

d) Let us show, that the function $S_\e^\nu(F)$ meets the Lipschitz' condition with the constant $\mathcal L,$ 
thus, ending the proof of $S_\e^\nu(F)$ belonging to $\mathcal F.$ 
Also we will prove, that the mapping $S_\e^\nu(F)$ is contracting on the $\mathcal F$ with the contraction rate equal to $1/2.$

For any \,$t_0,\psi_1^{t_0},\psi_2^{t_0}$\, and \,$F_1(t,\psi),F_2(t,\psi)\in \mathcal F$\, the following inequalities are correct:
$$\begin{matrix} |S_\e^\nu(F_2)(t_0,\psi_2^{t_0})-S_\e^\nu(F_1)(t_0,\psi_1^{t_0})|\buildrel (\ref{lip}_1) \over \le \\	
	\displaystyle \le \int_{t_0}^{\infty} e^{-K_{\star}^\nu (s-t_0)\e^{\nu+1}}\e^{\nu+1}  
	L\big(h^s+|F_2(s,\psi_2^s)-F_1(s,\psi_1^s)|\e^{1/2}\big)ds\! \buildrel \eqref{F2-F1} \over \le \quad \\
	\displaystyle \le \rho(F_1,F_2)L\int_{t_0}^{\infty}\! e^{-K_{\star}^\nu (s-t_0)\e^{\nu+1}}\e^{\nu+3/2} ds+
	L(1+\mathcal L\e^{1/2})\int_{t_0}^{\infty} e^{-K_{\star}^\nu (s-t_0)\e^{\nu+1}}\e^{\nu+1}h^s\,ds \buildrel \eqref{psih}\over \le \\	
	\displaystyle \le \big(\rho(F_1,F_2)L\e^{\nu+3/2}+L(1 + \mathcal L\e^{1/2})h^{t_0}\e^{\nu+1}\big)
	\int_{t_0}^{\infty} e^{(L(1+\mathcal L\e^{1/2})\e^{\nu+1/2}-K_{\star}^\nu)(s-t_0)\e^{\nu+1}}ds \buildrel (\ref{fineq}_2)\over \le \\	
	\le 2\big(\rho(F_1,F_2)L\e^{\nu+3/2}+L(1+\mathcal L\e^{1/2})h^{t_0}\e^{\nu+1}\big)/|K_{\star}^\nu\e^{\nu+1}|\!\buildrel (\ref{fineq}_2)\over \le     
	\rho(F_1,F_2)/2+\mathcal Lh^{t_0}. \end{matrix}$$

Choosing $F_2=F_1=F,$ we obtain the Lipschitz' inequation
$$|S_\e^\nu(F)(t_0,\psi_2^{t_0})-S_\e^\nu(F)(t_0,\psi_1^{t_0})|\le \mathcal Lh^{t_0}.$$ 

Choosing $\psi_2^{t_0}=\psi_1^{t_0}=\psi^{t_0}\ \ (h^{t_0}=0),$ we obtain the {\sl contracting} inequality 
$$|S_\e^\nu(F_2)(t_0,\psi^{t_0})-S_\e^\nu(F_1)(t_0,\psi^{t_0})|\le \max_{\,t,\psi\in \mathbb{R}}|F_2(t,\psi)-F_1(t,\psi)|/2.$$

According to the fixed-point principle, which can be found in  \cite[ch.\,9,\,\S\,7]{Zor} (contracting mapping principle) \, the function $S_\e^\nu(F)$ has a singular
fixed point $A_\e^\nu,$ i.\,e. $A_\e^\nu \in \mathcal F$ and $S_\e^\nu(A_\e^\nu)(t_0,\psi^{t_0})=A_\e^\nu(t_0,\psi^{t_0}).$ 

Therefore, for any \,$t_0,\psi^{t_0}\in \mathbb{R}$ we have:
\begin{equation} \label{is} 
	A_\e^\nu(t_0,\psi^{t_0})=-\!\int\nolimits_{t_0}^\infty 
	e^{-K_{\star}^\nu (s-t_0)\e^{\nu+1}}\e^{\nu+1}V_{\star}^\nu(s,\psi^s,A_\e^\nu(s,\psi^s)\e^{1/2},\e)ds, \end{equation} 
where, according to \eqref{sceq}, \,$\displaystyle \psi^s=\psi^{t_0}+
\int_{t_0}^s \!\big(\mathcal K_{\star\e}^\nu+\Psi_{\star}^\nu(\tau,\psi^\tau,A_\e^\nu(\tau,\psi^\tau)\e^{1/2},\e)\e^{3/2}\big)\e^\nu d\tau.$

5) Choose the arbitrary point $(t_{\text{\scriptsize A}},\psi_{\text{\scriptsize A}}^{t_A},z_{\text{\scriptsize A}}^{t_A}),$ located on the surface \,$z=A_\e^\nu(t,\psi),$ as an initial value for the system \eqref{scsys2}, 
i.\,e. \,$t_{\text{\scriptsize A}},\psi_{\text{\scriptsize A}}^{t_A}\in \mathbb{R},$ 
$z_{\text{\scriptsize A}}^{t_A}=A_\e^\nu(t_{\text{\scriptsize A}},\psi_{\text{\scriptsize A}}^{t_A}).$ 

Assume \,$\psi_{\text{\scriptsize A}}^t=\psi(t,t_{\text{\scriptsize A}},\psi_{\text{\scriptsize A}}^{t_A},\e)$ is a solution 
of the equation \eqref{sceq} with $F(t,\psi)=A_\e^\nu(t,\psi),$ 
i.\,e. \,$\displaystyle \psi_{\text{\scriptsize A}}^t=\psi_{\text{\scriptsize A}}^{t_A}+
\int_{t_{\text{\scriptsize A}}}^t\! \big(\mathcal K_{\star\e}^\nu+
\Psi_{\star}^\nu(s,\psi_{\text{\scriptsize A}}^s,A_\e^\nu(s,\psi_{\text{\scriptsize A}}^s)\e^{1/2},\e)\e^{3/2}\big)\e^\nu ds.$ Then the function 
$\psi_{\text{\scriptsize A}}^t$ satisfies the second equation of the system \eqref{scsys2} with $z=A_\e^\nu(s,\psi_{\text{\scriptsize A}}^s)\e^{1/2}.$

Consider the function
$$z_{\text{\scriptsize A}}^t=A_\e^\nu(t,\psi_{\text{\scriptsize A}}^t)\buildrel \eqref{is}\over =-\!\int\nolimits_t^\infty 
e^{-K_{\star}^\nu (s-t)\e^{\nu+1}}\e^{\nu+1}V_{\star}^\nu(s,\psi_{\text{\scriptsize A}}^s,A_\e^\nu(s,\psi_{\text{\scriptsize A}}^s)\e^{1/2},\e)ds.$$ 
Its graph is, obviously, located on the surface \,$z=A_\e^\nu(t,\psi).$ Additionally,  
$$\begin{matrix} \displaystyle \dot z_{\text{\scriptsize A}}^t=
	-K_{\star}^\nu\e^{\nu+1} \!\int\nolimits_t^\infty e^{-K_{\star}^\nu (s-t)\e^{\nu+1}}\e^{\nu+1}  
	V_{\star}^\nu(s,\psi_{\text{\scriptsize A}}^s,A_\e^\nu(s,\psi_{\text{\scriptsize A}}^s)\e^{1/2},\e)ds+ \\ 
	+V_{\star}^\nu(t,\psi_{\text{\scriptsize A}}^t,A_\e^\nu(t,\psi_{\text{\scriptsize A}}^t)\e^{1/2},\e)\e^{\nu+1}=
	(K_{\star}^\nu z_{\text{\scriptsize A}}^t+V_{\star}^\nu(t,\psi_{\text{\scriptsize A}}^t,z_{\text{\scriptsize A}}^t\e^{1/2},\e))\e^{\nu+1}. 
\end{matrix}$$ 

Therefore,
the vector-function $(z_{\text{\scriptsize A}}^t,\psi_{\text{\scriptsize A}}^t)$ is a solution of the initial value problem 
of the system \eqref{scsys2} with the initial values 
$t_{\text{\scriptsize A}},\psi_{\text{\scriptsize A}}^{t_A},A_\e^\nu(t_{\text{\scriptsize A}},\psi_{\text{\scriptsize A}}^{t_A}),$ 
hence, the surface \eqref{isep} is invariant for the system \eqref{scsys2}. \ $\Box$

\smallskip
Let us note, that, for $\nu=0,$ this theorem, as a theorem\;1, proven by Yu.\,N\,Bibikov in \cite[ch.\,1]{Bib3}, 
is the analogue of the Hale's lemma (see\,\cite{Hale}).

The three following lemmas describe the properties of the obtained surface.

\smallskip
{\bf Lemma 6.4.} \ 
{\it For any $\e\in(0,\e_*],$ $\tilde t,\tilde\psi^{\tilde t}\in\mathbb{R}$ and $\tilde z^{\tilde t}\!: |\tilde z^{\tilde t}|\le v_*\e^{-1/2},$ consider the functions 
$\tilde z^t =z(t,\tilde t,\tilde\psi^{\tilde t},\tilde z^{\tilde t},\e),\ \tilde\psi^t =\psi(t,\tilde t,\tilde\psi^{\tilde t},\tilde z^{\tilde t},\e),$ 
	which are the solutions of the initial value problem of the system \eqref{scsys2}. 
	Then 
	$$\begin{matrix}	\tilde z^t - A_\e^\nu(t,\tilde\psi^t) \xrightarrow[t\rightarrow -\infty]{} 0 \mbox{ when }K_{\star}^\nu > 0,\quad 
		\tilde z^t - A_\e^\nu(t,\tilde\psi^t) \xrightarrow[t\rightarrow +\infty]{} 0 \mbox{ when }K_{\star}^\nu < 0, \end{matrix}$$
	where $A_\e^\nu(t,\psi)$ is invariant surface of the system \eqref{scsys2} from theorem\;6.1. }

{\it Proof.}\, Assume $K_{\star}^\nu > 0.$  Consider the function
$$\tilde z_{\text{\scriptsize A}}^t=A_\e^\nu(t,\tilde\psi^t)\buildrel \eqref{is}\over =-\!\int\nolimits_t^\infty 
e^{-K_{\star}^\nu (s-t)\e^{\nu+1}}\e^{\nu+1}V_{\star}^\nu(s,\tilde \psi^s,A_\e^\nu(s,\tilde\psi^s)\e^{1/2},\e)ds.$$ 
Then \,$\dot{\tilde z}_{\text{\scriptsize A}}^t = 
(K_{\star}^\nu \tilde z_{\text{\scriptsize A}}^t+V_{\star}^\nu(t,\tilde\psi^t,\tilde z_{\text{\scriptsize A}}^t\e^{1/2},\e))\e^{\nu+1},$ 
and $\tilde z^t$ satisfies the first equation of the \eqref{scsys2}. 

Assuming \,$w^t=\tilde z^t-\tilde z_{\text{\scriptsize A}}^t$\, \,$(w^t\ne 0,$ if $\tilde z^{\tilde t}\ne A_\e^\nu(\tilde t,\tilde\psi^{\tilde t})),$ we obtain
$$\dot w^t = (K_{\star}^\nu w^t + V_{\star}^\nu(t,\tilde\psi^t,\tilde z^t \e^{1/2},\e) - 
V_{\star}^\nu(t,\tilde\psi^t, \tilde z_{\text{\scriptsize A}}^t \e^{1/2},\e))\e^{\nu+1} = W(t)w^t\e^{\nu+1},$$
additionally,
$|V_{\star}^\nu(t,\tilde\psi^t,\tilde z^t \e^{1/2},\e) - V_{\star}^\nu(t,\tilde\psi^t, \tilde z_{\text{\scriptsize A}}^t \e^{1/2},\e)|
\!\buildrel (\ref{lip}_2)\over \le L|w^t|\e^{1/2} \!\buildrel \eqref{fineq}\over\le K_{\star}^\nu|w^t|/ 4.$ 

Therefore \,$W(t) \ge 3K_{\star}^\nu/4$\, and
\,$\displaystyle \ln{w^t} = \ln{w^{\tilde t}} + \e^{\nu + 1}\int_{\tilde t}^t W(s)ds \xrightarrow[t\rightarrow -\infty]{} -\infty.$

Henceforth, \,$w^t = \tilde z^t - A_\e^\nu(t,\tilde\psi^t)\rightarrow 0$ when $t\rightarrow -\infty.$

If the constant $K_{\star}^\nu<0,$ then, substituting $t$ to $\tau=-t$ and performing the same actions, we obtain the equality for the $w^\tau.$ 
Substituting $\tau=-t$ in the equality, we obtain, that $w^t\rightarrow 0$ when $t\rightarrow +\infty.$ \ $\Box$

\smallskip
We will call
the invariant surface \eqref{isep} {\sl asymptotically stable} in a sense of previous lemma when $t\to -\infty,$ 
if the constant $K_{\star}^\nu>0$ in the nondegeneracy condition \eqref{ud}, or when $t\to +\infty,$ if $K_{\star}^\nu<0.$ 

\smallskip
{\bf Lemma 6.5.} \ 
{\it The function $A_\e^\nu(t,\psi)$ from \eqref{is} has the same degree of smoothness in $\psi,$ as a right-hand side of the system \eqref{scsys}. } 

The proof of this lemma can be found in \cite{Hale1}.  

Thus, in our case, $A_\e^\nu(t,\psi)\in C_{t,\psi}^{0,1}(\mathbb {R}^2).$

\smallskip
{\bf Lemma 6.6.}  
{\it The function $A_\e^\nu(t,\psi)$ from \eqref{is} is continuous in $\e$ when $\e\in (0,\e_*).$ }

{\it Proof.}\,  1) Let us show that the mapping $S_\e(F),$ introduced in \eqref{SF}, is continuous in $\e$ 
in complete metric space $(\mathcal F,\rho).$ 

Specifically, for any $\tilde \e\in(0,\e_*),$ by fixing sych constant $c>0,$ that 
the closed interval  $\overline P_c=[\tilde\e-c,\tilde\e+c]\subset (0,\e_*],$ it will be proved, that $S_{\e}(F)$ is continuous on the $\overline P_c.$ 

For any $F\in\mathcal F$ and any $t_0,\psi^{t_0}\in \mathbb{R}$ 
the mapping $\displaystyle S_\e(F)(t_0,\psi^{t_0})\buildrel \eqref{SF} \over =\int_{t_0}^{\infty} h(s,\e)ds,$ 
where \,$h=-e^{-K_{\star}^\nu (s-t_0)\e^{\nu+1}}\e^{\nu+1} V_{\star}^\nu(s,\psi^s,F(s,\psi^s)\e^{1/2},\e)$ is a 
continuous function on $R=\{(s,\e)\!:\,s\in [t_0,+\infty),\,\e\in \overline P_c\}.$ 
Indeed, $V_{\star}^\nu$ is continuous in totality of its arguments, and, according to the theorem on the integral continuity, 
the solution $\psi_\e^s=\psi(s,t_0,\psi^{t_0},\e)$ is continuous in a cylindrical neighbourhood of the graph of the solution 
$\psi_{\breve\e}^s,$ $s\in [t_0,s_0)$ with any $s_0>\breve s,$ therefore the function $\psi_\e^s=\psi(s,t_0,\psi^{t_0},\e)$ 
is continuous in a specific neighbourhood of an arbitrary point $(\breve s,\breve \e)\in R.$ 

\smallskip
Assume \,$H(s)=M(\tilde\e+c)^{\nu+1}e^{-K_{\star}^\nu(s-t_0)(\tilde\e-c)^{\nu+1}}\ \ (s\in [t_0,+\infty),\ M$ is from $(\ref{lip}_1)),$ 
then \,$h(s,\e)\le H(s)$\, for any \,$(s,\e)\in R$\, and
$\displaystyle \int_{t_0}^{\infty} H(s)ds=\frac{M(\tilde\e+c)^{\nu+1}}{K_{\star}^\nu(\tilde\e-c)^{\nu+1}}.$

According to the Weierstrass' majorant theorem (see \,\cite[ch.\,17,\,\S\,2,\,prop.\,2]{Zor}), the improper integral, 
that describes $S_{\e}(F)$ in \eqref{SF}, converges absolutely and uniformly on the closed interval $\overline P_c.$
According to the theorem on the continuity of the improper integral (see \,\cite[ch.\,17,\,\S\,2,\,prop.\,5] {Zor}), 
it is a function of the values of the continuous in $\overline P_c$ function $S_{\e}.$

Therefore, according to the proposition on the stability of the fixed point from \cite[ch.\,9,\,\S\,7]{Zor} 
the function  $A_{\e}^\nu$ is continuous in $\e$ when $\e\in(0,\e_*).$ \ $\Box$

\section{Theoretical results} 

\subsection{The existence of the two-dimensional invariant tori.\\}

\smallskip
{\bf Lemma 7.1.} \ {\it 
	For any $b_{kl}^{\star}$ and for any $\e\in (0,\e_*],$ the special polar system $(\ref{ps})$ has a two-dimensional  invariant surface 
	\begin{equation}\label{surr} \rho=\Upsilon_{\e}^\nu(t,\varphi)\quad (t,\varphi\in \mathbb{R}), \end{equation}
	which is paramterized by the continuous, $T$-periodic in $t,$ continuously differentiable and $\omega_{\star}$-periodic in $\varphi$ function 
	$$\Upsilon_{\e}^\nu=\alpha_{\star}^{-1}(\varphi)(Z_{\e}^\nu(t,\varphi)+\beta_{\star}(\varphi){Z_{\e}^\nu}^2(t,\varphi)),$$ 
	where 	$Z_{\e}^\nu=G_{\star}^\nu(t,\varphi,\e)\e+
	\big(A_\e^\nu(t,\varphi+\Omega_{\star}^\nu(t,\varphi,\e)\e)(1+H_{\star}^\nu(t,\varphi,\e)\e)+F_{\star}^\nu(t,\varphi,\e)\big)\e^2.$ 
	The function $\Upsilon_{\e}^\nu$ is bound and continuous in $\e,$ when $\e\in (0,\e_*),$ as a function of the argument $\e.$ It also converges uniformly with the respect to $t,\varphi$ to zero when $\e\to 0.$ }

\smallskip
{\it Proof.}\, The surface $r=Z_{\e}^\nu(t,\varphi)$ is obtained by substituting the surface  $(\ref{isep})$ into the composition of changes $(\ref{save}),(\ref{psi}),(\ref{scale}),(\ref{scale2}).$ Substituting just obtained surface into the change $(\ref{pave}),$ we obtain the surface \eqref{surr}.
Its boundedness, countinuousness and uniform convergence to zero is implied from the lemma\;6.6,
the boundedness of the functions $G_{\star}^\nu,H_{\star}^\nu,F_{\star}^\nu,$ included in the change \eqref{save} and the uniform boundedness of the function $A_\e^\nu,$ because $|A_\e^\nu(t,\varphi)|\le \mathcal M,$ where the constant $\mathcal M$ is introduced in \eqref{fineq}. \ $\Box$

\smallskip
{\bf Corollary.} \ {\it	The functions $\Upsilon_{\e}^0$ and $\Upsilon_{\e}^1$ from \eqref{surr} have the following asymptotic expansion: 
$$\begin{matrix} \Upsilon_{\e}^0=\alpha_{\star}^{-1}(\varphi)\big((\overline g_{\star}^0+\tilde g_{\star}^0(t,\varphi))\e+
	(A_\e^0(t,\varphi)+\beta_{\star}(\varphi)(\overline g_{\star}^0+\tilde g_{\star}^0(t,\varphi))^2+\tilde f_{\star}^0(t,\varphi))\e^2\big)+O(\e^3),\\ 
	\Upsilon_{\e}^1=\alpha_{\star}^{-1}(\varphi)\big((\overline g_{\star}^1+\hat g^1_{\star}(\varphi))\e+
	(A_\e^1(t,\varphi)+\beta_{\star}(\varphi)(\overline g_{\star}^1+\hat g^1_{\star}(\varphi))^2+\hat f_{\star}^1(\varphi))\e^2\big)+O(\e^3),\hfill
\end{matrix}$$  
	where $A_\e^\nu(t,\varphi)=O(1)$ in $\e.$ }

\smallskip
{\bf Theorem 7.1.} \ {\it Assume $b_{kl}^{\star}$ is an arbitrary admissable solution from the definition\;5.1, then:\, 
	$1)$ for any $\e\in (0,\e_*],$ the system $(\ref{sv})$ has two-dimensional invariant surface
	\begin{equation}\label{t1} T\!I\!S_{\e}^\nu=\{ (x,y,t)\!:\,x=x_\e(t,\varphi),\,y=y_\e(t,\varphi),\,t\in \mathbb{R}\ \ (\varphi\in \mathbb{R}) \},
	\end{equation}
	where \,$x_\e=C(\varphi)+(C(\varphi)-k)\Upsilon_{\e}^\nu(t,\varphi),\,y_\e=S(\varphi)+(S(\varphi)-l\,)\Upsilon_{\e}^\nu(t,\varphi)$ are 
	continuous, $T$-periodic in $t,$ continuously differentiable and $\omega_{\star}$-peridoic in $\varphi$ functions, 
	where $(C(\varphi),S(\varphi))$ is a solution of the system $(\ref{cs})$ with the initial values $C(0)=b_{kl}^{\star},\,S(0)=l,$ 
	and function $\Upsilon_{\e}^\nu$ is from $\eqref{surr};$\, 
	
	$2)$ $T\!I\!S_{\e}^\nu,$ which is homeomorphic to two-dimensional torus, if the time $t$ is factored by the period $T,$ is  {\sl asymptotically robust} in the sense of lemma\;6.4 when $t\to -\infty,$
	if $K_{\star}^\nu>0$ in \eqref{ud}, and when $t\to +\infty,$ if $K_{\star}^\nu<0;$\,
	
	$3)$ The continuous in $\e,$ when $\e\in (0,\e_*),$ functions $x_\e(t,\varphi)$ and $y_\e(t,\varphi)$ uniformly with the respect to $t,\varphi$ converge respectively to $C(\varphi)$ and $S(\varphi),$ when $\e\to 0,$ i.\,e. for sufficiently small $\e$ the surface $T\!I\!S_{\e}^\nu$ is located in infinitely small neighbourhood of the cylindrical surface $C\!I\!S_*$ from the definition \;$5.2.$ }

{\it Proof.}\, Torus-like invariant surface $(\ref{t1})$ is obtained by the substituting of the function \eqref{surr} into the change \eqref{pz}.
Its "uniform convergence"\ to $C\!I\!S_{\star}$ of the autonomous system \eqref{snv} is implied from the uniform convergence to zero of the function 
$\Upsilon_{\e}^\nu(t,\varphi),$ determined on lemma\;7.1. \ $\Box$

\subsection{The results in the autonomous case.}

Let us adapt the obtained results to the important particular case when the system \eqref{sv} is autonomous. The possible change of time variable allows to assume that $\nu=0.$ Moreover, the independence from $t$ of the perturbations allows to consider only those limitations, required in case $\nu=1.$

Thus, the system \eqref{sv} can be written in the following form
$$\dot x=\gamma(y^3-y)+\e X^a(x,y,\e),\ \ \dot y=-(x^3-x)+\e Y^a(x,y,\e),\eqno(\ref{sv}^a)$$  
where $X^a=X_0^a(x,y)+X_1^a(x,y,\e)\e,$ $Y^a=Y_0^a(x,y)+Y_1^a(x,y,\e)\e,$ 
$X_0^a,Y_0^a$ are real-analytic in $x,y$ function on $D_{\sigma,\sigma}^{x,y}=\{(x,y)\!:\,|x|,|y|<\sigma\}$\, 
when $\sigma>\sqrt{1+\gamma^{-1/2}}\ge \sqrt{2};$\,   
the continuous functions $X_1^a,Y_1^a\in C^1(G_{\sigma,\sigma,\e_0}^{x,y,\e}),$ 
\,$G_{\sigma,\sigma,\e_0}^{x,y,\e}=\{(x,y,\e)\!:\,x,y\in \mathbb{R},$ $|x|,|y|\le \sigma,\,\e\in [0,\e_0]\}.$ 

Let us consequently describe the changes in the denotions, definitions, formulae and results, related to the fact that the initial system is autonomous.

The superindex $\nu$ is not required now, as well as lemma\;5.1, The numeration of the formulae, analogous to the ones already obtained, now uses the symbol \,$a\ \ (autonomous)$ 
and for any \,$\omega$-periodic functions $\eta(\varphi)$ the extension \,$\eta(\varphi)=\overline\eta+\hat\eta(\varphi)$ is used, 
where \,$\displaystyle \overline\eta=\int_0^\omega \eta(\varphi)\,d\varphi$ is an average value of the function $\eta.$ 
Next,  
$$R_{kl}^{a\circ}=\overline {R_{kl}^{a\circ}}+\widehat {R_{kl}^{a\circ}}(\varphi)=
C'(\varphi)Y_0^a(C(\varphi),S(\varphi))-S'(\varphi)X_0^a(C(\varphi),S(\varphi)),\eqno(\ref{rio}^a)$$
where \,$(C(\varphi),S(\varphi))$ is the real-analytic $\omega_{kl}$-periodic solution of the initial value problem 
of the system $(\ref{cs})$ with the initial values $C(0)=b_{kl},\,S(0)=l,$ 
the function $R_{kl}^{a\circ}=R_{kl}(\varphi,0,0),$ and $R_{kl}=R_{kl}(\varphi,r,\e)$ is introduced in \eqref{ps5}; 
$$K_{kl}^a=\overline{(\!{R_{kl}}'_r)^{\!\circ}}-\overline{R_{kl}^{\circ} q_{kl}}\ne 0,\eqno (\ref{ud}^a)$$
and $(\!{R_{kl}}'_r)^{\!\circ}={R_{kl}}'_r|_{(\varphi,0,0)}$ is introduced in \eqref{pf2}, $q_{kl}(\varphi)$ is introduced \eqref{ps}.

\smallskip
{\bf Definition ${\bf 5.1^a.}$} \ 
{\it The solution of the generating (bifurcation) equation  
	$$\overline {R_{kl}^{a\circ}}(b_{kl})=\frac{1}{\omega_{kl}} \int_0^{\omega_{kl}} R_{kl}^{a\circ}(\varphi)\,d\varphi=0\eqno(\ref{pu}^a)$$
	is admissable for the system $(\ref{sv}^a)$ and is denoted as $b_{kl}^{\star},$ 
	if it satisfies the conditions \eqref{nd} and \eqref{gamma*} and the nondegeneracy condition $(\ref{ud}^a)$ holds when it is fixed.}   

\smallskip 
Let us fix the arbitrary admissable solution $b_{kl}^{\star}.$ All subindexes $kl$ in all denotions, except $b,$ are changed to $\star$ for brevity, 
for example, $K_{kl}^a=K_{\star}.$ 
Then in lemma\;6.1 the change \eqref{save} can be written as
$$r=u+g_{\star}(\varphi)\e+\hat h_{\star}(\varphi)u\e+\hat f_{\star}(\varphi)\e^2,\eqno(\ref{save}^a)$$
where \,$\displaystyle \hat g_{\star}(\varphi)=\int_{\varphi_g}^\varphi R_{\star}^{\circ}(s)ds,$ \ 
$\overline g_{\star}=\big(\overline{R_{\star}^{\circ}(\Phi_{\star}^{\circ}+\hat g_{\star}q_{\star}\big)}-
\overline{(\!{R_{\star}}'_\e)^{\!\circ}}-\overline{\hat g_{\star}(\!{R_{\star}}'_r)^{\!\circ}})/K_{\star},$ \ 
$\displaystyle \hat h_{\star}(\varphi)=\int_{\varphi_h}^\varphi \big((\!{R_{\star}}'_r)^{\!\circ}-R_{\star}^{\circ} q_{\star}-K_{\star}\big)ds,$ \ 
$\displaystyle \hat f_{\star}(\varphi)=\int_{\varphi_f}^\varphi 
\big(\overline g_{\star}K_{\star}-R_{\star}^{\circ}(\Phi_{\star}^{\circ}+\hat g_{\star}q_{\star})+
(\!{R_{\star}}'_\e)^{\!\circ}+\hat g_{\star}(\!{R_{\star}}'_r)^{\!\circ}\big)ds,$ \ 
and \,$(\!{R_{\star}}'_\e)^{\!\circ}=C'Y_1^a(C,S)-S'X_1^a(C,S),$ according to \eqref{pf2}. 

After that, in system \eqref{savesys} the function
\,$\Theta_{\star}(\varphi)=\Phi_{\star}^{\circ}(\varphi)+g_{\star}q_{\star}(\varphi).$
Formulate the proposition that includes theorem\;6.1 and lemmas 6.5,\,6.6.

\smallskip
{\bf Theorem ${\bf 6.1^a.}$} \ {\it 
	For any $\e\in (0,\e_*],$ the autonomous system \eqref{scsys2} with $\nu=0$ has two dimensional cylindrical invariant surface
	$$z=\mathcal A_\e(t,\psi)\quad (t,\psi\in \mathbb{R})\eqno (\ref{isep}^a)$$  
	where $\mathcal A_\e(t,\psi)\equiv A_\e(\psi),$ the function $A_\e(\psi)$ is continuous, $\omega_{\star}$-periodic, 
	satisfies the Lipschitz condtion with the global constant $\mathcal L,$ 
	and \,$\max |A_\e^\nu(\psi)|\le \mathcal M,$ the constants $\mathcal M$ and $\mathcal L$ are introduced in $(\ref{fineq}).$ } 

\smallskip
The proof is analogous to the proof of the theorem\;6.1, in which the metric function space $F(\psi)$ is considered as $\mathcal F.$

As a result the graph of the periodic function $z=A_\e^\nu(\psi)$ is a generatrix of the cylindrical surface $(\ref{isep}^a).$

\smallskip
{\bf Theorem ${\bf 7.1^a.}$} \ {\it Assume $b_{kl}^{\star}$ is an arbitrary admissable solution from the definition\,$5.1^a,$ then:\, 
	
	$1)$ for any $\e\in (0,\e_*],$ the system $(\ref{sv}^a)$ has an invariant cycle. usually called limit cycle,
	\begin{equation}\label{lc} LC_{\e}=\{ (x,y)\!:\,x=x_\e(\varphi),\ y=y_\e(\varphi) \}\quad (\varphi\in \mathbb{R}) , \end{equation}
	where the functions \,$x_\e=C(\varphi)+(C(\varphi)-k)\Upsilon_{\e}^\nu(\varphi),\,y_\e=S(\varphi)+(S(\varphi)-l\,)\Upsilon_{\e}^\nu(\varphi)$ are 
	continuously differentiable	and $\omega_{\star}$-periodic in $\varphi$; 
	
	$2)$ The cycle $LC_{\e}$ is {\sl asymptotically stable} when $t\to -\infty,$ if the constant $K_{\star}>0$ in the $(\ref{ud}^a),$
	and, when $t\to +\infty,$ 	if the constant $K_{\star}<0;$\,  
	
	$3)$ The continuous in $\e$ for any $\e\in (0,\e_*)$ functions $x_\e(\varphi),\,y_\e(\varphi)$ uniformly with the respect to $t,\varphi$ converge to the \,$C(\varphi),S(\varphi)$ respectively, when $\e\to 0,$ i.\,e. for small $\e$ the limit cycle $LC_{\e}$  is located in infinitely small neighbourhood of the cycle $GC_{\star}$ from the definition\;$5.2.$ } 

\section{Practical results}

\subsection{The analysis of the generating equation.} \ 
It is assumed that in the system \eqref{sv} the functions $X_0^\nu,Y_0^\nu$ can be written as
\begin{equation} \label{xya} 
	X_0^\nu(t,x,y)=\sum_{m,n=0}^\infty X_0^{(m,n)}(t)x^my^n,\ \ Y_0^\nu(t,x,y)=\sum_{m,n=0}^\infty Y_0^{(m,n)}(t)x^my^n 
\end{equation}
where the power series' with the real, continuous and $T$-periodic in $t$ coefficients converge absolutely, when $|x|,|y|<\sigma.$
Therefore, for any admissable solution $b_{kl},$ that satisfies the conditions \eqref{nd} and \eqref{gamma*}, in \eqref{rio}
$$R_{kl}^{\nu \circ}(t,\varphi)=\sum_{m,n=0}^\infty \left(Y_0^{(m,n)}(t)C'(\varphi)-X_0^{(m,n)}(t)S'(\varphi)\right)C^m(\varphi)S^n(\varphi).$$

The left-hand side of the generating equation \eqref{pu} can be written as
$$\overline {R_{kl}^{\nu \circ}}(b_{kl})={1\over \omega_{kl}}\sum_{m=0}^\infty\sum_{n=0}^\infty
\left(\overline{Y_0^{(m,n)}}\int_0^{\omega_{kl}} C^mS^nC'\,d\varphi-
\overline{X_0^{(m,n)}}\int_0^{\omega_{kl}} C^mS^nS'\,d\varphi\right),$$
and, when $n=0,$ the first integral is equal to zero, when $m=0,$ --- the second integral. 

Integrating the equality $(C^mS^{n+1})'=mC^{m-1}S^{n+1}C'+(n+1)C^mS^nS'$ over the period and collecting the terms, we obtain:
$$\overline {R_{kl}^{\nu \circ}}(b_{kl})={1\over T\omega_{kl}}\sum_{m=0}^\infty\sum_{n=1}^\infty
\!\left({m+1\over n}\overline{X_0^{(m+1,n-1)}}+\overline{Y_0^{(m,n)}}\right)I_{kl}^{mn},$$	
where $\displaystyle I_{kl}^{mn}=\int_0^{\omega_{kl}} C^m(\varphi)S^n(\varphi)C'(\varphi)\,d\varphi,$\,
moreover, for each class this formula can be simplified:
\begin{equation}\label{rfinal} \begin{matrix}
		\displaystyle \overline {R_{00}^{\nu \circ}}(b_{00})={4\over \omega_{00}}\sum_{m,n=0}^\infty P^{(2m,2n+1)}J_{00}^{mn},\ \ 
		\displaystyle \overline {R_{k0}^{\nu \circ}}(b_{k0})={2\over \omega_{k0}}\sum_{m,n=0}^\infty P^{(m,2n+1)}J_{k0}^{mn},\hfill 	\\
		\displaystyle \overline {R_{kl}^{\nu \circ}}(b_{kl})={1\over \omega_{kl}}\sum_{m,n=0}^\infty P^{(m,n+1)}J_{kl}^{mn},\hfill 
		P^{(m,n)}={m+1\over n}\overline{X^{(m+1,n-1)}}+\overline{Y^{(m,n)}}, \end{matrix} \end{equation}
where 
$$\begin{matrix}
	\displaystyle J_{00}^{mn}=\bigg\{ \int_0^{ r_0^{1e}} \varsigma^{2m}(S_+(\varsigma^2))^{2n+1}\,d\varsigma+
	\int_{r_0^{1e}}^{b_{00}^e} \varsigma^{2m}(S_-(\varsigma^2))^{2n+1}\,d\varsigma,\hbox{ if }\ b_{00}^e\in (r^e,r^\sigma),\\
	\displaystyle \int_{b_{00}^i}^0 \varsigma^{2m}(S_-(\varsigma^2))^{2n+1}d\varsigma,\hbox{ if }\ b_{00}^i\in (0,r^i)\bigg\},\\ 
	\displaystyle kJ_{k0}^{mn}=\int_{k l_1^0}^{k l_1^1}\varsigma^m(S_-(\varsigma^2))^{2n+1}d\varsigma+
	\int_{k l_1^1}^{k r_1^1}\varsigma^m(S_+(\varsigma^2))^{2n+1}\,d\varsigma+\int_{k r_1^1}^{b_{k0}}\varsigma^m(S_-(\varsigma^2))^{2n+1}\,d\varsigma,\\
	\displaystyle kJ_{kl}^{mn}=l^n \int_{k l_2^1}^{b_{kl}}\varsigma^m((S_+(\varsigma^2))^{n+1}-(S_-(\varsigma^2))^{n+1})\,d\varsigma\ \ (k,l=\pm 1); \quad S_\pm(\varsigma^2)\hbox{ from }\eqref{Spm}.\end{matrix}$$

Let us, for example, deduce the first formula in $(\ref{rfinal}),$ taking into account that the constants used in $(\ref{rfinal})$ are introduced in \eqref{ndb}.

In class $0^e],$ when $b_{00}^e\in (r^e,r^\sigma),$ the motion along the cycle is clockwise, and, when $|C(\varphi)|=r_{0e}^1,$ the function $S^2(\varphi)-1$ changes its sign in every quadrant.
Therefore
$$\begin{matrix}
	\displaystyle I_{00}^{mn}=\int_{b_{00}^e}^{r_0^{1e}}\varsigma^m(-S_-(\varsigma^2))^n\,d\varsigma+
	\int_{ r_0^{1e}}^0 \varsigma^m(-S_+(\varsigma^2))^n\,d\varsigma+\int_0^{-r_0^{1e}}\varsigma^m(-S_+(\varsigma^2))^n\,d\varsigma+ \\
	\displaystyle +\int_{-r_0^{1e}}^{-b_{00}^e}\varsigma^m(-S_-(\varsigma^2))^n\,d\varsigma+ 
	\int_{-b_{00}^e}^{-r_0^{1e}}\varsigma^m( S_-(\varsigma^2))^n\,d\varsigma+\int_{-r_0^{1e}}^0 \varsigma^m( S_+(\varsigma^2))^n\,d\varsigma+I_e,\quad  
\end{matrix}$$
where $\displaystyle I_e=
\int_0^{ r_0^{1e}}\varsigma^m( S_+(\varsigma^2))^n\,d\varsigma+\int_{ r_0^{1e}}^{ b_{00}^e}\varsigma^m( S_-(\varsigma^2))^n\,d\varsigma.$

\smallskip  
It is obvious, that $I_{00}^{mn}=4I_e,$ if $m$ is even, and $n$ is odd. For all other cases $I_{00}^{mn}=0.$ After the reindexing, the formula takes the required form.

In class $0^i],$ when $b_{00}^i\in (0,r^i),$ the motion along the cycle os counter-clockwise. Therefore \,$\displaystyle I_{00}^{mn}=
I_i+\int_0^{-b_{00}^i}\varsigma^m( S_-(\varsigma^2))^n\,d\varsigma+\int_{-b_{00}^i}^0\varsigma^m(-S_-(\varsigma^2))^n\,d\varsigma+
+\int_0^{b_{00}^i} \varsigma^m(-S_-(\varsigma^2))^n\,d\varsigma,$\, 
where $\displaystyle I_i=\int_{b_{00}^i}^0 \varsigma^m( S_-(\varsigma^2))^n\,d\varsigma.$
Next, $I_{00}^{mn}=4I_i,$ if the constant $m$ is even, and constant $n$ is odd. For all other cases $I_{00}^{mn}=0.$

\medskip
\subsection{The application of the obtained results.} \
Let us give an example of the system \eqref{sv}, in which the analytic extensions \eqref{xya} of the functions $X_0^\nu,\ Y_0^\nu$ do not contain terms of the power greater than three with zero average values,
and for any sufficiently small $\e>0,$ the system has an eleven two dimensional invariant tori \eqref{t1}.

\smallskip
{\bf Lemma 8.1.} \ {\it 
	Consider the system \eqref{sv} with \begin{equation}\label{coeff}	\begin{matrix} 
			\gamma=1/2,\ T=2\pi,\ \sigma> \sqrt{6};\\
			\nu=1\!:\ \ \overline {X_0^{(m,n)}}=0,\ \overline {Y_0^{(m,n)}}=0\ \ (m,n\in \mathbb{Z}_+),\ \hbox{ except}\hfill \\ 
			\quad  \overline {Y_0^{(0,1)}}=\tau_0=-3.314,\ \overline {Y_0^{(0,3)}}=\tau_1=-0.361,\ \overline {Y_0^{(2,1)}}=\tau_2=4.493; \\
			\nu=0\!:\ \ X_0^0=\cos t,\ \ Y_0^0=\tau_0(y + \sin t)+\tau_1y^3 + \tau_2(x^2y + \sin t).\hfill \end{matrix} \end{equation} 
	where $\overline{X_0^{(m,n)}},\overline{Y_0^{(m,n)}}$\, are the average values of the coefficients from extension \eqref{xya}. 
	Then the generating equation \eqref{pu} has an eleven solutions, denoted as $b_{kl}^j$
	and located in the following intervals: \,$b_{00}^0\in (1.795,1.815),$ $kb_{k0}^1\in (1.118, 1.148),$ $kb_{k1}^2,kb_{k,-1}^2\in (1.266,1.276),$ 
	$kb_{k1}^3,kb_{k,-1}^3\in (1.299,1.303)$ $(k=\pm 1).$ These solutions satisfy the conditions \eqref{nd} and \eqref{gamma*}.  
	Moreover, the periods $\omega_{kl}^j$ of the solutions of the system \eqref{cs} $CS_{kl}^j(\varphi)=(C(\varphi),S(\varphi))$ with the initial values $C(0)=b_{kl}^j,$ $S(0)=l$
	are located in the following intervals: \,$\omega_{00}^0\in (3.81,3.94),$ $\omega_{k0}^1\in (13.25,14.05),$ $\omega_{kl}^2\in (5.99,6.29),$ 
	$\omega_{kl}^3\in (7.72,8.49)$ (see fig.\,$3.4$), eleven constants $K_{kl}^{\nu j}$ from the condition 
	\eqref{ud} take the following values: 
	$K_{00}^{00}\approx 2.78,$ $K_{k0}^{01}\approx -0.34,$ $K_{kl}^{02}\approx -0.29,$ $K_{kl}^{03}\approx 2.07;$ 
	$K_{00}^{10}\approx 2.79,$ $K_{k0}^{11}\approx -0.32,$ $K_{kl}^{12}\approx -0.05,$ $K_{kl}^{13}\approx  0.10$ \,$(k,l=\pm 1).$
}

\smallskip
{\it Proof.}\, In the considered system \eqref{sv} 
\,$r^i=2^{-1/4}(\sqrt 2-1)\approx 0.348,$ $r^e=r_\gamma=2^{-1/4}(\sqrt 2+1)\approx 1.306$ according to \eqref{exc}
and $r^\sigma=(1+\sqrt{12})^{1/2}>2.113,$ according to \eqref{nd}, $b_d(0.5)=1.15-0.16\cdot 0.3^{1/2}\approx 1.062,\ b_u(0.5)=1.15+0.28\cdot 0.3^{1/2}\approx 1.303$ according \eqref{gamma*}.
Therefore, the intervals for $b_{kl}$ with the approximately calculated bounds are these: 
for class 0] $b_{00}^i\in (0,0.348),\ b_{00}^e\in (1.306,2.113),$ for class 1] $kb_{k0}\in (1.062,1.303),$ for class 2] $kb_{kl}\in (1,1.306).$ 
Each of them contains the corresponding intervals from the theorem's hypothesis,
for each bound of such interval, the bounds for approximate values of the periods $\omega_{kl}^j$ are found.

In turn, the formulas \eqref{rfinal} can be written as:

$\overline {R^{\nu \circ}_{00}}(b_{00})=
4(\overline {Y^{(0,1)}}J_{00}^{00}+\overline {Y^{(0,3)}}J_{00}^{01}+\overline {Y^{(2,1)}}J_{00}^{10})/\omega_{00},$ 

$\overline {R^{\nu \circ}_{k0}}(b_{k0})=
2(\overline {Y^{(0,1)}}J_{k0}^{00}+\overline {Y^{(0,3)}}J_{k0}^{01}+\overline {Y^{(2,1)}}J_{k0}^{20})/\omega_{k0},$ 

$\overline {R^{\nu \circ}_{kl}}(b_{kl})=
{\phantom 2}\,(\overline {Y^{(0,1)}}J_{kl}^{00}+\overline {Y^{(0,3)}}J_{kl}^{02}+\overline {Y^{(2,1)}}J_{kl}^{20})/\omega_{kl}.$

Let us notice, that these three formulae cover all eleven cases, 
because in \eqref{coeff} we consider as nonzero only such coefficients, that
 for class 1] \,$J_{10}^{mn}=J_{-10}^{mn}$ for even $m$ and odd $n,$  
for class 2] \,$J_{11}^{mn}=J_{-11}^{mn}=J_{-1,-1}^{mn}=J_{1,-1}^{mn}$ for even $m$ and $n.$  

We obtain the following results:

for class 0]\, $\overline {R^{\nu \circ}_{00}}(1.795)<-10^{-4},\,\overline {R^{\nu \circ}_{00}}(1.815 )>10^{-4};$ 

for class 1]\, $\overline {R^{\nu \circ}_{k0}}(1.118\cdot k)>10^{-4},\,\overline {R^{\nu \circ}_{k0}}(1.148\cdot k)<-10^{-4};$

for class 2]\, $\overline {R^{\nu \circ}_{kl}}(1.266\cdot k)>10^{-4},\,\overline {R^{\nu \circ}_{kl}}(1.276\cdot k)<-10^{-4};$

for class 2]\, $\overline {R^{\nu \circ}_{kl}}(1.299\cdot k)<-10^{-4},\,\overline {R^{\nu \circ}_{kl}}(1.303\cdot k)>10^{-4},$ \\
and the function $\overline {R^{\nu \circ}_{kl}}$ doesn't change its sign otherwise. 
Therefore, there are eleven values of the parameter $b,$ contained in the intervals from theorem's hypothesis,  
and these parameters are the solutions of the generating equation \eqref{pu}. 

\smallskip 
For $\nu=1,$ for each obtained $b_{\star}=b_{kl}^j$ and the corresponding period $\omega_{\star}=\omega_{kl}^j,$ calculated in \eqref{ome}, 
the approximate value of the constant $K_{\star}^1 = K_{kl}^{1j}$ from the nondegeneracy condition \eqref{ud} can be calculated, using the following formula:
$$K_{\star}^1=\overline{(\!{R_{\star}^1}'_r)^{\circ}}-\overline{\widehat {R_{\star}^{1 \circ}}q_{\star}}=
\frac{1}{\omega_{\star}}\int_0^{\omega_{\star}} \mathcal R_{\star}^1(\varphi)\,d\varphi,$$
where \,$\mathcal R_{\star}^1=\alpha_{\star}^{-1}C'(2\tau_2 C(C-k)S+(\tau_1+\tau_2 C^2+3\tau_3 S^2)(S-l))-
(\tau_1 S+\tau_2 C^2 S+\tau_3 S^3)(\alpha_{\star}^{-2}\alpha'_{\star}(C-k)+2\beta_{\star}C')-q_{\star}C'(\tau_1 S+\tau_2 C^2 S+\tau_3 S^3).$ 

\smallskip 
For $\nu=0,$ for each obtained $b_{\star}=b_{kl}^j$ and the corresponding period $\omega_{\star}=\omega_{kl}^j,$
the function $\tilde g_{\star}^{0}(t,\varphi)=\tilde g_{kl}^{0j}(t,\varphi)$ from the equation \eqref{gn} is evaluated. 
To achieve that, we will switch to the sine-cosine form of the Fourier series' in formulas \eqref{zigeq},\eqref{lnu},\eqref{chitphi} from lemma\;5.1.
Consider 
$$\eta_1^{(n)}=\frac{2}{\omega_{\star}}\int_0^{\omega_{\star}} 
\widetilde R_{\star}^{0 \circ}(t,\varphi)\cos \frac{2\pi n\varphi}{\omega_{\star}}\,d\varphi,\ \ 
\eta_2^{(n)}=\frac{2}{\omega_{\star}}\int_0^{\omega_{\star}} 
\widetilde R_{\star}^{0 \circ}(t,\varphi)\sin \frac{2\pi n\varphi}{\omega_{\star}}\,d\varphi,$$
then
\begin{equation}\label{furi}
	\begin{matrix} 
		\displaystyle \widetilde R_{\star}^{0 \circ}=\sum_{n=1}^\infty \left(\eta_1^{(n)}(t)\cos \frac{2\pi n\varphi}{\omega_{\star}}+
		\eta_2^{(n)}(t)\sin \frac{2\pi n\varphi}{\omega_{\star}}\right),\\ 
		\displaystyle \tilde g_{\star}^0=\sum_{n=1}^\infty \left(\breve\chi_1^{(n)}(t)\cos \frac{2\pi n\varphi}{\omega_{\star}}+
		\breve\chi_2^{(n)}(t)\sin \frac{2\pi n\varphi}{\omega_{\star}}\right),	\end{matrix}
\end{equation}
where, for any $n,$ the vector $\breve\chi^{(n)}(t)$ \, of $T$-periodic coefficients $(T=2\pi)$ satisfies the linear system
\begin{equation}\label{lns} \dot \chi^{(n)}=A^{(n)}\chi^{(n)}+\eta^{(n)}(t)\quad (n\in \mathbb{N}), \end{equation}
where $A^{(n)}=\begin{pmatrix} 0&-\alpha_n\\ \alpha_n&0\end{pmatrix},$ $\alpha_n=\cfrac{2\pi n}{\omega_{\star}},$ 
i.\,e. $\lambda_{1,2}=\pm i\alpha_n$ are eigenvalues of the matrix $A^{(n)}.$ 

The system \eqref{lns} is a real analogue of the equation \eqref{lnu} and  
has a singular \,$T$-periodic real solution $\breve\chi^{(n)}$ 
if conditions $\lambda_{\iota}T/(2\pi i)\not\in \mathbb{Z}\ \ (\iota=1,2),$ 
or \,$nT-m\omega_{\star}\ne 0\ \ (m,n\in \mathbb{N})$ hold. They hold due to the Siegel's condition \eqref{zig}. 
This solution can be written as:  
$$\breve\chi^{(n)}(t)=(E-e^{A^{(n)}T})^{-1}\int_{t-T}^t e^{A^{(n)}(t-s)}\eta^{(n)}(s)\,ds\quad (T=2\pi),$$ 
where $e^{A^{(n)}t}=\begin{pmatrix} \cos \alpha_nt & -\sin \alpha_nt\\ \sin \alpha_nt & \cos \alpha_nt \end{pmatrix}$ --- 
normalized fundamental matrix of the system $\dot \chi^{(n)}=A^{(n)}\chi^{(n)},$
\,$(E-e^{A^{(n)}T})^{-1}=(2-2\cos \alpha_nT)^{-1} \!\begin{pmatrix} \varpi_n & -\ae_n \\ \ae_n & \varpi_n \end{pmatrix}\!,$
$\varpi_n=\cos \alpha_nt-\cos \alpha_n(T-t),\ \ae_n=\sin \alpha_nt+\sin \alpha_n(T-t).$ 
Additionally, for any $m\in \mathbb{N},$ $\cos \alpha_nT=\cos(2\pi \omega_{\star}^{-1}(nT-m\omega_{\star})).$

\smallskip
The approximate values of the function $\tilde g_{\star}^0$ are calculated in points $(l_\varphi h_\varphi,l_t h_t),$ 
where $0\le l_\varphi,l_t \le 101,\ h_\varphi=0.01\omega_{\star},\ h_t=0.02\pi,$ 
by considering the finite sums of fifteen terms in the Fourier extensions \eqref{furi}.

Two-periodicity of the function $\tilde g_{\star}^0$ is confirmed by the comparison of its approximate values at each bound of both intervals, which length is equal to one of the periods: 
$$|\tilde g_{\star}^0(l_\varphi h_\varphi,0) - \tilde g_{\star}^0(l_\varphi h_\varphi,2\pi)|, \
|\tilde g_{\star}^0(0,l_t h_t) - \tilde g_{\star}^0(\omega_{\star},l_t h_t)| < 0.011\quad (l_\varphi,l_t=\overline{0,100}).$$

According to lemma 5.1, the function $\tilde g_{\star}^0$ is continuously differentiable, therefore, the approximate values 
of the function $\tilde g_{kl}^{0\,'}$ at the same points as $\tilde g_{\star}^0$ can be calculated by using the following formula
$\tilde g_{kl}^{0\,'}(l_\varphi h_\varphi,l_th_t) = 
100\big(\tilde g_{\star}^0((l_\varphi+1)h_\varphi,l_th_t) - \tilde g_{\star}^0(l_\varphi h_\varphi,l_th_t)\big).$

Then the approximate value of the constant $K_{\star}^0 = K_{kl}^{0j}$ from the nondegeneracy condition \eqref{ud} can be calculated by using formula:
$$K_{\star}^0=\overline{(\!{R_{\star}^0}'_r)^{\circ}-\tilde g_{\star}^{0\,'}q_{\star}}= 
\frac{1}{2\pi\omega_{\star}}\int_0^{\omega_{\star}} \int_0^{2\pi} \mathcal R_{\star}^0(\varphi)\,d\varphi,$$
where $\mathcal R_{\star}^0=\alpha_{\star}^{-1}C'(2\tau_2C(C-k)S+(\tau_0+3\tau_1S^3+\tau_2C^2)(S-l))+\alpha_{\star}'\alpha_{\star}^{-2}((S-l)\cos t- 
(C-k)(\tau_o(S+\sin t) + \tau_1S^3 + \tau_2(C^2S + \sin t)) - 2\beta_{\star}(C'(\tau_o(S+\sin t) + \tau_1S^3 + \tau_2(C^2S + \sin t)-S'\cos t) - 
\tilde g_{\star}^{0\,'}q_{\star}.$ \ $\Box$
 
\includegraphics[scale=0.36]{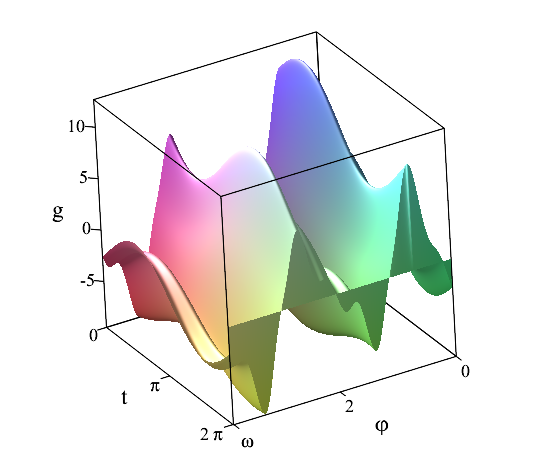}\qquad
\includegraphics[scale=0.36]{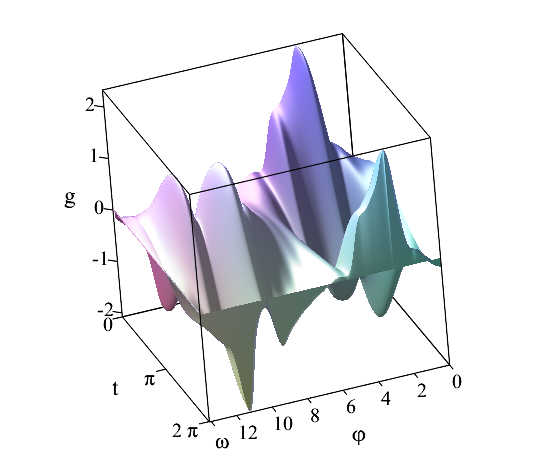}

\includegraphics[scale=0.36]{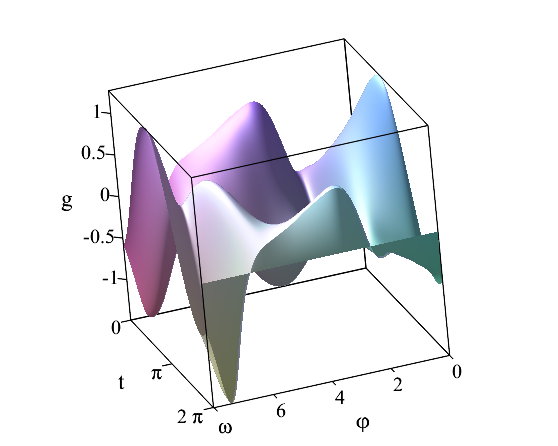}\qquad
\includegraphics[scale=0.36]{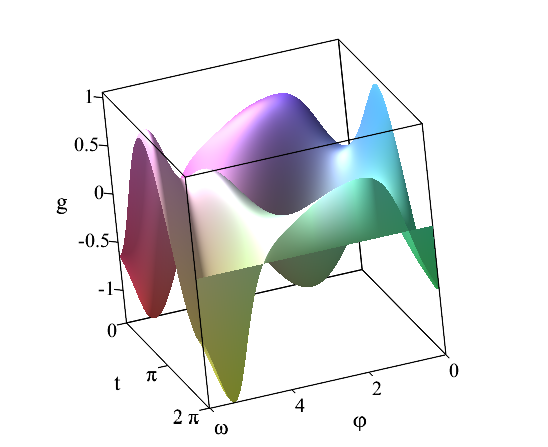}

{\small \bf Fig.\,8.1. 	The graphs of the two-periodic functions
	$\tilde g_{00}^{00}(t,\varphi),\tilde g_{10}^{01}(t,\varphi),\tilde g_{02}^{11}(t,\varphi),\tilde g_{11}^{03}(t,\varphi).$ }
 
 \smallskip
 Thus, in lemma it is established that, for any $\nu=1,$ each found $b_{kl}^j,$ by definition\;5.1, is admissable,   
 and for any $\nu=0,$ the solution $b_{kl}^j$ is admissable, if the periods $2\pi$ and $\omega_{kl}^j$ satisfy Siegel's condition
on the "small denominators"\ \eqref{zig}. Unfortunately, it is complicated to prove that Siegel's condition holds, because the accurate value of the period $\omega_{kl}^j$ is not known, but there is a following fact to be taken into account.

\smallskip
{\bf Proposition 8.1.} \ {\sl Siegel condition \eqref{zig}
	holds for almost every value of the period $\omega_{kl}^j$ with the respect to the Lebesgue measure. }

\smallskip
{\bf Remark 8.1.} \ {\sl Considering the proposition\;8.1, we will assume that Siegel's condition holds, therefore, for $\nu=0$ 
	each solution $b_{kl}^j$ is admissable. It is correct for almost every value of the coefficients $X_0^{(m,n)},Y_0^{(m,n)},$ 
	located in the infinitely small neighbourhoods of the coefficients chosen in \eqref{xya}.  }

\smallskip
{\bf Theorem 8.1.} {\it There is such constant $\e_*>0,$ that for any $\varepsilon\in (0,\e_*],$ the system \eqref{sv} from lemma\;8.1
	has an eleven invariant surfaces $T\!I\!S_{\e kl}^{\nu j},$ determined by the admissable solutions $b_{kl}^j$ and described in the theorem\;$7.1.$ }  

\smallskip
{\bf Corollary.} {\it Each surface $T\!I\!S_{\e kl}^{\nu j}$ takes form $\eqref{t1}$ and is homeomorphic to two-dimensional torus, 
	because it is parametrized by two-periodic functions $x=x_\e(t,\varphi),\,y=y_\e(t,\varphi),$ constructed from the solution $CS_{kl}^j(\varphi);$
	the surfaces $T\!I\!S_{\e 00}^{\nu 0},$ $T\!I\!S_{\e k,\pm1}^{\nu 3}$  are {\sl asymptotically stable} when $t\to -\infty,$ 
	and $T\!I\!S_{\e k0}^{\nu 1},$ $T\!I\!S_{\e k,\pm1}^{\nu 2}$ are {\sl asymptotically stable} when $t\to +\infty,$ in a sense of lemma\;$6.4$ 
	$(k=\pm 1);$ 
	the continuous in $\e$ functions $x_\e(t,\varphi)$ and $y_\e(t,\varphi)$ uniformly with the respect to $t,\varphi$ converge 
	respectively to $C(\varphi)$ and $S(\varphi),$ when $\e\to 0,$ hence, for any small $\e,$ each of the eleven surfaces $T\!I\!S_{\e kl}^{\nu j}$ is located in infinitely small neighbourhood of the cylindrical surface $C\!I\!S_{kl}^{\nu j}$ from the definition\;$5.2$ with the generating cycle, parameterized by the solution $CS_{kl}^j(\varphi),$ as generatrix \,(see fig.\,$8.2).$ }

\smallskip
Thus, the generating cycles of the system \eqref{sv} from the hypothesis of the theorem\;8.1 are illustrated on the figure 8.2, which shows the {\sl phase portrait} of the unperturbed system \eqref{snv} with $\gamma=1/2.$
These cycles pass respectively through the points $(x_0,0),$ $(\pm x_1,0),$ $(\pm x_j,1),$ $(\pm x_j,-1)$ $(j=2,3),$ 
where $x_0=b_{00}^0,$ $x_1=b_{10}^1,$ $x_j=b_{11}^j,$ and dash lines are used to denote separatrices.

\begin{center} \includegraphics[scale=0.64]{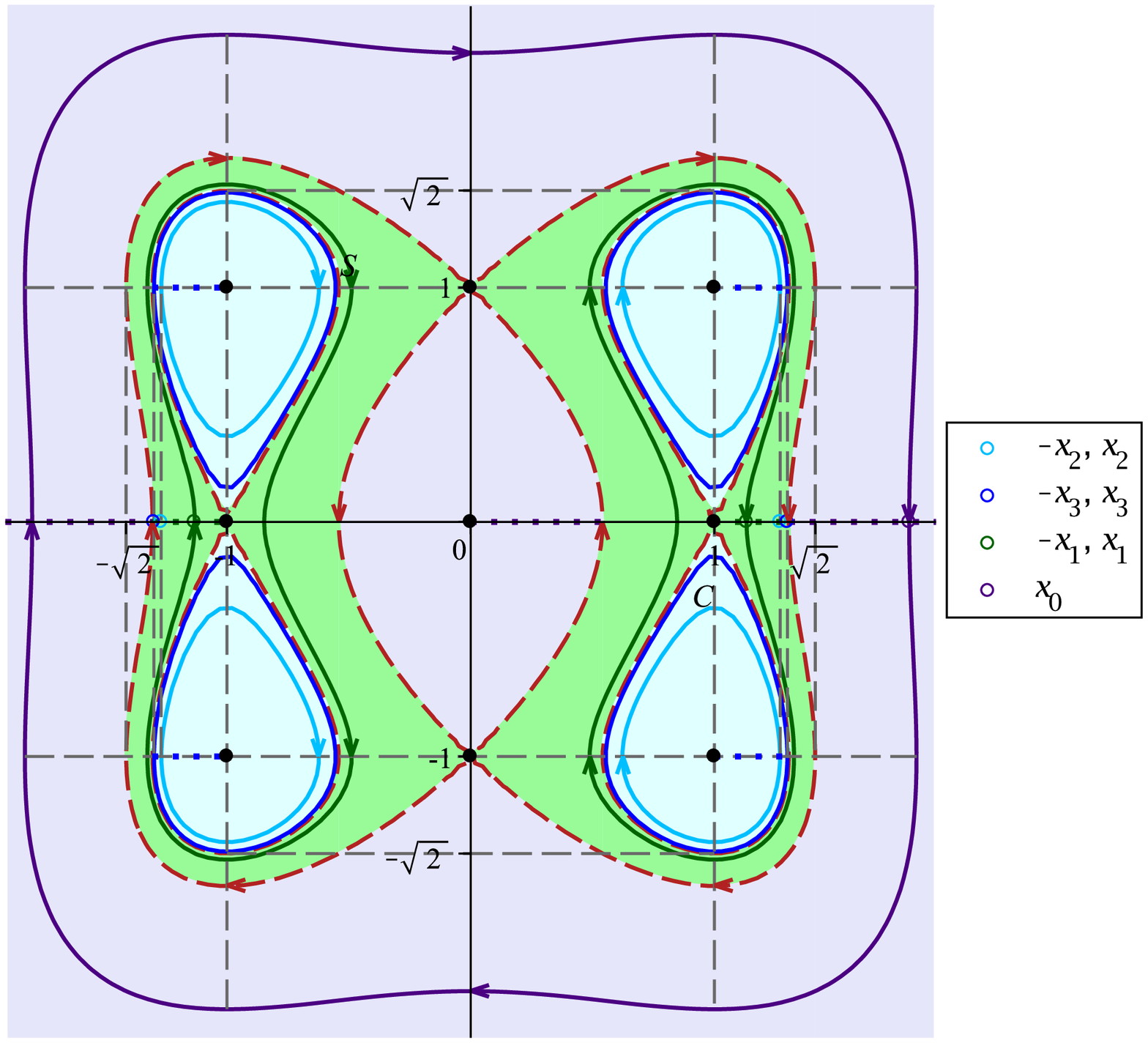} \end{center}
\begin{center} {\small \bf Fig.\,8.2. The generating cycles of the system (1.1) from theorem\;8.1. } \end{center}
	
\medskip
{\bf $8.3.$ The numeric confirmation of the obtain results.}
Theorem\;8.1 does not give the evaluation 
for the right bound $\e_*$ of the interval for the small parameter $\e,$ 
but it can be found numerically.

Consider only the autonomous case, because it allows to demonstrate the consequent appearance of all eleven limit cycles as the $\e_*$ decreases.

Thus, consider the system \eqref{sv} with the autonomous perturbation \eqref{coeff}
\begin{equation}\label{sa}
	\dot x=(y^3-y)/2,\ \ \dot y=-(x^3-x)+Y_0^a(x,y)\e\quad (\e\in[0,\e_*]), 
\end{equation}
where $Y_0^a=-3.314y+4.493x^2y-0.361y^3.$

The system \eqref{sa} is invariant with the respect to the change $x$ to $-\tilde x$ and $y$ to $-\tilde y.$ 
Therefore it is enough to locate the generating cycles from classes 1],\,2] only in the right semiplane, 
i.\,e. to choose parameters $b_{k0},b_{k,\pm 1}$\, with \,$k=1.$

First, let us notice that the abscissas of the equilibrium points of the system \eqref{sa} that don't lay on the abscissa axis have changed. 
Now the equilibrium points have the following coordinates:
\,$(\pm a_0^\e,\pm 1),$ $(\pm a_1^\e,\pm 1),$ $(\mp a_2^\e,\pm 1),$ 
where, for example, for $\e=0.05,$ $a_0^\e\approx 0.182,\,a_1^\e\approx 0.983,\,a_2^\e\approx 1.025.$ 

Next, the phase portrait of the orbits for the system \eqref{sa} with $\e=0.05$ is plotted (see fig.\,8.3).
The plot shows the limit cycles $LC_{\!\e 00}^{\, 0},$ $LC_{\!\e k0}^{\, 1}$ from the theorem\,$7.1^a,$ defined by the formula \eqref{lc} and located in the small neighbourhoods of the corresponding generating cycles $GC_{\star},$ related to classes 0],\,1]. These generating cycles are shown on the figure 8.1.

\begin{center} \includegraphics[scale=0.62]{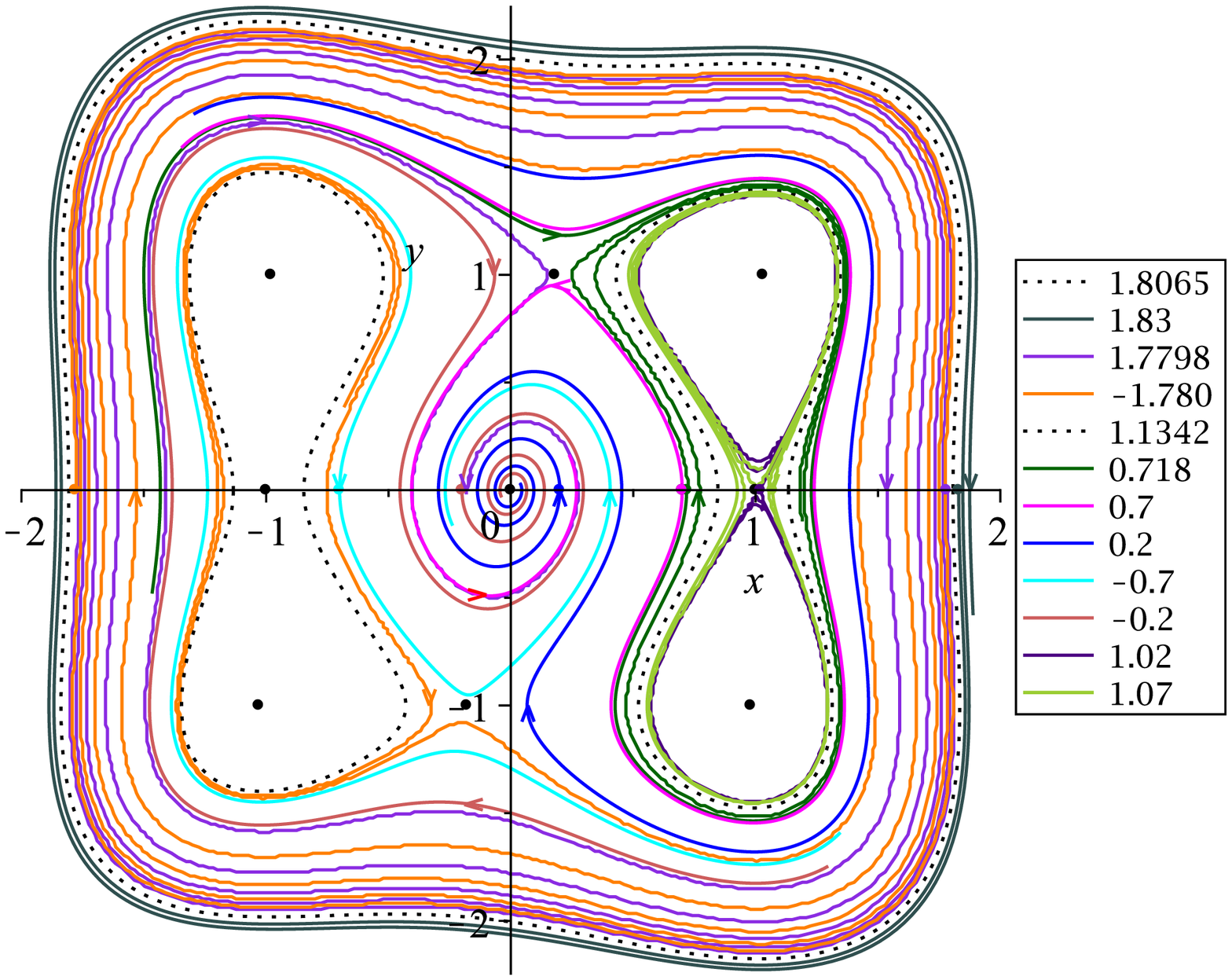} \end{center}
\begin{center} {\small \bf Fig.\,8.3. \ Three limit cycles of the system (8.4) with $\bf{\e=0.05.}$ } \end{center}

The behaviour of the separatrices, that adjoin the equilibrium points $(\pm a_0^\e,\pm 1),$ can be seen on the figure 8.3. 
One of the two separatrices, that "start"\ from the point $(a_0^\e,1),$ {\sl spirals onto} the equilibrium point $(0,0),$ as the time increases. 
The second separatrix spirals onto the right limit cycle $LC_{\!\e 10}^{\, 1},$ while both separatrices, that "end"\ at the point $(a_0^\e,1),$ 
{\sl spirals onto} the outer limit cycle $LC_{\!\e 00}^{\, 0}$ as the time decreases.

Unfortunately, for $\e=0.05$ the theorem\;8.1, that guarantees the existence of eight more limit cycles $LC_{\!\e kl}^{\, 2},$ $LC_{\!\e kl}^{\, 3},$ located in the small neighbourhoods of the corresponding generating cycles from class 2] (two in each quadrant), doesn't take into effect yet, and, for example, there are no other limit cycles inside the right limit cycle $LC_{\!\e 10}^{\, 1}.$ 
On top of the calculations proving it, it can be seen on the figure 8.4, where on top of the limit cycle $LC_{\!\e 10}^{\, 1}$ with $\e=0.05$ 
the generating cycles from class 2] are plotted. 
These generating cycles are parametrized by the solutions $CS_{1,\pm 1}^2(\varphi),CS_{1,\pm 1}^3(\varphi),$ which are independent of $\e,$ 
because the system \eqref{sa} is such, that all admissable solutions $b_{kl}^j$ are independent of $\e.$ 
\begin{center}\includegraphics[scale=0.62]{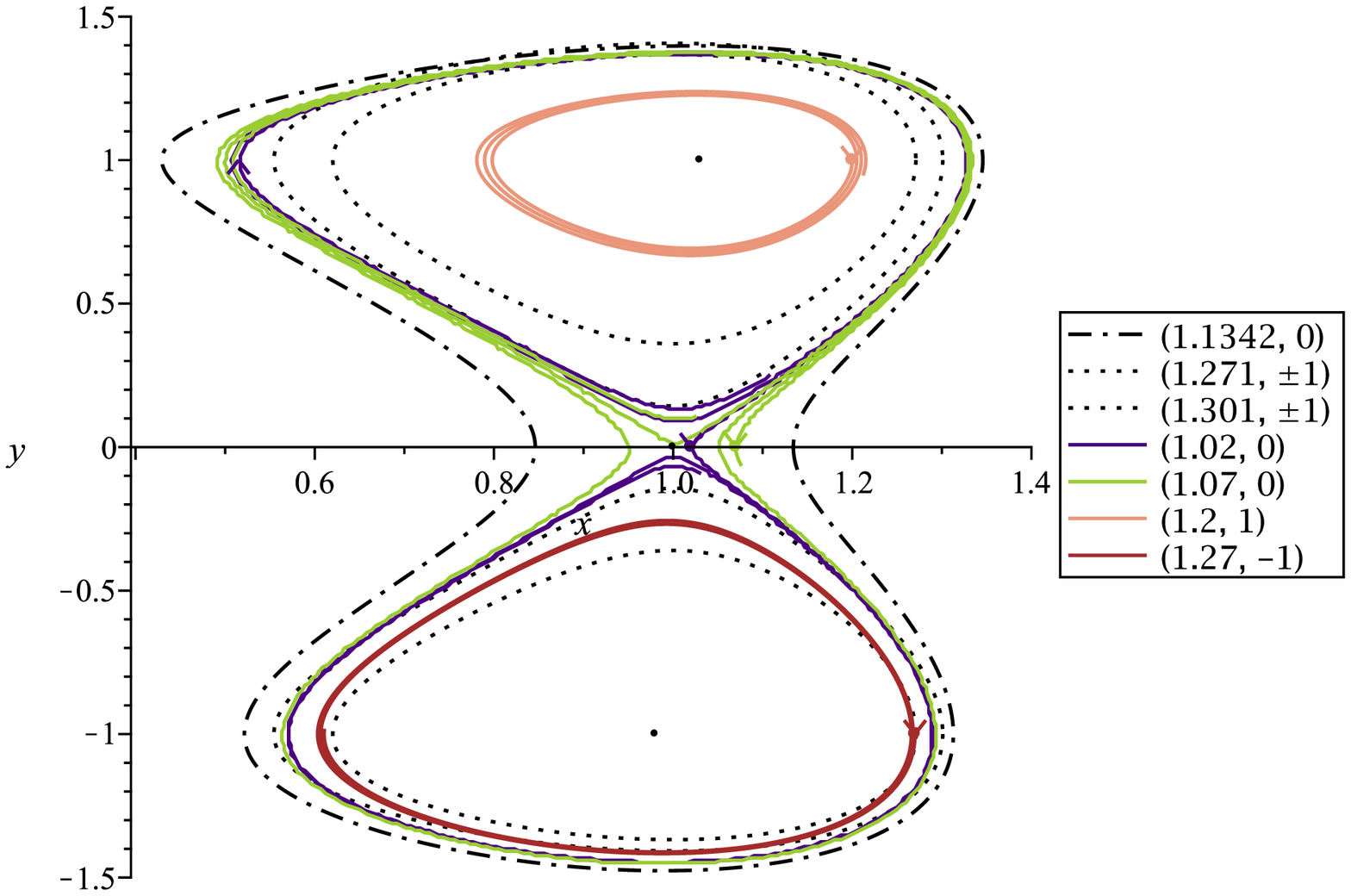} \end{center} 
\begin{center} {\small \bf Fig.\,8.4. The orbits inside of the right limit cycle with $\bf{\e=0.05.}$ } \end{center}
\begin{center} \includegraphics[scale=0.62]{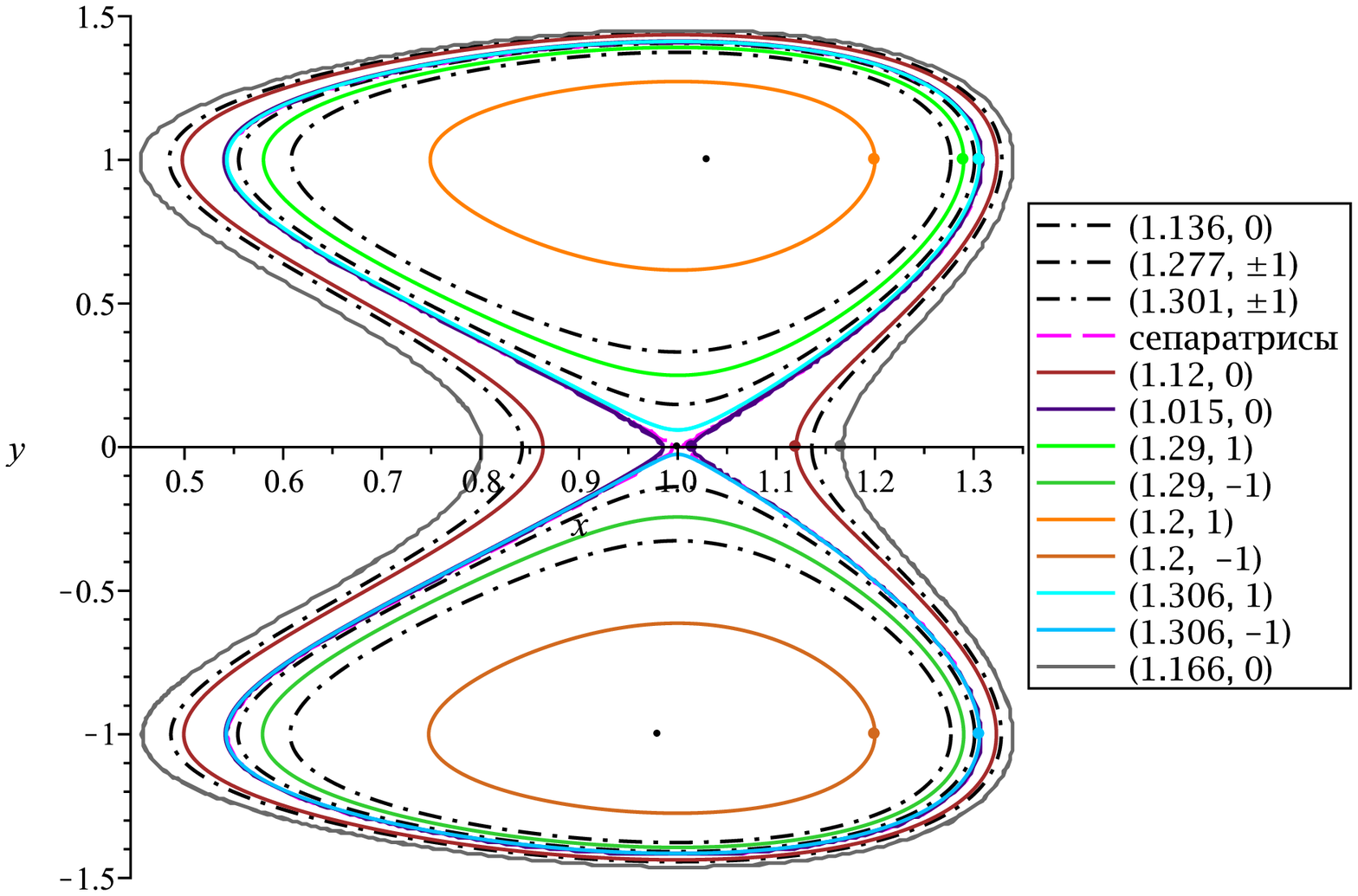} \end{center}
\begin{center} {\small \bf Fig.\,8.5. The orbits inside of the right limit cycle with $\bf{\e=0.001.}$ } \end{center}

As the value of the parameter $\e$ decreases, the limit cycle $LC_{\!\e 10}^{\, 1}$ changes its shape, and all four generating cycles from class 2] are caged inside of it along with the limit cycles, generated in their small neighbourhood (see fig.\;8.5).

The phase portrait 8.4 shows the behaviour of the orbits, which {\sl spiral onto} the equilibrium point $(a_0^\e,1)$ (unstable focus) as the time decreases, and some of them {\sl spiral onto} the equilibrium point $(a_0^\e,-1)$ (stable focus) as the time increases, other orbits --- it depends on the location relative to separatrices --- {\sl spiral onto} the stable limit cycle $LC_{\!\e 10}^{\, 1}.$ 

Let us decrease $\e.$ It is found out that, for  $\e=10^{-3}$ all eight limit cycles related to class 2], are generated, but the behaviour of the orbits (spirals) in their neighbourhoods can not be illustrated easily, because they appear on the phase portrait as cycles (see fig.\,8.5). 

The existence of the eleven limit cycles is implied from the table below. This table show the values of the abscissas of the intersection points of the spirals and the line $y\equiv 0$ in classes $0^e],1]$ or $y\equiv \pm 1$ in class~$2].$
These abscissa's values converge to each other as the time increases (decreases), which guarantees the existence of the stable (unstable) limit cycle between these spirals.  

Thus, for $\e=10^{-3}$ we obtain:

\medskip                  
{\small  
$\begin{matrix} 
	\phantom{aaaa} \text{0] Unstable limit cycle} \\ \phantom{aaaa} \text{ i.\,v. } (0,x_0^\e,0),\ \ x_0^\e\approx 1.808\ \end{matrix}  
\phantom{aaaaaaaaaaaa} 
\begin{matrix} 
	\text{1] Stable limit cycle} \\  \text{ i.\,v. } (0,x_1^\e,0),\ \ x_1^\e\approx 1.136\ \end{matrix}$
	
	\smallskip
	$\begin{matrix}        t& x(t)<x_0^\e\\
		\hfill 00.0& 1.804\hfill \\
		-03.9& 1.8040434\\
		-07.8& 1.8040870\\
		-11.7& 1.8041290\\
		-15.6& 1.8041712\\
		-19.4& 1.8042129\\
		-23.3& 1.8042540\\
		-27.2& 1.8042950\\
		-31.1& 1.8043357\end{matrix} $ \ \ \
	$\begin{matrix}         t& x(t)>x_0^\e\\
		\hfill 00.0& 1.810\hfill \\
		-03.8& 1.8099799\\
		-07.7& 1.8099597\\
		-11.6& 1.8099404\\
		-15.4& 1.8099205\\
		-19.3& 1.8099011\\
		-23.1& 1.8098819\\
		-26.0& 1.8098631\\
		-30.8& 1.8098447\end{matrix} $ \qquad
	$\begin{matrix}        t& x(t)<x_1^\e\\
		00.0& 1.132\hfill \\
		13.6& 1.1320219\\
		27.3& 1.1320438\\
		40.9& 1.1326554\\
		54.6& 1.1320873\\
		68.2& 1.1321088\\
		81.9& 1.1321304\\
		95.5& 1.1321518\\
		109.& 1.1321731\end{matrix} $ \ \ \ 
	$\begin{matrix}        t& x(t)>x_1^\e\\
		00.0& 1.138\hfill \\
		13.5& 1.1379961\\
		27.0& 1.1379918\\
		40.5& 1.1379874\\
		53.0& 1.1379833\\
		67.5& 1.1379793\\
		81.0& 1.1379753\\
		94.4& 1.1379714\\
		108.& 1.1379673\end{matrix} $ 
	
	\medskip                  
	$\begin{matrix} 
		\phantom{aaa} \text{2] Stable limit cycle} \\ \phantom{aaa} \text{ i.\,v. } (0,x_2^\e,1),\ \ x_2^\e\approx 1.277\ \end{matrix}
	\phantom{aaaaaaaaaaaaa} 
	\begin{matrix} 
		\text{2] Unstable limit cycle} \\  \text{ i.\,v. } (0,x_3^\e,1),\ \ x_3^\e\approx 1.301\ \end{matrix} $
	
	\smallskip
	$\begin{matrix}         t& x(t)<x_2^\e\\
		00.0& 1.272\hfill \\
		06.2& 1.2720017\\
		12.3& 1.2720037\\
		18.3& 1.2720054\\
		24.6& 1.2720075\\
		30.8& 1.2720095\\
		37.8& 1.2720114\\
		43.8& 1.2720133\\
		49.3& 1.2720149 \end{matrix} $ \ \ \
	$\begin{matrix}         t& x(t)>x_2^\e\\
		00.0& 1.282\hfill \\
		06.5& 1.2819995\\
		13.0& 1.2819988\\
		19.5& 1.2819980\\
		26.0& 1.2819979\\
		32.6& 1.2819974\\
		39.1& 1.2819967\\
		45.6& 1.2819960\\
		52.1& 1.2819955 \end{matrix} $ 
	\qquad
	$\begin{matrix}         t& x(t)<x_3^\e\\
		\hfill 00.0& 1.298\hfill \\
		-07.6& 1.2980020\\
		-15.1& 1.2980040\\
		-22.7& 1.2980059\\
		-30.2& 1.2980080\\
		-37.8& 1.2980101\\
		-45.4& 1.2980118\\
		-52.9& 1.2980140\\ 
		-60.5& 1.2980159 \end{matrix} $ \ \ \
	$\begin{matrix}         t& x(t)>x_3^\e\\
		\hfill 00.0& 1.304\hfill \\
		-08.7& 1.3039972\\
		-17.4& 1.3039944\\
		-26.0& 1.3039917\\
		-34.7& 1.3039887\\
		-43.4& 1.3039860\\
		-52.1& 1.3039833\\
		-60.7& 1.3039807\\ 
		-69.4& 1.3039774 \end{matrix} $ 
	
	\medskip                  
	$\begin{matrix} 
		\phantom{a}  \text{2] Stable limit cycle} \\ \phantom{a} \text{ i.\,v. } (0,x_{-2}^\e,1),\ \ x_{-2}^\e\approx -1.267\ \end{matrix}
	\phantom{aaaaaaaaa} 
	\begin{matrix} 
		\text{2] Unstable limit cycle} \\ \quad \text{ i.\,v. } (0,x_{-3}^\e,1),\ \ x_{-3}^\e\approx -1.301\ \end{matrix} $
	
	\smallskip
	$\begin{matrix}        t& x(t)>x_{-2}^\e\\
		00.0& -1.272\hfill \\
		06.2& -1.2719992\\
		12.3& -1.2719985\\
		18.5& -1.2719977\\
		24.7& -1.2719970\\
		30.8& -1.2719963\\
		37.0& -1.2719957\\
		43.2& -1.2719950\\
		49.3& -1.2719939\end{matrix} $ \ \ \
	$\begin{matrix}        t& x(t)<x_{-2}^\e\\
		00.0& -1.262\hfill \\
		05.9& -1.2620020\\
		11.8& -1.2620042\\
		17.7& -1.2620062\\
		23.6& -1.2620084\\
		29.5& -1.2620106\\
		35.4& -1.2620127\\
		41.3& -1.2620149\\
		47.2& -1.2620170 \end{matrix} $ 		
	\qquad
	$\begin{matrix}        t& x(t)>x_{-3}^\e\\
		\hfill 00.0& -1.304\hfill \\
		-09.0& -1.3039971\\
		-18.0& -1.3039950\\
		-26.9& -1.3039928\\
		-35.9& -1.3039904\\
		-44.8& -1.3039879\\
		-53.8& -1.3039857\\
		-62.8& -1.3039833\\
		-71.7& -1.3039808\end{matrix} $ \ \ \
	$\begin{matrix}        t& x(t)<x_{-3}^\e\\
		\hfill 00.0& -1.298\hfill \\										
		-07.6& -1.2980035\\
		-15.3& -1.2980069\\
		-22.9& -1.2980105\\
		-30.5& -1.2980141\\
		-38.2& -1.2980177\\
		-45.8& -1.2980212\\
		-53.4& -1.2980247\\
		-61.1& -1.2980283\end{matrix} $ }
	
	\smallskip
	Also, for each admissable solution $b_{kl}^j,l$ the constants $K_{kl}^{aj}$ from the nondegeneracy condition $(\ref{ud}^a)$ are calculated:
	$K_{00}^{a0}\approx 10.7,\,K_{k0}^{a1}\approx -4.3,$ $K_{kl}^{a2}\approx -0.28,\,K_{kl}^{a3}\approx 0.83$ \ $(k,l=\pm 1).$
	
	As a result, the following proposition has been proven.
	
	\smallskip
	{\bf Theorem ${\bf 8.1^a.}$} \ {\it For any $\e\in [0,0.05],$ the system \eqref{sa} has three limit cycles: 
		$LC_{\!\e 00}^{\, 0},$ which is stable, $LC_{\!\e k0}^{\, 1},$ which are unstable. For any $\e\in [0,10^{-3}],$ there are eight additional cycles: 
		the cycles $LC_{\!\e k,\pm 1}^{\, 3},$ which are stable, $LC_{\!\e k,\pm 1}^{\, 2},$ which are unstable \ $(k=\pm 1).$ }

\bigskip
\centerline{\large Conclusion}
 
\smallskip
{\quad} Let us mention possible areas of applications of the GTS method.
		
		1) The research of the periodic systems with more complicated Hamiltonians, for example, taken from \cite{Li}.
		
		2) The utilization of the GTS method to reinforce the existing results:
		
		$2_1)$ The most obvious instance of possible application is to assess the Hilbert number for system \eqref{loop} from \cite{Bas4} 
		due to it being a generalization of the systems from papers  \cite{Du} and \cite{ILY}. 
		Another possible application is the research of the systems with Hamiltonian from \cite{LLY} with periodic perturbations.
		
		$2_2)$ To consider already studied  higher dimensional autonomous systems with small parameter 
	 adding  the periodic or quasiperiodic perturbations to those systems (\cite{Bib2},\,\cite{Bas3}), 
		assuming basis frequencies satisfy the standard Diophantine-type condition. 
		For example, periodic perturbations can be considered in systems with unperturbed parts such as  
		$(x_2,-x_1^{2n-1},x_4,$ $-x_3^{2n+1})$ from \cite{Bas} or $(x_{1+d},-x_1^{2n-1},\ldots, x_{2d},-x_d^{2n-1})$ with $d\ge 2$ from \cite{Bas2}. 
		
		The special point of interest is the possibility to construct the invariant torus of the same dimension as a system. 
		To accomplish this, we consider the $2^d$-order system from \cite{Bas5} with periodic or quasiperiodic perturbation. 
		In mentioned paper the class of autonomous systems is constructed, including polynomial ones, 
		such that bifurcation scenario of a birth of invariant torus occurs with torus codimension equal to one.
		
		Notice, that the algorithm allowing to find invariant torus bifurcations of various dimensions that branch out 
		from the equilibrium point, was developed by Yu.\,N.\,Bibikov in \cite{Bib3}.
		
		3) The application of the GTS method for systems for which the critical periods are being found. For example, the system 
		$\dot x=y\prod\nolimits_{i=1}^k ((y-\alpha_i)^2+\e ),\ \dot y=-x\prod\nolimits_{i=1}^k ((x-\beta_i)^2 + \e)$ from \cite{XCen} can be studied. 

\bigskip
		{\bf Clarification.} The perturbations of the systems researched in papers \cite{Bas4},\,\cite{Bas3},\,\cite{BZ}, which were the start of the development and application of the GTS method used to find the invariant tori and limit cycles in autonomous case, require limitations similar to the ones made for the system \eqref{sv}.
		
\bigskip
		{\bf Acknowledgements.}\, I express the sincere thanks to my colleagues and friends S.\,A. Ivanov, Yu.\.A. Ilyin and S.\,G. Kryzhevich for the discussions of various questions, related to the complex and functional analysis. $\qquad$ V.\,V. Basov 


\newpage		

\bigskip\bigskip
{\small 
\begin{thebibliography}
\smallskip
\bibitem{Arn}
{\it Arnold V.\ I.}, Loss of stability of self-oscillations close to resonance and versal deformations of four-dimensional smooth manifolds, 
and the arithmetic of integral quadratic forms, Funct. Anal. Appl. {\bf 11} (1977), 85-92.

\bibitem{Bas4}  {\it Basov V.\,V.} 
"Invariant Surfaces of Standard Two-Dimensional Systems with Conservative First Approximation of the Third Order", 
{\it Differentsial'nye Uravneniya,} {\bf 44}(1), 3-18 (2008) \ [in Russian]. 
Eng. version: {\it Differential Equations,} 2008, 44:1, 1-18.

\bibitem{Lyap} {\it Lyapunov, A.M.,}  “Study of a special case of the motion stability problem,” in: 
Collection of Works, Vol.\,2, Moscow-Leningrad, Akad.\,Nauk SSSR (1956), pp. 272–331.

\bibitem{Hale}  {\it Hale J.\,K.} 
"Integral Manifolds of Perturbed Differential Systems",
{\it Annals of Mathematics Second Series,} {\bf 73}(3), 496-531 (1961).

\bibitem{LiQ} {\it Li J., Huang Q.}
Bifurcations of limit cycles forming compound eyes in the cubic system, Chin. Ann. Math, {\bf B8} (1987), 391-403.

\bibitem{Bas3}  {\it Basov V.\,V.}  
"Invariant surfaces of two-dimensional periodic systems with bifurcating rest points in the first approximation",  
{\it Journal of Mathematical Sciencies,} {\bf 147}(1), 6398-6415 (2007). 
Translated from Contemporary Mathematics and Its Applications, Vol. 38, Suzdal Conference-2004, Part\,3, 2006.
https://doi.org/10.1007/s10958-007-0474-x

\bibitem{Bib2} {\it Bibikov Yu.\,N.} "Bifurcation of the generation of invariant tori with infinitesimal frequency"\,,
{\it Algebra i Analiz,} {\bf 10}(2), 81-92 (1998) \ [in Russian].
Eng. version: {\it  St. Petersburg Mathematical Journal}, 1999, 10:2, 283-292.	

\bibitem{Var}
{\it Varchenko A.\ N.}, An estimate of number of zeros of an Abelian integral depending on a parameter and limiting cycles, Funct. Anal. Appl. {\bf 18}(1984), 98-108.

\bibitem{Ilyas}
{\it Ilyashenko Yu.\ S.}, Finiteness theorems for limit cycles, American Mathematical Society, Providence, RI, 1991.
Russian Text The Author, 1990, published in Uspekhi Mat. Nauk, {\bf 54}(2), pp. 143-200. 
	
\bibitem{Mel} 
{\it Mel'nikov, V.\ K.} On the stability of a center for time-periodic perturbations. 
Trudy Moskov. Mat. Ob\v s\v c. {bf 12} (1963) 3-52. [in Russian]. 

\bibitem{GH} {\it Guckenheimer J.,\ Holmes P.}
Nonlinear oscillations, dynamical systems, and bifurcations of vector fields, Springer-Verlag (1983).

\bibitem{Li} {\it Li J.}
Hilbert's 16th problem and bifurcations of planar polynomial vector fields,
{\it International Journal of Bifurcation and Chaos}, Vol. 13, No. 1 (2003), 47-106.

\bibitem{ILY} {\it Iliev I.\ D.\ Li C.\ Yu J.} On the cubic perturbations of the symmetric 8-loop Hamiltonian (2019),
arXiv:1909.09840v1. 

\bibitem{Ili1}  {\it Iliev I.\,D.,\ Perko L.} Higher order bifurcations of limit cycles,
{\it Diff Equat} {\bf 154} (1999), 339-363.

\bibitem{Du} {\it Dumortier F.,\ Li C.} Perturbation from an elliptic Hamiltonian of degree four --- IV figure eight-loop,
{\it Differential Equations} {\bf 188} (2003), 512-554.

\bibitem{WTX} {\it Wei L.,\ Tian Y.\ Xu Y.} 
The Number of Limit Cycles Bifurcating from an Elementary Centre of Hamiltonian Differential Systems.
Mathematics 10, (2022), 1483-1496. https: doi.org/10.3390/math10091483

\bibitem{China} {\it Wei M.,\ Cai J.,\ Zhu H.}
Poincare Bifurcation of Limit Cycles from a Lienard System with a Homoclinic Loop Passing through a Nilpotent Saddle,
{\it Discrete Dynamics in Nature and Society}, Hindawi, vol. 2019

\bibitem{CL} {\it Christopher C.\,J., Lloyd N.\,G.}
Polynomial systems: A lower bound for the Hilbert numbers, {\it Proc. Royal Soc. London Ser.}, {\bf A450} (1995), 219-224. 

\bibitem{LLY} {\it Li C.,\ Liu C.,\ Yang J.}
A cubic system with thirteen limit cycles,
{\it Journal of Differential Equations} {\bf 246}(2009), 3609-3619.

\bibitem{Han} {\it Han M.-A.} Bifurcations of Invariant Tori and Subharmonic Solutions for Periodic Perturbed Systems. 
{\it Science in China}, {\bf 37}(11), (1994). 

\bibitem{BZ} {\it Basov V.\,V.,\ Zhukov A.\,S.}
"Invariant Surfaces of Periodic Systems with Conservative Cubic First Approximation", 
{\it Vestnik Sankt-Peterburgskogo Universiteta: Matematika, Mekhanika, Astronomiya,} {\bf 64}(3), 376-393 (2019) \ [in Russian].
Eng. version: {\it Vestnik St.\,Petersburg University, Mathematics,} 2019, 52:3, 244-258.

\bibitem{Bib3} Bibikov, Yu.N., Lokal’nye problemy teorii mnogochastotnykh nelineinykh kolebanii (Local Problems of the Theory of Multifrequency Nonlinear Oscillations), St. Petersburg, 2003.

\bibitem{Hale1}  {\it Chow S.-N., Hale J.\,K.} 
Methods of bifurcation theory, {\it N.Y., Springer-Verlag,} (1982), 515\,p.

\bibitem{Moz}  {\it Mozer J.} A rapidly convergent iteration method and non-linear differetial equations, 
{\it Ann.\,Scuola Norm.\,Sup.\,Pisa Ser. III,} {\bf 20}, (1966), 265-315.

	\bibitem{Bas}  {\it Basov V.\,V.}
	"Bifurcation of the Equilibrium Point in the Critical Case of Two Pairs of Zero Characteristic Roots",
	Differential equations and dynamical systems, Collected papers.
	Dedicated to the 80th anniversary of academician Evgenii Frolovich Mishchenko,
	Tr. Mat. Inst. Steklova, 236, Nauka, Moscow, 2002, 45-60 \ [in Russian].  (http://mi.mathnet.ru/eng/tm275).
	Eng. version: Proceedings of the Steklov Institute of Mathematics, 2002, 236, 37-52.	
	
	\bibitem{Bas2}  {\it Basov V.\,V.} 
	"Bifurcation of the Point of Equilibrium in Systems with Zero Roots of the Characteristic Equation",  
	{\it Mat. Zametki,} {\bf 75}(3), 323-341 (2004) \ [in Russian]. 
	Eng. version: {\it Mathematical Notes}, 2004, 75:3, 297-314.
	
	\bibitem{Bas5}  {\it Basov, V.V.} Bifurcation of Invariant Tori of Codimension One. Mathematical Notes 69:1, 3–16 (2001). 
	https://doi.org/10.1023/A:1002888709793
	
	\bibitem{XCen} {\it Xiuli C.} New lower bound for the number of critical periods for planar polynomial systems (2019), 
	doi: 10.13140/rg.2.2.34326.98885
	
	\bibitem{Zor}  {\it Zorich V.\,A.} Matematicheskiy analiz, Vol. 2, 9th edition, MCCME, Moscow, 2019
	
\end{thebibliography} }
\end{document}